\documentstyle[sprocl,twoside,amssymb,amsmath,eucal]{article}

\newcommand {\supplus}{\mathop{{\supset}\llap{\raise 
0.5pt\hbox{\normalfont\small+}\hskip 0.5pt}}} 
\newcommand {\subplus}{\mathop{{\subset}\llap{\raise 
0.5pt\hbox{\normalfont\small+}\hskip 0.5pt}}}  
\input cyracc.def
\newfam\cyrfam
\font\tencyr=wncyr10
\font\sevencyr=wncyr7

\def\cyr{\fam\cyrfam\tencyr\cyracc}
\textfont\cyrfam=\tencyr \scriptfont\cyrfam=\sevencyr
  \scriptscriptfont\cyrfam=\sevencyr

\font\tencyr=wncyr10

\newcommand {\Cee}    {{\mathbb  C}}
\newcommand {\Kee}    {{\mathbb  K}}
\newcommand {\Ree}    {{\mathbb  R}}
\newcommand {\Zee}    {{\mathbb  Z}}

\newcommand {\fa}     {{\mathfrak{a}}}

\newcommand {\fab}    {{\mathfrak{ab}}} 

\newcommand {\fas}    {{\mathfrak{as}}}
 
\newcommand {\fb}     {{\mathfrak{b}}}
\newcommand {\fc}    {{\mathfrak{c}}}

\newcommand {\fcvect}   {{\mathfrak{cvect}}}
\newcommand {\fder}   {{\mathfrak{der}}}   %

\newcommand {\fg}     {{\mathfrak{g}}}    %
\newcommand {\fgl}    {{\mathfrak{gl}}}  %
\newcommand {\fh}     {{\mathfrak{h}}}
\newcommand {\fhei}   {{\mathfrak{hei}}}
\newcommand {\fk}     {{\mathfrak{k}}}

\newcommand {\fL}     {{\mathfrak{L}}}
\newcommand {\fle}    {{\mathfrak{le}}}
\newcommand {\fm}     {{\mathfrak{m}}}

\newcommand {\fo}     {{\mathfrak{o}}}
\newcommand {\fosp}   {{\mathfrak{osp}}}

\newcommand {\fs}     {{\mathfrak{s}}}

\newcommand {\fsl}    {{\mathfrak{sl}}}
\newcommand {\fsle}   {{\mathfrak{sle}}}

\newcommand {\fsp}    {{\mathfrak{sp}}}
\newcommand {\fspe}   {{\mathfrak{spe}}}

\newcommand {\fsvect} {{\mathfrak{svect}}}
\newcommand {\fsu}    {{\mathfrak{su}}}

\newcommand {\fu}     {{\mathfrak{u}}}
\newcommand {\fv}     {{\mathfrak{v}}}     %
\newcommand {\fvect}  {{\mathfrak{vect}}}   %

\newcommand {\cD}     {{\cal D}}

\newcommand {\cF}     {{\cal F}}

\newcommand {\cL}     {{\cal L}}
\newcommand {\cM}     {{\cal M}}

\newcommand {\cO}     {{\cal O}}

\def \opname#1#2%
  {\expandafter\newcommand \csname #1\endcsname {{\mathop{#2}\nolimits}}}

\newcommand{\rmname}[1]
  {\expandafter\newcommand \csname #1\endcsname {{\operatorname{#1}}}}

\newcommand{\rmnameii}[2]
  {\expandafter\newcommand \csname #1\endcsname {{\operatorname{#2}}}}

\rmname{act}
\rmname{Ad}
\rmname{Add}
\rmname{ad}
\rmname{Alt}
\rmname{alt}
\rmname{Ann}
\rmname{antidiag}
\rmname{Ber}
\rmname{ber}
\rmname{Br}
\rmname{card}
\rmname{ch}
\rmname{Char}
\rmname{cem}
\rmname{cj}
\rmname{Cliff}
\rmname{cntr}
\rmname{codim}
\rmname{coind}
\rmname{const}
\rmname{col}
\rmname{cork}
\rmname{cpr}
\rmname{diag}
\rmnameii{Div}{div}
\rmname{Def}
\rmname{Der}
\rmname{Dim}
\rmname{End}
\rmname{Even}
\rmname{Ext}
\rmname{gr}
\rmname{Hom}
\rmname{HT}
\rmnameii{Ht}{ht}
\rmname{hwt}
\rmname{Id}
\rmname{id}
\rmname{ind}
\rmname{Ind}
\rmname{Inf}
\rmname{irr}
\rmname{Le}
\rmname{Lie}
\rmname{lwt}
\rmname{mult}
\rmname{Mor}
\rmname{nm}
\rmname{Ob}
\rmname{Odd}
\rmname{Osc}
\rmname{per}
\rmname{Pic}
\rmname{pr}
\rmname{pro}
\rmname{Prime}
\rmname{Proj}
\rmname{prt}
\rmname{pt}
\rmname{Q}
\rmname{qet}
\rmname{qtr}
\rmname{rd}
\rmname{rk}
\rmname{row}
\rmname{Res}
\rmname{salt}
\rmname{Sch}
\rmname{SBr}
\rmname{scalar}
\rmname{Ser}
\rmname{sign}
\rmname{Smbl}
\rmname{spin}
\rmname{ssym}
\rmname{str}
\rmname{st}
\rmname{sgn}
\rmname{sq}
\rmname{symm}
\rmname{supp}
\rmname{Supp}
\rmname{St}
\rmname{Spec}
\rmname{Spm}
\rmname{tr}
\rmname{weyl}
\rmname{Weyl}
\rmname{Witt}

\opname{vvol}  {{v\hspace{-0.1ex}o\hspace{-0.02ex}l\/}}
\opname{pnt}  {\text{\normalfont pt}}
\opname{Span} {{Span}}
\opname{slim} {\overline{\lim}}
\opname{Vol}  {{V\hspace{-0.55ex}o\hspace{-0.02ex}l\/}}
\opname{Par}  {{P\hspace{-0.3ex}a\hspace{-0.05ex}r\/}}

\rmname{Mat}
\rmname{Bil}
\rmname{Diff}
\rmname{Ker}
\rmname{Herm}
\rmname{Coker}
\rmname{Conn}
\rmname{Covect}
\rmname{Vect}
\rmname{Int}

\rmnameii {IM} {Im}
\rmnameii {RE} {Re}

\opname{Aut} {{A\hspace{-0.2ex}u\hspace{-0.1ex}t\/}}
\opname{GL} {{G\hspace{-0.3ex}L}}
\opname{SL} {{S\hspace{-0.3ex}L}}
\opname{Exp} {{E\hspace{-0.2ex}x\hspace{-0.1ex}p\/}}
\opname{GQ} {{G\hspace{-0.2ex}Q}}
\opname{OSp} {{O\hspace{-0.25ex}S\hspace{-0.15ex}p\/}}
\opname{Out} {{O\hspace{-0.25ex}u\hspace{-0.15ex}t\/}}
\opname{Spp} {{S\hspace{-0.2ex}p\/}}
\opname{SpO} {{S\hspace{-0.2ex}p\hspace{-0.02ex}O\/}}
\opname{Pe} {{P\hspace{-0.25ex}e\/}}
\opname{SPe} {{S\hspace{-0.25ex}P\hspace{-0.25ex}e\/}}
\opname{Spin} {{S\hspace{-0.25ex}p\hspace{-0.05ex}i\hspace{-0.1ex}n\/}}
\opname{Iso} {{I\hspace{-0.25ex}s\hspace{-0.1ex}o\/}}
\opname{SSPe} {{S\hspace{-0.25ex}S\hspace{-0.15ex}P\hspace{-0.25ex}e\/}}
\opname{PeU} {{P\hspace{-0.25ex}e\hspace{-0.1ex}U\/}}
\opname{QU} {{Q\hspace{-0.15ex}U\/}}
\opname{U} {{U\/}}

\opname{cGQ} {{\cal G \hspace{-0.2em} Q \/}}
\opname{cSL} {{\cal S \hspace{-0.2em} L \/}}
\opname{cGL} {{\cal G \hspace{-0.2em} L \/}}
\opname{cOSp} {{\cal O \hspace{-0.2em} S \hspace{-0.3em} \it p\/}}
\opname{cPe} {{\cal P \hspace{-1.5pt} \it e\/}}
\opname{cVect} {{\cal V \hspace{-1.5pt} \it e\hspace{-0.1ex}c\hspace{-0.1ex}t\/}}
\opname{cVol} {{\cal V \hspace{-1.5pt} \it o\hspace{-0.1ex}l\/}}
\opname{cAut} {{\cal A \hspace{-0.2em} \it u\hspace{-0.1em}t\/}}
\opname{cCovect} {{\cal C \hspace{-1.5pt}
     \it o\hspace{-0.1ex}v\hspace{-0.1ex}e\hspace{-0.1ex}c\hspace{-0.1ex}t\/}}
\opname{CW} {{C\hspace{-0.15ex}W}}

\newcommand {\ev} {{\bar0}}
\newcommand {\od} {{\bar1}}

\newcommand {\tto} {\longrightarrow}
\newcommand {\pder}[1] {{\frac{\partial}{\partial {#1}}}}
\newcommand {\pderf}[2] {{\frac{\partial {#1}}{\partial {#2}}}}

\newcommand {\bcdot}   {\mathbin{\hbox{\raise.4ex\hbox{\bf.}}}} % bold \cdot

\def\be{\begin{equation}}
\def\ee{\end{equation}}
\def\bea{\begin{eqnarray}}
\def\eea{\end{eqnarray}}
\def\citebk#1{\hspace{0.9mm}\raisebox{-1.85mm}[0mm][0mm]
  {\Large\cite{#1}}\hspace{-0.1mm}}
\begin{document}
\sloppy

%\title{Invariant operators on the exceptional supermanifolds and the
%Standard Models}
\title{INVARIANT OPERATORS ON SUPERMANIFOLDS AND STANDARD MODELS}

\author{PAVEL GROZMAN, DIMITRY LEITES}

\address{Department of Mathematics, University of Stockholm,
Kr\"aftriket hus 6,\\ SE-106 91 Stockholm, Sweden;
mleites@matematik.su.se} 

\author{IRINA SHCHEPOCHKINA}

\address{The Independent University of Moscow,
Bolshoj Vlasievsky per, dom 11,\\ Moscow RU-121 002, Russia;
Ira@Pa\-ra\-mo\-no\-va.mccme.ru}


\maketitle

%\vspace{1cm}

\abstracts{Here we continue to list the differential operators
invariant with respect to the 15 exceptional simple Lie superalgebras
$\fg$ of polynomial vector fields.  A part of the list (for operators
acting on tensors with finite dimensional fibers) was earlier obtained
in 2 of the 15 cases by Kochetkov and in one more instance by Kac and
Rudakov.  Broadhurst and Kac conjectured that some of these structures
pertain to the Standard Models of elementary particles and the Grand
Unified Theories.  So, GUT, if any exists, will be formulated in terms
of operators we found, or their $r$-nary analogs to be found. 
Calculations are performed with the aid of Grozman's {\it
Mathematica}-based SuperLie package.  When degeneracy conditions are
violated (absence of singular vectors) the corresponding module of
tensor fields is irreducible.  We also verified some of the earlier
findings.}

\vspace*{0.3cm}

\tableofcontents

\vspace*{0.3cm}

%\newpage
\subsubsection{To the memory of Misha Marinov (by D.L.)}  
Only one of us knew Misha.  %I deleted ``At the time'' added by Ed 
%because this specification is wrong
It was at the flat of one of my teachers,
Felix Berezin, that we became acquainted.  I was a student then and
was bubbling with enthusiasm on account of fantastic visions,
applications of supersymmetry, which at that time was not so clear to
most.  Misha shared the awe %why AWE? neither of us was afraid on 
%account of SUSY ...I originally
%wrote ``sensation'', but ``feeling'' or some such is also OK
of the beauty of the mathematics related with SUSY. Besides this,
he was a charming person.  Moreover, he was writing a review on
superstrings.  So Misha fascinated me.  

His review\,\cite{M} illuminated my research for many years ahead
(Ref.\,\protect\citebk{GLS} is just a part of the story).  Later
Misha, being a man of passion, bitterly criticized the whole trend as
useless waste of many talents of several generations.  I did not take
his vitriolic words at face value.  I knew that however strong
language he used, it was merely rhetoric.  He similarly criticized me
for not being able to quantize the ``odd mechanic'' (and
Ref.\,\protect\citebk{LSH2} is a part of our answer to Misha) nor
demonstrate the ``usefulness'' of what was rediscovered several years
later by Batalin and Vilkovisky under the name antibracket together
with the manifest demonstrations of its usefulness (cf.\
Refs.\,\protect\citebk{L0} and \protect\citebk{GPS}).  When the
antibracket was rediscovered Misha came to me and said, ``Do you
remember that stuff about the two types of mechanics, you told me some
time ago?  Recently Batalin talked a lot about what looks precisely
like it, and he is one of the rare few who only speaks about what he
knows.  If I had not known him for a long time I'd think he had become
crazy.  He plans to do almost all gauge field theory by means of this
new mechanic.  So, let me get you acquainted, you two talk similar
\ldots''

Later, before he applied for emigration and while
waiting for permission (more likely for a refusal, at that time)  
Misha frequented an informal seminar I hosted
at home being also ready to quit the regular ways.  Some
mathematicians having applied for emigration became {\it personae non
grata} at official seminars; in my flat they could discuss new and old
results with those who dared.  Misha was our physics teacher, he tried
to convey to us (J.~Bernstein, A.~Beilinson, me, and occasional
``guests'', from students to S.~Gindikin) some ideas of modern
physics.

Misha visited me at Stockholm recently.  I took him to look at the frescos
of T\"aby kyrkan, one of the most interesting churches in Stockholm. 
The most famous of the frescos is currently called ``The Knight playing
chess with Death''.  The original title, if it ever existed, is
forgotten; the modern title is hardly correct: depicted are a
chess-board ($5\times 5$, to simplify Death's task even further, I
presume), Death, and the Knight, smiling his best (as his PR advisor,
no doubt, taught him) politely inclined towards the ruling party.  The
Knight stands on the same side of the board as Death, so manifestly
both are playing together against {\it us}.  Misha, who was restless during
the day being anxious not to miss his train, suddenly relaxed having
seen this fresco.  This was the last fresco we saw and on the way back
we discussed something neutral.  He said he enjoyed the frescoes.  On
our way back we got into a traffic jam and missed his train by a few
minutes.  Naturally, I paid the penalty.  Misha wrote to me later
insisting to repay or split the cost, so I answered that I did not change that
much after emigration, and he should know better.  He retorted as in
the good old days when we all were younger saying, ``I regret that Israel is so
small, and so penalties for missing the train are low.  Anyway, I hope to
see you one day and promise to make you miss the train from Haifa
to wherever.''  He also gave me a reprint of his superstring paper
(somebody ``borrowed'' my copy of the reprint he gave me in Moscow)
and wrote in Russian: ``To DL, as a memory of our friendship that I hope not
 interrupt''.  The inscription looked very strange to me.  Now I
understand it: he already knew about his cancer.

I think all who knew Misha love him.  One can not use past tense with
this verb.  And therefore it is infinitely sad not to be able, ever,
to miss the train in Misha's company, this side of the Chess-Board we
all try to describe before we meet the ultimate opponent.

\section{Introduction}

This is a part of an expanded transcript of two talks at the
International Workshops ``Supersymmetry and Quantum Symmetry'': (1) at
Dubna, July 22-26, 1997, we described the case of $\fk\fa\fs$ and
$\fk(1|6)$, sketched the case $\fk(1|n)$ and mentioned earlier results
of Kochetkov on $\fv\fle(4|3)$ and $\fm\fb(4|5)$, (2) at Karpacz,
September 23, 2001, we considered $\fm\fb(3|8)$ and
$\fk\fs\fle(5|10)$.

Broadhurst and Kac observed\,\cite{Ka1} that some of the exceptional
Lie superalgebras (listed in Refs.\,\protect\citebk{S4}, \protect\citebk{CK}) 
might pertain to a GUT or 
the Standard Model, their linear parts being isomorphic to $\fsl(5)$
or $\fsl(3)\oplus\fsl(2)\oplus\fgl(1)$.  Kac demonstrated\,\cite{Ka2}
that for the Standard Model with $\fsu(3)\oplus\fsu(2)\oplus\fu(1)$ as
the gauge group a certain remarkable relation between $\fv\fle(3|6)$
and some of the known elementary particles does take place; it seems
that for $\fm\fb(3|8)$ there is even better correspondence.

The total lack of enthusiasm from the physicists' community concerning
these correspondences is occasioned, perhaps, by the fact that no real
form of any of the simple Lie superalgebras of vector fields with
polynomial coefficients has a unitary Lie algebra as its linear part. 
Undeterred by this, Kac and Rudakov calculated\,\cite{KR1} some
$\fv\fle(3|6)$-invariant differential operators.  They calculated the
operators for finite dimensional fibers only.  This restriction makes
calculations a sight easier but strikes out many operators.  The
amount of calculations for $\fm\fb(3|8)$ is too high to be performed
by hands.

The problem we address --- calculation of invariant differential
operators acting in tensor fields on manifolds and supermanifolds with
various structures --- was a part of our Seminar on Supermanifold's
agenda since mid-70's.  Here we use Grozman's code
SuperLie\,\cite{GL1} to verify and correct earlier results and obtain
new ones, especially when bare hands are inadequate.  The usefulness
of SuperLie was already demonstrated when we calculated the left-hand
side of $N$-extended SUGRA equations for any $N$, cf.\
Ref.\,\protect\citebk{GL2}.  We review the whole field with its open
problems and recall interesting Kirillov's results and problems buried
in the VINITI collection\,\cite{Ki}
%not very accessible VINITI collection.\cite{Ki} 
which is not very accessible.  Another nice (and accessible) review we
can recommend in addition to Ref.\,\protect\citebk{Ki} is
Ref.\,\protect\citebk{KMS}.

{\bf What is done} \quad  

(1) We list (degeneracy conditions) all differential operators, or
rather the corresponding to them singular vectors, of degrees 1, and 2
and, in some cases, of all possible degrees (which often are $\leq
2$), invariant with respect to several exceptional Lie superalgebras
of vector fields.  When degeneracy conditions are violated (absence of
singular vectors) the corresponding induced and coinduced modules are
irreducible.  For some exceptions EVERY module $I(V)$ has a singular
vector.  This is a totally new feature.  For $\fv\fle(3|6)$ (and
finite dimensional fibers) the answer coincides with Kac-Rudakov's
one.

(2) We observe that the linear parts of two of the W-regradings of
$\fm\fb$ are Lie superalgebras strictly greater than
$\fsl(3)\oplus\fsl(2)\oplus\fgl(1)$.  They (or certain real forms of
them) are natural candidates for the algebras of The {\it would be}
Standard Models; modern ``no-go'' theorems do not preclude them.

\subsection{Veblen's problem} The topology of differentiable manifolds
has always been related with various geometric objects on them and, in
particular, with operators invariant with respect to the group of
diffeomorphisms of the manifold, operators which act in the spaces of
sections of ``natural bundles",\cite{KMS} whose sections are tensor
fields, or connections, etc.  For example, an important invariant of
the manifold, its {\it cohomology}, stems from the de Rham complex
whose neighboring terms are bridged by an invariant differential
operator --- the exterior differential.

The role of invariance had been appreciated already in XIX century in
relation with physics; indeed, differential operators invariant with
respect to the group of diffeomorphisms preserving a geometric
structure are essential both in formulation of Maxwell's laws of
electricity and magnetism and in Einstein--Hilbert's formulation of
relativity.

Simultaneously, invariance became a topic of conscious interest for
mathematicians: the representation theory flourished in works of
F.~Klein, followed by Lie and \' E.~Cartan, to name the most important
contributors; it provided with the language and technique adequate in
the study of geometric structures.  Still, it was not until
O.~Veblen's talk in 1928 at the Mathematical Congress in
Bologna\,\cite{V} that invariant operators (such as, say, Lie
derivative, the exterior differential, or integral) became the primary
object of the study.  In what follows we rule out the integral and
other non-local operators; except in Kirillov's example, we only
consider local operators.

Schouten and Nijenhuis tackled Veblen's problem: they reformulated it
in terms of modern differential geometry and found several new
bilinear invariant differential operators.  Schouten conjectured that
there is essentially one {\it un}ary invariant differential operator:
the exterior differential of differential forms.  This conjecture had
been proved in particular cases by a number of people, and in full
generality in 1977--78 by A.~A.~Kirillov and, independently,
C.~L.~Terng (see Refs.\,\protect\citebk{Ki}, \protect\citebk{T}).

Thanks to the usual clarity and an enthusiastic way of Kirillov's
presentation he drew new attention to this problem, at least, in
Russia.  Under the light of this attention it became clear (to
J.~Bernstein) that in 1973 A.~Rudakov\,\cite{R1} also proved this
conjecture by a {\it simple algebraic} method which reduces Veblen's
problem for differential operators to a \lq\lq computerizable" one. 

Thus, a tough analytic problem reduces to a problem formally
understandable by any first year undergraduate: a series of systems of
linear equations in small dimensions plus induction.  The only snag is
the volume of calculations: to list all unary operators in the key
cases requires a half page; for binary operators one needs about 50
pages and induction becomes rather nontrivial; for $r$-nary operators
with $r>2$ only some cases seem to be feasible.

Later Rudakov and, for the contact series, I.~Kostrikin, classified
unary differential operators in tensor fields on manifolds invariant
with respect to each of the remaining three simple infinite
dimensional $\Zee$-graded Lie algebras $\fL$ of vector fields.  In
passing, the definition of the tensor field was generalized and
primitive forms came to foreground.

\subsection{Rudakov's breakthrough (following Bernstein$^{5}$)}
Hereafter the ground field $\Kee$ is $\Cee$ or $\Ree$.  Without going
into details which will be given later, observe that the spaces in
which invariant operators act fall into two major cases: spaces of
tensor fields (transformations depend on the 1-jet of diffeomorphism)
and spaces depending on higher jets, called {\it HJ-tensors} for
short.  We will only study tensors here, not HJ-tensors.

1) Instead of considering $\Diff (U)$-invariant operators, where $U$
is a local chart, let us consider $\fvect (U)$-invariant operators,
where $\fvect (U)$ is the Lie algebra of vector fields on $U$ with
polynomial coefficients, or its formal completion.  ({\it A posteriori}
one proves that the global and the local problems are equivalent, 
cf.\  Ref.\,\protect\citebk{BL3}).  Accordingly, instead of tensor fields with
smooth coefficients, we consider their formal version: $T(V)=V\otimes
\Kee[[x]]$, where $x=(x_{1}, \dots, x_{n})$ and $n=\dim U$.

2) We assume here that $V$ is an {\it irreducible $\fgl(n)$-module
with lowest weight}.  Observe that while the requirement of lowest
weight seems to be ``obviously'' reasonable, that on irreducibility is
not, unless we confine ourselves to finite dimensional modules $V$. 
In super setting we are forced, in the absence of complete
reducibility, to consider {\it indecomposable} representations even
for finite dimensional modules.  Irreducible modules is just the
simplest first step.

3) Instead of the coinduced module, $T(V)$, consider the dual induced
module, $I(V^*)=\Kee[\partial]\otimes V^*$, where
$\partial_{i}=\pder{x_{i}}$.  The reason: formulas for $\fvect
(U)$-action are simpler for $I(V^*)$ than for $T(V)$. (The results,
contrariwise, are more graphic in terms of tensor fields.)

Observe that each induced module is a ``highest weight one'' with
respect to the whole $\fg=\fvect (U)$, i.e., the vector of the most
highest weight with respect to the linear vector fields from
$\fg_{0}=\fgl(n)$ is annihilated by $\fg_{+}$, the subalgebra of $\fg$
consisting of all operators of degree $>0$ relative the standard
grading ($\deg x_{i}=1$ for all $i$).  

In what follows the vectors annihilated by $\fg_{+}$ will be called
{\it singular} ones.

4) To every $r$-nary operator $D: T(V_{1})\otimes \dots \otimes
T(V_{r})\tto T(V)$ the dual operator corresponds
\vspace*{-2mm}
$$
D^*: I(W)\tto I(W_{1})\otimes \dots \otimes
I(W_{r}), \quad {\rm  where }\quad W=V^{*},\;  W_{i}=W_{i}^*,
$$
and, since (for details see Ref.\,\protect\citebk{R1}) each induced
module is a highest weight one, {\sl to list all the $D$'s, is the
same as to list all the $\fg_{0}$-highest singular vectors $D^*\in
I(W_{1})\otimes \dots \otimes I(W_{r})$.} In what follows $r$ is
called the {\it arity} of $D$.

5) In super setting, as well as for non-super but infinite dimensional
one, the above statement is not true: the submodule generated by a
singular vector does not have to be a maximal one; it could have
another singular vector of the same degree due to the lack of complete
reducibility.

For unary operators and Lie algebras this nuisance does not happen;
this was one of the (unreasonable) psychological motivations to stick
to the finite dimensional case even for Lie superalgebras,\cite{BL4}
cf.\  Ref.\,\protect\citebk{Sm4}.

6) Rudakov's paper\,\cite{R1} contains two results: 

\hskip 0.5 cm (A) description of $\fvect(n)$-invariant operators in
tensor fields (only the exterior differential exists) and (the main 
bulk of the paper)

\hskip 0.5 cm (B) proof of the fact that between the spaces of
HJ-tensors there are no unary invariant operators.

{\bf Problems} \quad 1) Describe $r$-nary invariant operators in the
spaces  of HJ-tensors for $r>1$. 

2) Describe  $r$-nary invariant operators in the {\it super}spaces  of
HJ-tensors. 

{\bf The dual operators} \quad Kirillov noticed\,\cite{Ki} that by
means of the invariant pairing (we consider fields on $M$ with compact
support and tensoring over the space $\cF$ of functions)
\vspace*{-2mm}
$$
B: T(V)\times (T(V^*)\otimes_{\cF}
\Vol(M))\tto \Ree, \quad (t, t^*\otimes\vvol)\mapsto \int (t, t^*)\vvol
$$

\vspace{-2mm}

\noindent
one can define the duals of the known invariant operators.  For the
fields with formal coefficients we consider there is, of course, no
pairing, but we consider a {\it would be} pairing induced by smooth
fields with compact support.  So the {\it formal dual} of $T(V)$ is
not $T(V^*)$ because the pairing returns a function instead of a
volume form to be integrated to get a number, and not $T(V)^*$ because
$T(V)^*$ is a highest weight module while we need a lowest weight one. 
Answer: the formal dual of $T(V)$ is $T(V^*)\otimes_{\cF}\Vol(M)$.

Possibility to dualize, steeply diminishes the number of cases to
consider in computations and helps to check the results.  Indeed, with
every invariant operator $D: T(V)\tto T(W)$ the dual operator $D^*:
T(W^*)\otimes_{\cF}\Vol(M)\tto T(V^*)\otimes_{\cF}\Vol(M)$ is also
invariant.  For example, what is the dual of $d:\Omega^k\tto
\Omega^{k+1}$?  Clearly, it is the same $d$ but in another
incarnation: $d: \Omega^{n-k-1}\tto \Omega^{n-k}$.  Though, roughly
speaking, we only have one operator, $d$, the form of singular vectors
corresponding to $d$ differs with $k$ and having found several ``new''
singular vectors we must verify that the corresponding operators are
indeed distinct.  This might be not easy.

Observe that these arguments do not work when we allow infinite
dimensional fibers (dualization sends the highest module into a lowest
weight one, so it is unclear if a highest weight module with a
singular vector always correspond to this lowest weight one). 
Sometimes, being tired of calculations, or when the computer gave up,
we formulated the description of singular vectors ``up to
dualization''; sometimes even the computer became ``tired''.  We will
mention such cases extra carefully; we intend to reconsider these
cases on a more powerful computer.  \vspace{-2mm}

\subsection{Further ramifications of Veblen's problem}
Rudakov's arguments\,\cite{R1} show that the fibers of HJ-tensors have
to be of infinite dimension; the same holds for Lie superalgebras,
though arguments are different.  Traditionally, fibers of tensor
bundles were only considered to be of finite dimension, though even in 
his first paper on the subject Rudakov\,\cite{R1} digressed from 
traditions.

1$^o$.  In the study of invariant operators, one of the ``reasons''
for confining to tensors, moreover, the ones corresponding to finite
dimensional fibers, is provided by two of Rudakov's results:\,\cite{R1}
(1) there are no invariant operators between HJ-tensors, (2) starting
with any highest weight modules $I(V)$, Rudakov unearthed singular
vectors only for fundamental (hence, finite dimensional) representations. 
Though (1) only applies to {\it unary} operators, researchers were
somewhat discouraged to consider HJ-tensors even speaking about binary
operators.

Consider invariant operators of arity $>1$ between the spaces of
HJ-tensors (see Problem above).  Is it true that in this case there
are no invariant operators either?

2$^o$.  Kirillov proved that (having fixed the dimension of the 
manifold and arity) {\sl the degree of invariant (with respect to $\fvect(n)$)
differential operators is bounded, even $\dim$ of the space of
invariant operators is bounded.}

There seems to be no doubt that a similar statement holds on
supermanifolds ...  but Kochetkov's examples reproduced below and our
own ones show that these expectations are false in some cases.

{\bf Problem} \quad Figure out the conditions when the dimension of the
space of invariant operators is bounded.  (We conjecture that this is
true for all the series of simple vectorial Lie superalgebras in the
standard grading.)

3$^o$.  On the line, all tensors are $\lambda$-densities and every 
$r$-linear differential operator is of the form
$$
L: (f_{1}dx^{\lambda_{1}}, \dots,  f_{r}dx^{\lambda_{r}})\tto 
P_{L}(f_{1}, \dots,  f_{r})dx^{\lambda}.
$$
Kirillov shows (with ease and elegance) that invariance of $L$ is
equivalent to the system
$$
\sum_{s=1}^r\left[t_{s}\frac{\partial ^{j+1}}{\partial 
t_{s}^{j+1}}+(j+1)\lambda_{s}\frac{\partial ^{j}}{\partial 
t_{s}^{j}}\right ] P_{L}(t)=\left\{\begin{matrix}\lambda 
P_{L}(t)&{\rm  for }\; j=0\\
0&{\rm  for }\; j>0\end{matrix} \right .\eqno{(*)}
$$
Clearly, differential operators correspond to polynomial solutions
$P_{L}(t)$ and in this case $\lambda=\sum_{s=1}^r \lambda_{s}-\deg
P_{L}$.  Kirillov demonstrated that nonpolynomial solutions do exist: for
$r=2$ and $\lambda_{1}=\lambda_{2}=0$ the function
$$
P_{L}(t)=\frac{t_{1}-t_{2}}{t_{1}+t_{2}}
$$
satisfies $(*)$ for $\lambda=0$.

{\bf Problem} \quad What invariant operator corresponds to this
solution?  Describe all (any) of the nonpolynomial solutions of $(*)$
and the corresponding operators.

4$^o$.  To select a reasonable type of $r$-nary operators is a good
problem.  Symmetric and skew-symmetric operators, as well as operators
on $\lambda$-densities are the first choices but even in such simple
cases there are few results.  These results, though scanty, are rather
interesting: quite unexpectedly, some of them are related to
calculation of the N.~Shapovalov determinant for the Virasoro algebra,
cf.\  Ref.\,\protect\citebk{F}.

5$^o$.  Since the real forms of simple vectorial Lie algebras are only
trivial ones (in the natural polynomial basis replace all complex
coefficients with reals), the results for $\Ree$ and $\Cee$ are
identical.  In super cases for nontrivial real forms some new
operators might appear; we will discuss this in the detailed version 
of the text.

\subsection{Arity $>1$} Grozman added a new dimension to Rudakov's
solution of Veblen's problem: in 1978 he described all {\it bi}nary
invariant differential operators.\cite{G1,G2} It turned out that there
are plenty of them but not too many: modulo dualizations and
permutations of arguments there are eight series of first order
operators and several second and third order operators all of which
are compositions of first order operators with one exception: the 3rd
order irreducible Grozman operator on the line.  There are no
invariant bilinear operators of order $>3$.

Miraculously, the 1st order differential operators determine, bar a
few exceptions, a Lie superalgebra structure on their domain.  (Here
Lie superalgebras timidly indicated their usefulness in a seemingly
nonsuper problem.  Other examples, such as Quillen's proof of the
index theorem, and several remarkable Witten's super observations
followed soon.)

{\bf Limits of applicability of Rudakov's method} \quad Though fans of
Rudakov's method,\footnote{~Interplay between restriction and
induction functors goes back to Frobenius, but discovery of each
instance deserves an acknowledgement, we presume.} let us point out
that its application to simple finite dimensional subalgebras of the
algebras of vector fields is extremely voluminous computational job;
therefore, it is ill applicable, say, to isometries of a Riemannian
manifold or the group preserving the Laplace operator.

Fortunately, when Rudakov's method fails, one can usually apply other
methods (Laplace-Casimir operators, N.~Shapovalov determinant, etc.).

\subsection{Generalized tensors and primitive forms} Rudakov
considered also operators invariant with respect to the Lie algebra of
Hamiltonian vector fields on the symplectic manifold $(M^{2n},
\omega)$.\cite{R2} Thanks to nondegeneracy of $\omega$ we can identify
$\Omega^i$ with $\Omega^{n-i}$.  So the operator $d^*:
\Omega^{n-i-1}\tto \Omega^{n-i}$, dual to the exterior differential
$d: \Omega^i\tto \Omega^{i+1}$ looks like a new operator, $\delta:
\Omega^{i+1}\tto \Omega^{i}$, the co-differential.  There are also
(proportional to each other) compositions
$\delta\circ\omega\circ\delta$ and $d\circ\omega^{-1}\circ d$, where
$\omega^{-1}$ is the convolution with the bivector dual to $\omega$.

A novel feature is provided by the fact that ``tensors'' now are
sections of the representation of $\fsp(V)$, not $\fgl(V)$.  Since
various representations of $\fsp(V)$ can not be extended to
representations of $\fgl(V)$ these ``tensors'' are, strictly speaking,
new notions.

Another novel feature we encounter considering subalgebras $\fg$ of
$\fvect$ are {\it primitive forms}.  If the $\fvect$-module $I(V)$
contains a singular vector with respect to $\fvect$, so it does with
respect to $\fg$.  But the irreducible $\fvect_{0}$-module $V$ does
not have to remain irreducible with respect to submodule of the
$\fg_{0}$.  The $\fg_{0}$-irreducible component with the biggest
highest weight in $V$ is called the $\fg_{0}$-{\it primitive}
(usually, just primitive) component.  Examples: the primitive
components appeared in symplectic geometry (we encounter their
counterparts in finite dimensional purely odd picture as spherical
harmonics,\cite{Sh2,LSH1}) contact analogues of primitive forms are
described in Ref.\,\protect\citebk{L1}.  Spaces of primitive
differential and integrable forms are just restrictions of the
``usual'' domain of the exterior differential; but other types of
primitive tensors are domains of really new invariant operators.

{\bf Further examples} A.~Shapovalov and Shmelev considered the Lie
{\it super}algebras of Hamiltonian vector fields and (following
Bernstein\,\cite{B} who considered the non-super case) their central
extension, the Poisson Lie superlagebra, see
Refs.\,\protect\citebk{Sh1}, \protect\citebk{Sh2} and
\protect\citebk{Sm2}-\protect\citebk{Sm4}.  Shmelev also considered
the operators invariant with respect to the funny exceptional
deformation $\fh_{\lambda}(2|2)$ of the Lie superalgebra $\fh(2|2)$ of
Hamiltonian vector fields.\cite{Sm1} See also
Ref.\,\protect\citebk{LSH2} for further description of
$\fh_{\lambda}(2|2)$.

By that time I.~Kostrikin described singular vectors for the contact
Lie algebras and found a ``new'' 2nd order invariant operator.  This
operator was actually well-known in differential geometry as an {\it
Euler operator} (for its description see Ref.\,\protect\citebk{Ly}; 
here we just
briefly observe that it is not $\sum x_{i}\partial_{i}$, this is
another Euler operator); it is needed for invariant formulation of
Monge-Amp\'ere equations, cf.\ Ref.\,\protect\citebk{LRC}. 
Leites\,\cite{L1} generalized I.~Kostrikin's calculations to contact
Lie {\it super}algebras and found out that there seem to be no
analogue of Euler's operator in supersetting.  This makes one
contemplate on the following:

{\bf Problem}\quad What are superanalogs of Monge-Amp\'ere equations,
if any?

In 1977 ``odd'' analogs of the hamiltonian and contact series were
discovered.\cite{L0,ALSh} Batalin and Vilkovisky rediscovered the
antibracket related to these series and showed its importance, cf.\
Ref.\,\protect\citebk{GPS}.  Kochetkov\,\cite{Ko1}$^-$\cite{Ko5}
undertook the task of calculating the corresponding singular vectors. 
He digressed to consider two of the three known at that time
Shchepochkina's exceptions\,\cite{Ko0} (and named after her with the
first Cyrillic letter of her name), one of which was recently
reconsidered in another realization in Ref.\,\protect\citebk{KR1}.  At
the last moment we have found out more singular vectors ($=$ invariant
operators) than in Ref.\,\protect\citebk{Ko0}, even with finite
dimensional fibers; SuperLie is indeed indispensable.  After the last
moment we have considered {\cyr shch}${}_{2}$; we reproduce
Kochetkov's result for comparison.  For the lack of resources we left
out some possible cases of singular vectors, but we are sure they are
improbable.  Though for $\fv\fle(4|3)$ and $\fm\fb(4|5)$ all degrees
can occure, we are sure induction \`a la Kochetkov (complete list of
singular vectors) can be performed.  Various regradings of
$\fk\fs\fle$ seem to be a tougher problem.

\subsection{Superization leads to new developments} The study of
invariant differential operators on supermanifolds began in 1976 as a
byproduct of attempts to construct an integration theory on
supermanifolds similar to the integration theory of differential forms
on manifolds.  Bernstein and Leites became interested in Veblen's
problem when they tried to construct an integration theory for
supermanifolds containing an analog of the Stokes
formula.\cite{BL1,BL2} At that time there were only known the
differential forms which are impossible to integrate and the volume
forms of the highest degree.  Unlike the situation on manifolds, no
volume form coincides with any differential form and there was known
no analogs of volume forms of lesser degrees.

Having discovered {\it integrable} forms\,\cite{BL1,BL3} (i.e., the
forms that can be integrated; Deligne calls them {\it integral}
forms\,\cite{D}) Bernstein and Leites wanted to be sure that there
were no other tensor objects that can be integrated.  Observe several
points of this delicate question.

(1) The conventional Stokes formula on a manifold exists due to the
fact that there is an invariant operator on the space of differential
forms.  The uniqueness of the integration theory with Stokes formula
follows then from the above result by Rudakov and its superization due
to Bernstein and Leites.

Since there are several superanalogs of the determinant, it follows
that on supermanifolds, there are, perhaps, several analogs of
integration theory, see Ref.\,\protect\citebk{LInt}, some of them
without Stokes formula.  Still, if we wish to construct an integration
theory for supermanifolds containing an analog of the Stokes formula,
and, moreover, coinciding with it when the supermanifold degenerates
to a manifold, we have to describe all differential operators in
tensor fields on supermanifolds.

(2) Bernstein and Leites confined themselves to {\it finite
dimensional} representations $\rho$ owing to tradition which says that
{\sl a tensor field is a section of a vector bundle with a finite
dimensional fiber on which the general linear group acts}.  Even in
doing so Bernstein and Leites had to digress somewhat from the
conventions and consider, since it was natural, ALL finite dimensional
irreducible representations $\rho$ of the general linear Lie
superalgebra.  Some of such representations can not be integrated to a
representation of the general linear supergroup.

Inspired by Duflo, Leites used calculations of
Ref.\,\protect\citebk{BL3} to describe invariant differential operators
acting in the superspaces of tensor fields with {\it infinite}
dimensional fibers, see Ref.\,\protect\citebk{LKV}.  These operators of
order $>1$ are totally new, though similar to fiberwise integration
along the odd coordinates.  The operators of order 1 are also not bad:
though they are, actually, the good old exterior differential $d$, the
new domain is that of semi-infinite forms, certain class of
pseudodifferential forms.  Observe that quite criminally (in
V.~I.~Arnold's words) no example of the corresponding new type
homology is calculated yet, except some preliminary (but important)
results of Shander, see Ref.\,\protect\citebk{L2} v.\ 31, Ch.\ 4, 5.

(3) Even under all the restrictions Bernstein and Leites imposed, to
say that ``the only invariant differential operator is just the
exterior differential'' would be to disregard how drastically they
expanded its domain (even though they ignored semi-infinite
possibilities).  It acts in the superspace of differential forms and
in the space of {\it integrable forms}, which is natural, since the
space of integrable forms is just the dual space to the superspace of
differential forms.  Though Bernstein and Leites did not find any new
invariant differential operator ({\it this proves that an integration
theory on supermanifolds containing an analog of the Stokes formula
can only be constructed with integrable forms}), they enlarged the
domain of the exterior differential to the superspace of
pseudodifferential and pseudointegrable forms.  These superspaces are
not {\it tensor fields} on $\cM^{m, n}$ unless $n=1$, but they are always
tensor fields on the supermanifold $\widehat{\cM}$ whose structure
sheaf $\cO_{\widehat{\cM}}$ is a completion of the sheaf of
differential forms on $\cM$; namely, the sections of
$\cO_{\widehat{\cM}}$ are arbitrary functions of differentials, not
only polynomial ones.

(4) Bernstein and Leites did not consider indecomposible
representations $\rho$ which are more natural in both the supersetting
and for infinite dimensional fibers.  The first to consider
indecomposible cases was Shmelev;\,\cite{Sm4} his result was, however,
``not interesting'': there are no totally new operators, just
compositions of the known ones with projections.  For a review of
indecomposible representations of simple Lie superalgebras see 
Ref.\,\protect\citebk{LInd}.

{\bf Integration and invariant differential operators for
infinite dimensional fibers}\quad  There are new operators invariant with
respect to the already considered (super)groups of diffeomorphisms or,
equivalently, their Lie superalgebras, if we let them act in the
superspaces of sections of vector bundles with infinite dimensional
fibers.  These operators of high order have no counterparts on
manifolds and are versions of the Berezin integral applied fiber-wise. 
(A year after the talk with these results\,\cite{LKV} was delivered,
I.~Penkov and V.~Serganova interpreted some of these new operators as
acting in the superspaces of certain tensor fields on \lq\lq curved"
superflag and supergrassmann supervarieties.\cite{PS})

We hope to relate with some of these operators new topological
invariants (or perhaps old, like cobordisms, but from a new
viewpoint).  Recall that since the de Rham cohomology of a
supermanifold are the same as those of its underlying manifold, the
\lq\lq old type" operators are inadequate to study \lq\lq topological"
invariants of supermanifolds.  The operators described here and
related to vector bundles of infinite rank lead to new (co)homology
theories (we prefix them with a ``pseudo'').  These pseudocohomologies
provide us with invariants different from de Rham cohomology; 
regrettably, never computed yet.

The approach adopted here for the operators in the natural bundles
with {\it infinite dimensional} fibers on supermanifolds prompts us to
start looking for same on manifolds. From the explicit
calculations in Grozman's thesis,\cite{G1} it is clear that there are
some new bilinear operators acting in the spaces of sections of tensor
fields with infinite dimensional fibers.

\subsection{An infinitesimal version of Veblen's problem} 
Denote $\cF=\Kee [[x]]$, where $x=(u_{1}, \dots, u_{n},
\xi_{1} , \dots, \xi_{m})$ so that $p(u_{i})=\ev$ and
$p(\xi_{j})=\od$.  Denote by $(x)$ the maximal ideal in $\cF$
generated by the $x_{i}$.  Define a topology in $\cF$ so that the
ideals $(x)^r$, $r=0, 1, 2, \dots$ are neighborhoods of zero, i.e.,
two series are $r$-{\it close} if they coincide up to order $r$.  We see
that $\cF$ is complete with respect to this topology.

Denote by $\fvect(n|m)$ the Lie superalgebra of formal vector fields,
i.e., of continuous derivations of $\Kee[[x]]$.  By abuse of notations
we denote $\fder \Kee[x]$, the Lie superalgebra of {\it
polynomial} vector fields, also by $\fvect(n|m)$.

Define partial derivatives $\partial_{i}=\pder{x_{i}}\in \fvect(n|m)$
by setting $\partial_{i}(x_{j})=\delta_{ij}$ with super-Leibniz rule. 
Clearly, $p(\partial_{i})=p(x_{i})$ and $[\partial_{i},
\partial_{j}]=0$.  Any element $D\in \fvect(n|m)$ is of the form
$D=\sum f_{i}\partial_{i}$, where $f_{i}=D(x_{i})\in\cF$.  We will
denote $\fvect(n|m)$ by $\cL$.  In $\cL$, define a filtration of the
form $\cL=\cL_{-1}\supset\cL_{0}\supset\cL_{1}\supset\dots$ setting
$$
 \cL_{r}=\{D\in \fvect(n|m)\mid D(\cF)\subset
(x)^{r+1}\}.
$$
This filtration defines a topology on $\cL$, the
superspaces $\cL_r$ being the base of the topology, open
neighborhoods of zero.

Denote by $L=\oplus L_r$, where $L_{r}=\cL_{r}/\cL_{r+1}$, the
associated graded Lie superalgebra.  Clearly, $L_0\simeq\fgl(n|m)$ 
with $E_{ij}\longleftrightarrow x_{j}\partial_i$.

Let $\rho$ be an irreducible representation of the Lie superalgebra
$L_{0}=\fgl(n|m)$ with lowest weight in a superspace $V$.  Define a
$\fvect(n|m)$-module $T(\rho)$ also denoted by $T(V)$ by setting
$T(V)=\cF\otimes_{\Kee}V$.  The superspace $T(V)$ evidently inherits
the topology of $\cF$.  To any vector field $D$, assign the operator
$L_{D}: T(V)\tto T(V)$ --- the Lie derivative --- such that for
$f\in\cF$ and $v\in V$
$$
L_{D}(fv)=D(f)v+(-1)^{p(D)p(f)} \sum D^{ij}\rho(E_{ij})(v), 
\eqno{(1.7.1)}
$$
where $D^{ij}=(-1)^{p(x_i)(p(f)+1)}\partial_i
f_j$.  We will usually write just $D$ instead of $L_{D}$.

The elements $t\in T(V)$ will be called {\it tensor fields of type
$V$}.  The modules $T(V)$ are topological; their duals are spaces with
discrete topology.  

Observe that even if $V$ is finite dimensional, the elements of $T(V)$ 
are generalized tensors as compared with the classical notion: the 
space $V$ might not be realized in the tensor product of co- and 
contra-variant tensors, only as a subquotient of such; e.g., unlike 
the determinant (or trace, speaking on the Lie algebra level), the
supertrace is not realized in tensors and we have to introduce new 
type of ``tensors'' --- the $\lambda$-densities.

For any $L_0$-module $V$ with highest weight and any $L_0$-module $W$
with lowest weight set
$$
I(V)=U(\cL)\otimes_{U(\cL_0)}V; \quad T(W)=\Hom_{U(\cL_0)}(U(\cL),  
W),\eqno{(1.7.2)}
$$ 
where we have extended the action of $L_0$ to a
$U(\cL_0)$-action by setting $\cL_1V=0$ and $\cL_1W=0$.  Clearly,

a) $I(V)$ is an $\cL$-module with discrete topology.

b) $I(V)^*\cong T(V^*)$

c) definition of the tensor fields with $\cL$-action (1.7.1) is
equivalent to the one given by (1.7.2).

Thus, instead of studying invariant maps $T(W_1)\tto T(W_2)$ (or
$T(W_1)\otimes T(W_2)\tto T(W_3)$, etc.)  we may study submodules ---
or, equivalently, singular vectors ---
of
$I(V)$ (resp.  of $I(V_1)\otimes I(V_2)$, etc.). They are much
easier to describe.

{\bf Further generalization of tensors. The highest 
weight theorem} Let 
$$
\cL=\cL_{-d}\supset\dots \supset\cL_{0}\supset\cL_{1}\supset\dots
$$
be a Lie superalgebra of vector fields with
formal or polynomial coefficients and endowed with a Weisfeiler
filtration/grading described in what follows (for the time being
consider a ``most natural'' grading, like that in $\fvect$ above).  We
define the space of {\it generalized} tensor fields and its dual by
the same formula (1.7.2) as for the usual tensor fields given any
$L_0$-module $V$ with highest weight and any $L_0$-module $W$ with
lowest weight such that $\cL_1V=0$ and $\cL_1W=0$.

Observe that for the Lie algebra of divergence-free vector fields the
spaces $T(W)$ are the same as for $\fvect$.  For some other Lie
superalgebras the notion of tensors we give is different because there
are representations of $L_{0}$ distinct from tensor powers of the
identity one.  For example, for the Lie superalgebra $\cL$ of
Hamiltonian vector fields $\fh(2n|m)$ such is the spinor
representation (for $n=0$); if we consider infinite dimensional fibers
such is the oscillator representation (for $m=0$), and in the general
case such is the spinor-oscillator representation, cf.\  Ref.\,\protect\citebk{LSH1}.

Thus, the first step in the study of $\cL$-invariant operators is a
description of irreducible $L_{0}$-modules, at least in terms of the
highest/lowest weight.  For the majority of the $L_{0}$'s this is not
a big deal, but the catch is that for some $L_{0}$'s there is no easy
highest/lowest weight theorem, even for finite dimensional modules. 
We will encounter this phenomenon with $\fa\fs$, the linear part of
$\fv\fa\fs$.

An aside remark: being interested not only in representations of
vectorial algebras (with polynomial coefficients) but in their stringy
analogs (with Laurent coefficients), too, observe that vacuum over
$L_{0}$ can be degenerate.

{\bf Problem} \quad For all Weisfeiler gradings of simple vectorial 
superalgebras $\cL$ describe conditions for the highest (lowest) weight
under which the irreducible quotient of the Verma module over $L_{0}$
is finite dimensional and describe the corresponding module (say, in 
terms of a character formula, cf.\  Ref.\,\protect\citebk{PS}). 

{\bf Examples of generalized tensor fields} \quad Clearly, for $\cL=
\fvect(n|m)$ we have $\cL\equiv T(\id)$, where
$\id=\Span(\partial_{i}\mid 1\leq i\leq n+m)$ is the (space $V$ of
the) identity representation of $L_{0}=\fgl(V)=\fgl(n|m)$.  The spaces
$T(E^i(\id^*))$ are denoted by $\Omega^i$; their elements are called
{\it differential $i$-forms} and the right dual elements to
$\partial_{i}$ are denoted by $\widehat{x_{i}}=dx_{i}$, where
$p(\widehat{x_{i}})=p(x_{i})+\od$.  In particular, let $\cF=\Omega^0$
be the algebra of functions.

The algebra $\widehat{\Omega}$ of arbitrary, not only polynomial,
functions in $\widehat{x_{i}}=dx_{i}$ is called the algebra of {\it
pseudodifferential forms}.  An important, as Shander showed in
Ref.\,\protect\citebk{L2} v.~31, Ch.~5, subspace
$\widehat{\Omega}^{(\lambda)}$ of homogeneous {\it pseudodifferential
forms of homogeneity degree} $\lambda\in\Kee$ is naturally defined as
functions of homogeneity degree $\lambda$ with respect to the hatted
indeterminates.

Define the space of volume forms $\Vol$ to be $T(\str)$; denote the
volume element by $\vvol(x)$ or $\vvol(u|\theta)$.  (Observe again
that it is a bad habit to denote, as many people still do, $\vvol$ by
$d^nud^m\theta$: their transformation rules are totally different,
see, e.g., Refs.\,\protect\citebk{BL3}, \protect\citebk{D}.)

The space of {\it integrable} $i$-forms is $\Sigma_{i}=
\Hom_{\cF}(\Omega^i, \Vol)$.  In other words, integrable forms are
$\Vol$-valued polyvector fields.  Pseudointegrable forms are defined
as elements of $\widehat{\Sigma}= \Hom_{\cF}(\widehat{\Omega}, \Vol)$;
the subspace $\widehat{\Sigma}_{(\lambda)}=
\Hom_{\cF}(\widehat{\Omega}^{(\lambda)}, \Vol)$ of homogeneous forms
is also important.

{\bf Particular cases:}  
a) $m=0$.  We see that $\Omega^i=0$ for $i>n$
and $\Sigma_i=0$ for $i<0$.  In addition, the mapping
$\vvol\mapsto\hat x_1\cdots\hat x_n$ defines an isomorphism of
$\Omega^i$ with $\Sigma_i$ preserving all structures.

b) $n=0$.  In this case there is an even $\cL$-module morphism $\int:
\Sigma_{-m}\tto\Kee$ called the {\it Berezin integral}.  It is
defined by the formula 
$$
\int\xi_1\cdots\xi_m\,\vvol=1, {\rm  and } \int\xi^{\nu_1}_1
\cdots\xi^{\nu_m}_m\, \vvol=0 \quad{\rm if}\quad \prod\nu_i=0.
$$

We will also denote by $\int$ the composition $\int: \Sigma_{-m}\to
\Kee\hookrightarrow\Omega^0 $ of the Berezin integral and the natural
embedding.

c) $m=1$.  We generalize $\Omega^i$ and $\Sigma_j$ to the spaces
$\Phi^\lambda$ of pseudodifferentional and pseudointegral forms
containing $\Omega^i$ and $\Sigma_j$, where $\lambda\in \Kee$.  Let
$x=(u_1, \cdots, u_n, \xi)$.  Consider a $\Kee$-graded $\Omega$-module
$\Phi=\oplus\Phi^\lambda$ (we assume that $\deg\hat x_i=1\in \Kee$)
generated by $\hat \xi^\lambda$, where $\deg\hat \xi^\lambda=\lambda$
and $p(\hat \xi^\lambda)=\ev$, with relations
$\hat\xi\cdot\hat\xi^\lambda=\hat\xi^{\lambda+1}$.  Define the action
of partial derivatives $\partial_i$ and $\hat\partial_j$ for $1\leq i,
j\leq n+1$ via $\hat\partial_j(x_i)=0$,
$\partial_i(\hat\xi^\lambda)=0$, $\partial_{\hat
u_i}(\hat\xi^\lambda)=0$ and
$\partial_{\hat\xi}(\hat\xi^\lambda)=\lambda\hat\xi^{\lambda-1}$.

On $\Phi$, the derivations $d$, $i_{\cD}$ and $L_\cD$ consistent
with the exterior derivation $d$, the inner product $i_{\cD}$ and the Lie
derivative $L_\cD$ on $\Omega$ are naturally defined.

It is easy to see that $\Phi=\oplus \Phi^{\lambda}$ is a
supercommutative superalgebra. 

Clearly, $\Phi$ is a superspace of tensor fields and for
$\Phi^{\Zee}=\oplus_{r\in\Zee}\Phi^r$ we have a sequence
$$
0\tto\Omega\stackrel{\alpha}{\tto}\Phi^\Zee \stackrel{\beta}{\tto}\Sigma\tto
0\eqno{(*)}
$$
where the maps $\alpha$ and $\beta$ are defined by
$$
\alpha(\omega)=\omega\hat\xi^0, \quad \beta(\hat u_1\cdots\hat
u_n\hat\xi^{-1})=\vvol. 
$$
Clearly, the homomorphisms $\alpha$ and $\beta$ are consistent with
the $\Omega$-module structure and the operators $d$, $i_\cD$ and
$L_\cD$.  The explicit form of the $\cF$-basis in $\Omega$, $\Sigma$
and $\Phi$ easily implies that $(*)$ is exact.

\vspace*{-2mm}

\subsection{Operators invariant with respect to nonstandard
realizations} At the moment the $\cL$-invariant differential operators
are described for all but one series of simple vectorial Lie
superalgebras in the standard realization.  Contrariwise, about
operators invariant with respect to same in nonstandard realizations
almost nothing is known, except for $\fvect(m|n; 1)$, see
Ref.\,\protect\citebk{LKV}.

For series, the standard realization is the one for which $\dim
\cL/\cL_{0}$ is minimal; for exceptional algebras the notion of the
standard realization is more elusive, and since there are 1 to 4
realizations, it is reasonable and feasible to consider all of them. 
It is also natural to consider $\fh_{\lambda}(2|2)$ and
$\fh_{\lambda}(2|2; 1)$ as exceptional algebras, especially at
exceptional values of $\lambda$.

\vspace*{-2mm}

\subsection{On description of irreducible $\cL$-modules} Having
described $\fvect(n|m)$-invariant differential operators in tensor fields
with finite dimensional fibers (answer: only $d$, and $\int$ if $n= 0$,
$m\neq 0$), we consider the quotients of $T(V)$ modulo the image of
the invariant operator.  It could be that the quotient also contains a
submodule.  In the general case there are no such submodules
(Poincar\'e lemma), in other cases anything can happen, see 
Refs.\,\protect\citebk{Ko4} and \protect\citebk{Ko3,Ko5}.

Observe that to describe irreducible $\cL$-modules, it does not always
suffice to consider only one realization of $\cL$.  It is like
considering generalized Verma modules induced or co-induced from
distinct parabolic subalgebras.  Similarly, the description of
invariant operators must be performed from scratch in each
realization. 

Here we do not specifically consider the irreducible $\cL$-modules; so
far, the answers are known for tensors with finite dimensional fibers and
in two cases only: Ref.\,\protect\citebk{KR2} ($\fv\fle(3|6)$) and
Ref.\,\protect\citebk{L1} ($\fk(1|n)$; weights of singular vectors are
corrected below).

\vspace*{-2mm}

\section{Brief description of the exceptional algebras}

To save space, we do not reproduce the details of definitions of
exceptional superalgebras, see Ref.\,\protect\citebk{S4}. 
Thus,\cite{R1,Ko3,BL3} suffice to grasp the details of the
theory;\,\cite{S4} to catch on with the list of exceptions we
consider.

Some of the exceptional algebras $\fg=\mathop{\sum}_{i\geq -d}\fg_{i}$
are isomorphic as abstract ones; there are five abstract families
altogether.  We realize them as Lie superalgebras of polynomial vector
fields with a particular, Weisfeiler (shortly W-), grading or
filtration.  In any W-grading (a) the sum of the terms of positive
degree is a maximal subalgebra of finite codimension and (b) the
linear part, $\fg_{0}$ irreducibly acts on $\fg_{-1}$.  If depth
$d=1$, then each $\fg$ is constructed as the Cartan prolong (or its
generalization) of the pair $(\fg_{-1}, \fg_{0})$.  To construct
$\fk\fa\fs$, still another generalization of the Cartan
prolongation is applied.  For these generalizations, first described
in Refs.\,\protect\citebk{S1} and \protect\citebk{S2}, see 
Ref.\,\protect\citebk{S4}.  In
Ref.\,\protect\citebk{CK} there is given a 15-th regrading we missed in
Ref.\,\protect\citebk{S4}, and a nice interpretation of the exceptional
algebras in terms of $\fg_{\ev}$ and the $\fg_{\ev}$-module,
$\fg_{\od}$.  This interpretation is convenient in some problems, but
in interpretations of our results the realization of $\fg$ as
(generalized) Cartan prolong is often more useful.

We only recall that the classical Cartan prolong $(\fg_{-1},
\fg_{0})_*$ is defined inductively, as a subalgebra of vector fields
in $\dim\fg_{-1}$ indeterminates with given $\fg_{-1}$ and the linear
part $\fg_{0}$, and where for $i>0$
$$
\fg_{i+1}=\{D\in\fvect(\dim\fg_{-1})\mid [D, \fg_{-1}]\subset \fg_{i}.
$$

Hereafter $\fvect(m|n)$ is the general vectorial Lie superalgebra on $m$ 
even and $n$ odd indeterminates; $\fsvect(m|n)$ is its divergence-free or special subalgebra; $\fk(2m+1|n)$ the contact algebra that 
preserves the Pfaffian equation $\alpha=0$, where 
$$
\alpha_1 =dt+\mathop{\sum}\limits_{1\leq i\leq
m}(p_idq_i-q_idp_i)+\mathop{\sum}\limits_{j\leq r}
(\xi_jd\eta_j+\eta_jd\xi_j)+\mathop{\sum}\limits_{k\geq
n-2r}\theta_kd\theta_k.
$$
For $f\in\Cee [t, p, q, \Theta]$, where $\Theta=(\xi, \eta, \theta)$, define 
the contact field $K_f$ by setting:
$$
K_f=(2-E)(f)\pder{t}-H_f + \pderf{f}{t} E,
$$
where $E=\sum\limits_i y_i \pder{y_{i}}$ (here the $y_{i}$ are all the 
coordinates except $t$), and $H_f$ is the hamiltonian field 
with Hamiltonian $f$ that preserves $d\alpha_1$:
$$
H_f=\sum\limits_{i\leq n}(\pderf{f}{p_i} \pder{q_i}-\pderf{f}{q_i}
\pder{p_i}) -(-1)^{p(f)}\left(\sum\limits_{j\leq m}\pderf{
f}{\theta_j} \pder{\theta_j}\right ).
$$
Replacement of some of the $\theta$'s in the above formula with $\xi,
\eta$ leads to obvious modifications.

The ``odd'' analog of the contact structure is given by the {\it even} form
$$
\alpha_0=d\tau+\sum\limits_{1\leq i\leq n}(\xi_idq_i+q_id\xi_i)
$$
and formula for the vector field $M_f$ --- the analogs of $K_f$ ---
generated by the function of $\tau, \xi, q$ are similar. The fields 
$M_f$ span the ``odd'' contact Lie superalgebra, $\fm(n)$.

The regradings are given after the dimension of the supermanifold in
the standard realization ($r=0$, optional, marked by $(*)$) after
semicolon.  Observe that the codimension of ${\cal L}_0$ attains its
minimum in the standard realization.
  
\tiny
$$
\renewcommand{\arraystretch}{1.3}
\begin{tabular}{|c|c|}
\hline
Lie superalgebra & its $\Zee$-grading \\ 
\hline
$\fvect (n|m; r)$, & $\deg u_i=\deg \xi_j=1$  for any $i, j$
\qquad
$(*)$\\ 
\cline{2-2}
$ 0\leq r\leq m$ & $\deg \xi_j=0$ for $1\leq j\leq r;$\\
&$\deg u_i=\deg \xi_{r+s}=1$ for any $i, s$ \\ 
\hline
\hline$\fm(n; r),$& $\deg \tau=2$, $\deg q_i=\deg \xi_i=1$  for any $i$ 
\qquad $(*)$\\ 
\cline{2-2}
$\; 0\leq r\leq n$& $\deg \tau=\deg q_i=2$, $\deg \xi_i=0$ for $1\leq i\leq r
<n$;\\ 
$r\neq n-1$& $\deg q_{r+j}=\deg \xi_{r+j}=1$ for any $j$\\ 
\hline
$\fm(1|n; n)$ & $\deg \tau=\deg q_i=1$, $\deg \xi_i=0$ for $1\leq i\leq n$ \\
\hline
\hline
\cline{2-2}
$\fk (2n+1|m; r)$, & $\deg t=2$\,, \\
$0\leq r\leq [\frac{m}{2}]$& $\deg p_i=\deg q_i=
\deg \xi_j=\deg \eta_j=\deg \theta_k=1$ for any $i, j, k$ \qquad $(*)$\\ 
\cline{2-2}
$r\neq k-1$ for $m=2k$& $\deg t=\deg \xi_i=2$, $\deg 
\eta_{i}=0$ for $1\leq i\leq r\leq [\frac{m}{2}]$; \\
%\cline{2-2}
and $n=0$&$\deg p_i=\deg q_i=\deg \theta_{j}=1$ for $j\geq 1$ and all $i$\\ 
\hline
$\fk(1|2m; m)$ & $\deg t =\deg \xi_i=1$, $\deg 
\eta_{i}=0$ for $1\leq i\leq m$ \\ \hline
\end{tabular}
$$
%\vskip 0.2 cm

\normalsize 

Now the last notations: $\fhei(2n|m)$ is the {\it Heisenberg algebra},
it is isomorphic to $\fk(2n+1|m)_-$, the negative part of
$\fk(2n+1|m)$; similarly, $\fab(m)$ is the {\it antibracket algebra},
isomorphic to $\fm(m)_-$.

A.~Sergeev central extension, $\fas$, of $\fspe(4)$: we represent an 
arbitrary element $A\in\fas$ as a pair $A=x+d\cdot z$, where 
$x\in\fspe(4)$, $d\in{\Cee}$ and $z$ is the central element.  We 
define the
bracket in $\fas$ in its matrix realization in the spinor 
representation:
%\vspace*{-2mm}
$$
\left[\begin{pmatrix} a & b \cr 
c & -a^t \end{pmatrix}\!+\!d\cdot z\,,\, \begin{pmatrix}
a' & b' \cr 
c' & -a'{}^t \end{pmatrix}\! +\!d'\cdot z\right]\!=\!
\left[\begin{pmatrix} a & b \cr 
c & -a^t \end{pmatrix}\,,\, \begin{pmatrix}
a' & b' \cr 
c' & -a'{}^t \end{pmatrix}\right]\!+\tr~c\tilde c'\cdot z\,,
$$
where $\, \tilde{}\, $ is extended via linearity from matrices
$c_{ij}=E_{ij}-E_{ji}$ on which $\tilde c_{ij}=c_{kl}$
for any even permutation $(1234)\mapsto(ijkl)$.

The five types of exceptional W-graded Lie superalgebras are given
above in their realizations as Cartan's prolongs $(\fg_{-1},
\fg_{0})_{*}$ or generalized Cartan's prolongs $(\fg_{-},
\fg_{0})_{*}^{mk}$ for $\fg_{-}=\mathop{\oplus}\limits_{-d\leq i\leq
-1}\fg_{i}$ expressed for $d=2$ as $(\fg_{-2}, \fg_{-1},
\fg_{0})_{*}^{mk}$ together with one of the serial Lie superalgebras
as an ambient which contains the exceptional one as a maximal
subalgebra. The regradings $R(r)$ of the ambients given
below are sometimes not of Weisfeiler type.
$$
\renewcommand{\arraystretch}{1.4}
\begin{tabular}{|l|}
\hline
$\fv\fle(4|3; r)=(\Pi(\Lambda(3))/\Cee\cdot 1, 
\fc\fvect(0|3))_{*}\subset\fvect(4|3; R), \quad r= 0, 1, K$\cr
\hline
$\fv\fas(4|4)=(\spin, \fas)_{*}\subset \fvect(4|4)$\cr
\hline
$\fk\fas^\xi(1|6; r)\subset \fk(1|6; r),\quad r=0, 1\xi, 
3\xi$\cr 
$\fk\fas^\xi(1|6; 3\eta)=(\Vol_{0}(0|3), 
\fc(\fvect(0|3)))_{*}\subset \fsvect(4|3)$\cr 
\hline
$\fm\fb(4|5; r)=(\fab(4), \fc\fvect(0|3))_{*}^m\subset \fm(4|5; R), \quad 
r=0, 1, K$\cr
\hline
$\fk\fsle(9|6; r)=(\fhei(8|6), \fsvect(0|4)_{3, 4})_{*}^k\subset 
\fk(9|6; r), \quad r=0, 2$, CK\cr 
$\fk\fsle(9|6; K)=(\id_{\fsl(5)}, 
\Lambda^2(\id^*_{\fsl(5)}), \fsl(5))_{*}^k\subset \fsvect(5|10; R)$\cr 
\hline
\end{tabular}
$$

Certain regradings $R(r)$ of the ambients are so highly nonstandard 
that even the homogeneous fibers are of infinite dimension:

\noindent 1) $\fv\fle(4|3; r)=(\Pi (\Lambda(3)/{\Cee}\cdot 1), \fcvect
(0|3))_*\subset \fvect(4|3), \quad r= 0, 1, K$;

\noindent ~~\underline{$r= 0$}: $\deg y= \deg u_i= \deg\xi_i=1$;

\noindent ~~\underline{$r= 1$}: $\deg y =\deg\xi_1=0$, $ \deg u_2=\deg 
u_3=\deg\xi_2=\deg\xi_3=1$, $\deg u_1=2$;

\noindent ~~\underline{$r= K$}: $\deg y=0$, $\deg u_i=2$; $\deg\xi_i=1$.

\noindent 2) $\fv\fas(4|4)=(\spin, \fas)_*\subset \fvect(4|4)$;

\noindent 3) $\fk\fas\subset \fk(1|6; r)$,\; $r=0, 1, 3\xi$; \quad 
$\fk\fas(1|6; 3\eta)\subset \fsvect(4|3)$;

\noindent ~~\underline{$r=0$}: $\deg t=2$, \; $\deg \eta_i=1$; \; 
$\deg\xi_i=1$; 

\noindent ~~\underline{$r=1$}: $\deg\xi_1\!=\!0$,\, 
$\deg\eta_1\!=\!\deg t\! =\!2$,\, 
$\deg\xi_2\!=\!\deg\xi_3\!=\!\deg\eta_2\!=\! \deg\eta_3=1$; 

\noindent ~~\underline{$r=3\xi$}: $\deg\xi_i=0$, \; $\deg\eta_i=\deg t=1$; 

\noindent ~~\underline{$r=3\eta$}: $\deg\eta_i=0$, $\deg\xi_i=\deg t=1$.

\noindent 4) $\fm\fb(4|5; r)=(\fab(4), \fcvect (0|3))_*^{m}\subset \fm(4), \quad 
r=0, 1, K$;

\noindent ~~\underline{$r\!=\! 0$}: $\deg \tau\!=\!2$, \; 
$\deg u_i\!=\! \deg\xi_i=1$ for   
$i=0$, 1, 2, 3;

\noindent ~~\underline{$r\!=\! 1$}: $\deg \tau\!=\!\deg \xi_0\!=\!\deg u_1\!
=\!2$,
\, $\deg u_2\!=\deg u_3\! =\!\deg\xi_2\!=\!\deg\xi_3\!=\!1$;

\noindent ~~$\deg\xi_1\!=\!\deg u_0\!=\!0$;

\noindent ~~\underline{$r= K$}: $\deg \tau=\deg \xi_0=3$, \; $\deg u_0=0$, \; 
$\deg u_i=2$; \; $\deg\xi_i=1$ for $i>0$.

\noindent 5) $\fk\fsle(9|6; r)=(\fhei(8|6), \fsvect _{3, 4}(4))_*^{k} \subset 
\fk(9|6), \quad r=0$,  2,  K, CK;

\noindent ~~\underline{$r\!=\! 0$}: $\deg t\!=\!2$, \; $\deg p_i\!=\! \deg q_i\!=\! \deg\xi_i\!=\! 
\deg\eta_i=1$; 

\noindent ~~\underline{$r\!=\!2$}: $\deg t\!=\!\deg q_3\!=\!\deg q_4\!=\!\deg\eta_1\!=\! 2$,\, $\deg 
q_1\!=\!\deg q_2\!=\! \deg p_1\!=\!\deg p_2$

\noindent ~~$\!=\!\deg\eta_2\!=\!\deg\eta_3=\deg\zeta_2\!=\!\deg\zeta_3\!=\!1$; \, $\deg p_3\!=\!\deg p_4\!=\!\deg\zeta_1\!=\!0$; 

\noindent ~~\underline{$r\!=\! K$}: $\deg t\!=\!\deg q_i=2$,\, $\deg p_i\!=\!0$;  \,
$\deg\zeta_i\!=\!\deg\eta_i=1$; 

\noindent ~~\underline{$r\!=\!CK$}: $\deg t\!=\!\deg q_1=3$,\,$\deg p_1\!=\!0$;\,$\deg q_2\!=\!\deg q_3\!=\! \deg q_4\!=\!\deg\zeta_1$

\noindent ~~$=\!\deg\zeta_2\!=\!
\deg\zeta_3\!=\!2$;\,$\deg p_2\!=\!\deg p_3\!=\! \deg p_4\!=\!
\deg\eta_1\!=\!\deg\eta_2\!=\!\deg\eta_3\!=\!1$.

To determine the minimal ambient  is important for
our problem: every operator invariant with respect to an algebra is,
of course, invariant with respect to any its subalgebra.  

Here is the list of nonpositive terms of the exceptional algebras.
Notations: $\fc(\fg)$ denotes the trivial central
extension with the 1-dimensional even center generated by $z$;
$\Cee[k]$ is the trivial $\fg_{0}$-module on which the central element
$z$ from $\fg_{0}$ chosen so that $z|_{\fg_i}=i\cdot \id_{\fg_i}$ acts
as multiplcation by $k$; $\fa\subplus \fb$ denotes the semidirect
product in which $\fa$ is the ideal; let $d$ determine the
$\Zee$-grading of $\fg$ and not belong to $\fg$, we shorthand
$\fg \; +\hspace{-3.6mm}\supset \Cee(az+bd)$ to $\fg_{a, b}$; $\Vol(n|m)$ is the space of
volume forms (densities) on the superspace of superdimension indicated,
subscript $0$ singles the subspace with integral $0$; $T^{1/2}$ is the
representation in the space of half-densities and $T^{0}_0$ is the
quotient of $\Vol_0$ modulo constants (over divergence-free algebra).

Observe that none of the simple W-graded vectorial Lie superalgebras
is of depth $>3$ and only two algebras are of depth 3: $\fm\fb(4|5;
K)$, for which we have
$$
\fm\fb(4|5; K)_{-3}\cong \Pi(\id_{\fsl(2)}),
$$
and another one, $\fk\fsle(9|6; CK)=\fc\fk(11|9)$, for which we have :
$$
\fc\fk(11|9)_{-3}\simeq \Pi(\id_{\fsl(2)}\otimes\Cee[-3]). 
$$ 
Here are the other terms $\fg_{i}$ for $i\leq 0$ of the 15 exceptional
W-graded algebras.  

\tiny
$$
\renewcommand{\arraystretch}{1.3}
\begin{tabular}{|c|c|c|c|c|}
\hline
$\fg$&$\fg_{-2}$&$\fg_{-1}$&$\fg_0$&$r$\cr
\hline
\hline
$\fv\fle(4|3)$&$-$&$\Pi(\Lambda(3)/\Cee\cdot 1)$&$\fc(\fvect(0|3))$&$0$\cr 
\hline
$\fv\fle(5|4)$&$\Cee[-2]$&$\id\otimes\Lambda (2)$& 
$\fc(\fsl(2)\otimes\Lambda(2) \; +\hspace{-2.4mm}\subset  
T^{1/2}(\fvect(0|2))$&$1$\cr 
\hline
$\fv\fle(3|6)$&$\id_{\fsl(3)}$&$\id^*_{\fsl(3)}\otimes \id_{\fsl(2)}\otimes 
\Cee[-1]$&$\fsl(3)\oplus\fsl(2)\oplus\Cee z$&$K$\cr 
\hline
\hline
$\fv\fas(4|4)$&$-$&$\spin$&$\fas$&$0$\cr 
\hline
\hline
$\fk\fas(1|6)$&$\Cee[-2]$&$\Pi(\id)$& $\fc\fo(6)$&$0$\cr 
\hline
$\fk\fas(5|5)$ 
&$\Lambda(1)$&$\id_{\fsl(2)}\otimes\id_{\fgl(2)}\otimes\Lambda(1)$& 
$(\fsl(2)\oplus\fgl(2)\otimes\Lambda(1))\; \; +\hspace{-2.4mm}\subset \fvect(0|1)$&$1\xi$\cr 
\hline
$\fk\fas(4|4)$&$-$& $\Lambda(3)$&$\Lambda(3)\oplus\fsl(1|3)$&$3\xi$\cr 
\hline
$\fk\fas(4|3)$&$-$&$\Vol_{0}(0|3)$& 
$\fc(\fvect(0|3))$&$3\eta$\cr 
\hline
\hline
$\fm\fb(4|5)$&$\Pi(\Cee[-2])$&$\Vol (0|3)$&$\fc(\fvect(0|3))$&$0$\cr  
\hline
$\fm\fb(5|6)$&$\Lambda(2)/\Cee 1$ 
&$\id_{\fsl(2)}\otimes\Lambda(2)$ 
&$\fc(\fsl(2)\otimes\Lambda(2)\; \; +\hspace{-2.4mm}\subset  T^{1/2}(\fvect(0|2))$&$1$\cr 
\hline
$\fm\fb(3|8)$&$\id_{\fsl(3)}$&$\Pi(\id^*_{\fsl(3)}\otimes \id_{\fsl(2)}\otimes 
\Cee[-1])$&$\fsl(3)\oplus\fsl(2)\oplus\Cee z$&$K$\cr  
\hline 
\hline
$\fk\fsle(9|6)$&$\Cee[-2]$&$\Pi(T^0_{0})$&$\fsvect(0|4)_{3, 
4}$&$0$\cr    
\hline
$\fk\fsle(11|9)$&$\Pi(\id_{\fsl(1|3)})$&$\id_{\fsl(2)}\otimes\Lambda(3)$&
$\left (\fsl(2)\otimes\Lambda(3)\right) \; \; +\hspace{-2.4mm}\subset  \fsl(1|3)$&$2$\cr    
\hline
$\fk\fsle(5|10)$&$\id$&$\Pi(\Lambda^2(\id^*))$&$\fsl(5)$&$K$\cr    
\hline
$\fk\fsle(9|11)$&$\id_{\fsl(3)}^*\otimes\Lambda(1)$&$\id_{\fsl(2)}\otimes\left
(\id_{\fsl(3)}\otimes\Lambda(1)\right)$&$\fsl(2)\oplus\left
(\fsl(3)\otimes\Lambda(1) \; \; +\hspace{-2.4mm}\subset \fvect(0|1)\right)$&$CK$\cr
\hline
\end{tabular} 
$$
\normalsize

In what follows, we write just $f$ instead of $M_{f}$ or $K_{f}$, and
$I$ instead of $M_{1}$.  So, $f\cdot g$ denotes $M_{f}M_{g}$ or
$K_{f}K_{g}$, not $M_{fg}$ or $K_{fg}$.  We shorthand $D: I(\chi)\tto
I(\psi)$ to $\chi=(a, b, c, d)\tto \psi=(e, f, g, h)$. Having selected 
a basis of Cartan subalgebra, we use it for every regrading.

The negated degree of the singular vector, i.e, the degree of the
corresponding operator, is denoted by $d$).  We are sure that in some
cases there are singular vectors (operators) of degree higher than
listed ones; in such cases we write ``in degrees indicated the
singular vectors are ...''  rather than ``the following are all
possible singular vectors''.  In what follows $m_{1}$ is the
nonzero highest weight vector of the $\fg_{0}$-module $V$.

\section{Singular vectors for {\boldmath $\fg=\fv\fle(3|6)$}}

We denote the indeterminates by $x$ (even) and $\xi$ (odd); the
corresponding partial derivatives by $\partial$ and $\delta$. The 
Cartan subalgebra is spanned by\footnotesize 
$$
\renewcommand{\arraystretch}{1.4}
\begin{array}{ll}
h_1\!=\!-2\,x_{4}\!\otimes \partial_{4}\! -\! \xi_{1}\!\otimes {\delta }_1
\! -\!
\xi_{2}\!\otimes {\delta}_2\! -\!
\xi_{3}\!\otimes {\delta }_3,&\!\!
h_2=-x_{2}\!\otimes \partial_{2} - x_{3}\!\otimes \partial_{3} -
\xi_{1}\!\otimes {\delta }_1,\\ h_3\!=\!-x_{1}\otimes \partial_{1} -
x_{3}\otimes \partial_{3} -
\xi_{2}\otimes {\delta}_2,&\!\! h_4=-x_{1}\!\otimes \partial_{1} -
x_{2}\otimes
\partial_{2} - \xi_{3}\otimes {\delta }_3
\end{array}
$$\normalsize
We consider the following negative operators from $\fg_{0}$:
\tiny
$$
\renewcommand{\arraystretch}{1.4}
\begin{array}{ll}
        a_{1} = {{\partial}_4},&a_{12} = -\,{{{x_4}}^2}{{\partial}_4} 
     - {x_4} {{\xi }_1}{{\delta }_1} 
     - {x_4} {{\xi }_2}{{\delta }_2} 
     - {x_4} {{\xi }_3}{{\delta }_3} 
     + {{\xi }_1} {{\xi }_2}{{\partial}_3} 
     - {{\xi }_1} {{\xi }_3}{{\partial}_2} 
     + {{\xi }_2} {{\xi }_3}{{\partial}_1} \cr
\end{array}
$$
$$
\renewcommand{\arraystretch}{1.4}
\begin{array}{ll}
        a_{2} = -\,{x_1}{{\partial}_1} 
     - {x_2}{{\partial}_2} - {x_3}{{\partial}_3} 
     + {x_4}{{\partial}_4},  &
        a_{3} = {x_2}{{\partial}_2} + {x_3}{{\partial}_3} 
     + {{\xi }_1}{{\delta }_1}  \cr
a_{4} = -\,{x_2}{{\partial}_1} 
     + {{\xi }_1}{{\delta }_2},  &
a_{5} = -\,{x_3}{{\partial}_1} 
     + {{\xi }_1}{{\delta }_3}  \cr
a_{6} = -\,{x_1}{{\partial}_2} 
     + {{\xi }_2}{{\delta }_1},  &
a_{7} = {x_1}{{\partial}_1} + {x_3}{{\partial}_3} 
     + {{\xi }_2}{{\delta }_2}  \cr
a_{8} = -\,{x_3}{{\partial}_2} 
     + {{\xi }_2}{{\delta }_3},  &
a_{9} = -\,{x_1}{{\partial}_3} 
     + {{\xi }_3}{{\delta }_1}  \cr
a_{10} = -\,{x_2}{{\partial}_3} 
     + {{\xi }_3}{{\delta }_2},  &
a_{11} = {x_1}{{\partial}_1} + {x_2}{{\partial}_2} 
     + {{\xi }_3}{{\delta }_3}  \cr
\end{array}
$$\normalsize
and the operators from $\fg_{-}$:
\vspace*{-2mm}\tiny
$$
\renewcommand{\arraystretch}{1.4}
\begin{array}{lllll}
n_{1} = {{\delta }_1},&\; &n_{4} = {x_4}{{\delta }_3} 
     - {{\xi }_1}{{\partial}_2} + {{\xi }_2}{{\partial}_1},&\; &n_{7} = {{\partial}_1}\cr
n_{2} = {{\delta }_2},&\; &n_{5} = \-\,{x_4}{{\delta }_2} 
     - {{\xi }_1}{{\partial}_3} + {{\xi }_3}{{\partial}_1},&\; &n_{8} = {{\partial}_2}\cr
n_{3} = {{\delta }_3},&\; &
n_{6} = {x_4}{{\delta }_1} 
     - {{\xi }_2}{{\partial}_3} + {{\xi }_3}{{\partial}_2},&\; &
n_{9} = {{\partial}_3}\cr
\end{array}
$$\normalsize
The $m_{i}$ are the following elements of the irreducible
$\fg_{0}$-module $V$:
\tiny
$$
\begin{matrix}
\renewcommand{\arraystretch}{1.4}
\begin{array}{l}
m_{1}\text{ is the highest weight vector}\cr
m_{2} =  a_{12}\cdot m_{1} \cr
m_{3} =  a_{4}\cdot m_{1} \cr
m_{4} =  a_{8}\cdot m_{1} \cr
m_{5} =  a_{12}\cdot a_{12}\cdot m_{1} \cr
m_{6} =  a_{12}\cdot a_{4}\cdot m_{1} \cr
m_{7} =  a_{12}\cdot a_{8}\cdot m_{1} \cr
m_{8} =  a_{4}\cdot a_{4}\cdot m_{1} \cr
\end{array}
&
\renewcommand{\arraystretch}{1.4}
\begin{array}{l}
    m_{9} =  a_{4}\cdot a_{8}\cdot m_{1} \cr
m_{10} =  a_{8}\cdot a_{8}\cdot m_{1} \cr
    m_{11} =  a_{5}\cdot m_{1} \cr
m_{12} =  a_{12}\cdot 
     a_{12}\cdot a_{12}\cdot m_{1} \cr
m_{13} =  a_{12}\cdot 
     a_{12}\cdot a_{4}\cdot m_{1} \cr
m_{14} =  a_{12}\cdot 
     a_{12}\cdot a_{8}\cdot m_{1} \cr
m_{15} =  a_{12}\cdot 
     a_{4}\cdot a_{4}\cdot m_{1} \cr
m_{16} =  a_{12}\cdot 
     a_{4}\cdot a_{8}\cdot m_{1} \cr
\end{array}\end{matrix}
$$
\normalsize
\noindent
{\bf Theorem}\quad In $I(V)$, there are only the following singular
vectors of degree $d)$:

\noindent 1a) $(k, k, l, l) \tto(k+1,k+1, l, l)$: $n_{1}\otimes
m_{1}$;\\[1mm]
%\noindent
1b) $(k, k, k-1, k-1) \tto(k+1, k, k, k-1)$:
$n_{2}\otimes m_{1} + n_{1}\otimes m_{3}$;\\[1mm]
%\noindent
1c) $(k-3, k, k, k-1) \tto(k-2, k, k, k)$:
$n_{3}\otimes m_{1} + n_{2}\otimes m_{4} + n_{1}\otimes
m_{9}$;\\[1mm]
%\noindent
1d) $(-k, k-2, l, 1) \tto(-k-1, k-1, l, l)$, where $k\neq 0$
\vspace*{-1mm}
$$
k\, n_{6}\otimes m_{1} + n_{1}\otimes m_{2};
$$
\vspace*{-1mm}
1e) $(-k, k-2, k-1, k-1) \tto(-k-1, k-2, k, k)$, where $k\neq 0, -1$
\vspace*{-1mm}
$$
k\, (n_{5}\otimes m_{1} + n_{6}\otimes m_{3}) - n_{2}\otimes m_{2} +
n_{1}\otimes m_{6};
$$
\vspace*{-1mm}
1f) $(-k, k+1, k+1, k) \tto(-k-1, k+1, k+1, k+1)$, where $k\neq 0$
\vspace*{-1mm}
$$
  k\, (n_{4}\otimes m_{1} - n_{5}\otimes m_{4} + n_{6}\otimes m_{9}) +
  n_{3}\otimes m_{2} + n_{2}\otimes m_{7} + n_{1}\otimes m_{16};
$$
\vspace*{-1mm}
2a) $(0, -2, k, k) \tto(0, 0, k-1, k-1)$:
$\left( n_{6}{\cdot}n_{1} \right) \otimes m_{1}$;\\[1mm]

\noindent 
2b) $(0, -2, 0, 0) \tto(0, -1, 1, 0)$:
\vspace*{-1mm}
$$
  \left( n_{5}{\cdot}n_{1} - n_{6}{\cdot}n_{2} \right) \otimes m_{1} +
  \left( n_{6}{\cdot}n_{1} \right) \otimes m_{3}
$$     
2c) $(1, 0,0,-1) \tto(1, 1, 0, 0)$:
\vspace*{-1mm}
$$
-\!\left(n_{8}\! +\! n_{4}{\cdot}n_{1} \right) \otimes m_{1} + \left(
n_{3}{\cdot}n_{1} \right) \otimes m_{2} + \left( n_{9} \! +\!
n_{5}{\cdot}n_{1} \right) \otimes m_{4} + \left( n_{2}{\cdot}n_{1}
\right) \otimes m_{7} - \left( n_{6}{\cdot}n_{1} \right) \otimes m_{9}
$$
2d) $(-3, 0, 0, -1) \tto(-3, 1, 0, 0)$
\vspace*{-1mm}
\footnotesize 
$$
\renewcommand{\arraystretch}{1.4}
\begin{array}{l}
        \left( n_{8} + n_{4}{\cdot}n_{1} - 2\, \left( n_{6}{\cdot}n_{3}
\right) \right) \otimes m_{1} + \left( n_{3}{\cdot}n_{1} \right)
\otimes m_{2} \\-
\left(n_{9} +n_{5}{\cdot}n_{1} +2\, \left(
n_{6}{\cdot}n_{2} \right) \right) \otimes m_{4} + \left(
n_{2}{\cdot}n_{1} \right) \otimes m_{7} - \left( n_{6}{\cdot}n_{1}
\right) \otimes m_{9}\end{array}
$$ \normalsize
2e) $(0, k, -1, k+1) \tto(0, k, 1, k+1)$ (The dual cases were not calculated.)
\vspace*{-1mm}
$$
-k\,\left( 1 + k \right) \, \left( n_{5}{\cdot}n_{2}
\right) \otimes m_{1} -k\, \left( n_{5}{\cdot}n_{1} -
k\, n_{6}{\cdot}n_{2} \right) \otimes m_{3} + \left(
n_{6}{\cdot}n_{1} \right) \otimes m_{8}
$$

\noindent 3a) $(k-2, k, k, k) \tto(k+1, k+1, k+1, k+1)$:
$\left( n_{3}{\cdot}n_{2}{\cdot}n_{1} \right) \otimes m_{1}$;\\[1mm]
%\noindent 
3b) $ (-3, -1, -1, -1) \tto(-2, 0, 0, 0)$ (The dual cases were not calculated.)
\vspace*{-1mm}
$$
\left( n_{7}{\cdot}n_{1}\! +\! n_{8}{\cdot}n_{2} \! +\! n_{9}{\cdot}n_{3} 
\!-\!
n_{4}{\cdot}n_{2}{\cdot}n_{1} \!-\! n_{5}{\cdot}n_{3}{\cdot}n_{1} \!-\!
n_{6}{\cdot}n_{3}{\cdot}n_{2} \right) \otimes m_{1} - \left(
n_{3}{\cdot}n_{2}{\cdot} n_{1} \right) \otimes m_{2}
$$      
%\end{Theorem}

\section{Singular vectors for {\boldmath $\fg=\fv\fle(4|3)$}}

Here $\fg=\fv\fle(4|3)$, former {\cyr shch}${}_{1}$.  In
$\fg_{0}=\fc(\fvect(0|3))$ considered in the standard grading, we take
the usual basis of Cartan subalgebra, $\xi_{i}\pder{_{i}}$, and $z$ we
identify the $\fg_{0}$-module $\fg_{-1}$ with $\Pi(\Lambda
(\xi)/\Cee\cdot 1)$, by setting
$$
\renewcommand{\arraystretch}{1.4}
\begin{array}{l}
        \partial_{i}= \Pi(\xi_{i});\quad \partial_{0}=\Pi(\xi_1\xi_2\xi_3),\\
\delta_{i}=\sign (ijk)\Pi(\xi_j\xi_k)  \text{ for }(i, j, k)\in S_3. 
\end{array}
$$
We consider the following negative operators from $\fg_{0}$:
\tiny
$$
\renewcommand{\arraystretch}{1.4}
\begin{array}{ll}
    y_1 = - x_1\partial_2 + \xi_2\delta_1&
y_2 = - x_2\partial_3 + \xi_3\delta_2\cr
y_3 = - x_0\delta_3 - \xi_1\partial_2 + \xi_2\partial_1&
y_4 = x_1\partial_3 - \xi_3\delta_1\cr
y_5 = x_0\delta_2 - \xi_1\partial_3 + \xi_3\partial_1&
y_6 = x_0\delta_1 + \xi_2\partial_3 - \xi_3\partial_2 
\end{array}
$$\normalsize
and the basis of Cartan subalgebra:
\tiny
$$
\renewcommand{\arraystretch}{1.4}
\begin{array}{l}
h_0 = x_1\partial_1 + x_2\partial_2 + x_3\partial_3 + x_0\partial_0 
  + \xi_1\delta_1 + \xi_2\delta_2 + \xi_3\delta_3\\
h_1 = x_2\partial_2 + x_3\partial_3 + \xi_1\delta_1\\
h_2 = x_1\partial_1 + x_3\partial_3 + \xi_2\delta_2\\
h_3 = x_1\partial_1 + x_2\partial_2 + \xi_3\delta_3
\end{array}
$$\normalsize
The $m_{i}$ are the following elements of the irreducible
$\fg_{0}$-module $V$:
\tiny
$$
\vspace*{-2mm}
\begin{matrix}
\renewcommand{\arraystretch}{1.4}
\begin{array}{l}
m_{2} = y_{1}\,m_{1} \cr
m_{3} = y_{2}\,m_{1} \cr
m_{7} = y_{3}\,m_{1} \cr
m_{8} = y_{4}\,m_{1} \cr
m_{12} = y_{1}\,y_{3}\,m_{1}\cr
m_{15} = y_{2}\,y_{3}\,m_{1} \cr
\end{array}&
\renewcommand{\arraystretch}{1.4}\begin{array}{l}
m_{17} = y_{5}\,m_{1} \cr
m_{26} = y_{1}\,y_{5}\,m_{1} \cr
m_{31} = y_{3}\,y_{4}\,m_{1} \cr
m_{32} = (y_{4})^2\,m_{1} \cr
m_{33} = y_{6}\,m_{1}\cr
m_{46} = y_{1}\,y_{2}\,y_{5}\,m_{1}\cr
\end{array}&
\renewcommand{\arraystretch}{1.4}\begin{array}{l}
m_{54} = y_{2}\,y_{3}\,y_{4}\,m_{1}\cr
m_{57} = y_{3}\,y_{5}\,m_{1} \cr
m_{82} = y_{1}\,y_{3}\,y_{5}\,m_{1} \cr
m_{94} = y_{3}\,y_{6}\,m_{1} \cr
m_{148} = y_{3}\,y_{4}\,y_{5}\,m_{1} \cr
m_{150} = y_{5}\,y_{6}\,m_{1} \cr
\end{array}\end{matrix}
$$
\normalsize
and $m_{320} = y_{3}\,y_{5}\,y_{6}\,m_{1}$. Observe that our choice 
of ordering obscures the fact that the vectors $m_{129}$, $m_{148}$, 
and $m_{150}$ are proportional.
%\normalsize

\noindent {\bf Theorem} In $I(V)$ in degrees indicated, there are only
the following singular vectors:

\vspace*{1mm} 1a) $\lambda \tto \lambda +(-1, 1,1,1)$, where
$2\lambda_1 = \lambda_2 + \lambda_3 + \lambda_4$: \,
$\partial_0\,m_{1}$

\vspace*{1mm} 1b) $(k, l, l, k-l+1) \tto (k-1, l+1, l+1, k-l+1)$:
\, $\partial_0\,m_{7} + (k - l) \delta_3)\,m_{1}$

\vspace*{1mm} 1c) $(k, k-1, 1, k) \tto  (k-1, k, 1, k+1)$:
\vspace*{-2mm}
$$
\partial_0\,(m_{15} + (k-2) m_{17}) 
 - (k-1) \delta_2\,m_{1} 
 + \delta_3\,m_{3}
$$

\vspace*{1mm} 1d) $(k, l, k-l, k-l) \tto  (k-1, l, k-l+1, k-l+1)$:
\vspace*{-2mm}
$$
\partial_0\,(m_{26} + m_{31} + (1 - k + 2 l) m_{33}) 
 + ((k - 2 l) (1 + l)) \delta_1\,m_{1} 
 - (l+1) \delta_2\,m_{2} 
 + (l+1) \delta_3\,m_{8}
$$

\vspace*{1mm} 1e) $(1, 1, 0, 0) \tto  (0, 1, 1, 0)$:
\vspace*{-2mm}
$$
 -2 \partial_1 \,m_{2} + 2 \partial_2\,m_{1} + \partial_0\,(m_{82} + 2
 m_{94}) + 2 \delta_1\,m_{7} + 2 \delta_2\,m_{12} + \delta_3\,(m_{26}
 - m_{31} + 2 m_{33})
$$

\vspace*{1mm} 1f) $(0, 0, 0, -1) \tto  (-1, 0, 0, 0)$:
\vspace*{-2mm}\footnotesize 
$$
\renewcommand{\arraystretch}{1.4}
\begin{array}{l}
- \partial_1\,m_{8} 
 + \partial_2\,m_{3} 
 - \partial_3\,m_{1}\\ 
 + \partial_0\,(m_{129} + m_{148}) 
 + \delta_1\,(2 m_{15} - m_{17}) 
 + \delta_2\,(2 m_{31} - m_{33}) 
 + \delta_3\,(m_{46} + m_{54})
\end{array}
$$
\normalsize

\vspace*{1mm} 2a) $\lambda \tto \lambda +(-2,2,2,2)$, where $\lambda_1
= \lambda_2 + \lambda_3 + \lambda_4 + 1$: \, $\partial_0^2\,m_{1}$

\vspace*{1mm} 2b) $(k, k-2, k-2, 3) \tto (k-2, k, k, 4)$: \,
$\partial_0^2\,m_{7} + 2 \partial_0\,\delta_3\,m_{1}$

\vspace*{1mm} 2c) $(k, k-2, 2, k-1) \tto  (k-2, k, 3, k+1)$: \vspace*{-2mm}
$$
\partial_0^2\,(m_{15} + (k-4) m_{17}) 
 - 2 (k-3) \partial_0\,\delta_2\,m_{1} 
 + 2 \partial_0\,\delta_3\,m_{3}
$$

\vspace*{1mm} 2d) $(k, 1, k-1, k-1) \tto  (k-2, 2, k+1, k+1)$: \vspace*{-2mm}
$$
\partial_0^2\,(m_{26} + m_{31} + (3 - k) m_{33}) 
 + 2 (k-2) \partial_0\,\delta_1\,m_{1} 
 - 2 \partial_0\,\delta_2\,m_{2} 
 + 2 \partial_0\,\delta_3\,m_{8}
$$

\vspace*{1mm} 2e) $(2, 1, 1, 1) \tto  (0, 2, 2, 2)$: \vspace*{-2mm}
\footnotesize 
$$
\renewcommand{\arraystretch}{1.4}
\begin{array}{l}
    \partial_0^2\,m_{320} 
 + 2 \partial_1\,\partial_0\,(m_{26} + m_{31}) 
 + 2 \partial_1\,\delta_1\,m_{1} 
 - 2 \partial_1\,\delta_2\,m_{2} \\
 + 2 \partial_1\,\delta_3\,m_{8} 
 - 2 \partial_2\,\partial_0\,m_{15} 
 + 2 \partial_2\,\delta_2\,m_{1} 
 - 2 \partial_2\,\delta_3\,m_{3} \\
 + 2 \partial_3\,\partial_0\,m_{7} 
 + 2 \partial_3\,\delta_3\,m_{1} 
 - 2 \partial_0\,\delta_1\,m_{57} 
 + \partial_0\,\delta_2\,(- m_{82} + m_{94}) \\
 + \partial_0\,\delta_3\,(m_{148} + 2 m_{150}) 
 - 2 \delta_1\,\delta_2\,m_{7} 
 - 2 \delta_1\,\delta_3\,m_{17} 
 + 2 \delta_2\,\delta_3\,m_{33}
 \end{array}
$$
\normalsize

\vspace*{1mm} 3a) $\lambda \tto \lambda+(-3,3,3,3)$, where $\lambda_1
= \lambda_2 + \lambda_3 + \lambda_4 + 2$: \, $\partial_0^3\,m_{1}$

\vspace*{1mm} 3b) $(k, k-3, k-3, 4) \tto (k-3, k, k, 6)$: \,
$\partial_0^3\,m_{7} + 3 \partial_0^2\,\delta_3\,m_{1}$

\vspace*{1mm} 3c) $(k, k-3, 3, k-2) \tto  (k-3, k, 5, k+1)$: \vspace*{-2mm}
$$
\partial_0^3\,(m_{15} + (k-6) m_{17}) 
 - 3 (k-5) \partial_0^2\,\delta_2\,m_{1} 
 + 3 \partial_0^2\,\delta_3\,m_{3}
$$

\vspace*{1mm} 3d) $(k, 2, k-2, k-2) \tto  (k-3, 4, k+1, k+1)$: \vspace*{-2mm}
$$
\partial_0^3\,(m_{26} + m_{31} + (5 - k) m_{33}) 
 + 3 (k-4) \partial_0^2\,\delta_1\,m_{1} 
 - 3 \partial_0^2\,\delta_2\,m_{2} 
 + 3 \partial_0^2\,\delta_3\,m_{8}
$$

\vspace*{1mm} 3e) $(2k, k, k, k) \tto  (2k-3, k+2, k+2, k+2)$: \vspace*{-2mm}
\footnotesize 
$$
\renewcommand{\arraystretch}{1.4}
\begin{array}{l}
    \partial_0^3\,m_{320} 
 + (k+1) \partial_1\,\partial_0^2\,(m_{26} + m_{31}) 
 + k(k+1) \partial_1\,\partial_0\,\delta_1\,m_{1}\\ 
 - k(k+1) \partial_1\,\partial_0\,\delta_2\,m_{2} 
 + k(k+1) \partial_1\,\partial_0\,\delta_3\,m_{8}\\ 
 - (k+1) \partial_2\,\partial_0^2\,m_{15}
 + k(k+1) \partial_2\,\partial_0\,\delta_2\,m_{1} \\
 - k(k+1) \partial_2\,\partial_0\,\delta_3\,m_{3} 
 + (k+1) \partial_3\,\partial_0^2\,m_{7} \\
 + k(k+1) \partial_3\,\partial_0\,\delta_3\,m_{1} 
 - (k+1) \partial_0^2\,\delta_1\,m_{57} \\
 + \partial_0^2\,\delta_2\,(- m_{82} + k m_{94}) 
 + \partial_0\,\delta_3\,(m_{148} + (k+1) m_{150})\\ 
 - k(k+1) \partial_0\,\delta_1\,\delta_2\,m_{7} 
 - k(k+1) \partial_0\,\delta_1\,\delta_3\,m_{17}\\ 
 + k(k+1) \partial_0\,\delta_2\,\delta_3\,m_{33}
 - (k-1)k(k+1) \delta_1\,\delta_2\,\delta_3 \,m_{1}
\end{array}
$$
\normalsize {\bf Remarks} Cases a) have an obvious generalization to
any degree, cf.  Ref.\,\protect\citebk{Ko0}.  Some expressions can be
shortened by an appropriate ordering of the elements of the enveloping
algebra, in other words, some vectors represent zero, e.g., $m_2$,
$m_3$, $m_8$, $m_{15}$ in cases 2e) and 3e).  Being way behind the
deadline, we did not always perform such renormalization; the cases
2e) and 3e) are left as they are to entertain the reader. The 
following case --- $\fg=\fm\fb(4|5)$ --- is strikingly similar.

\section{Singular vectors for 
{\boldmath $\fg=\fm\fb(4|5)$} (after Kochetkov)}

Here $\fg=\fm\fb(4|5)$, former {\cyr shch}${}_{2}$.  Recall that in
terms of generating functions we identify the $\fg_{0}$-module
$\fg_{-2}$ with $\Pi(\Cee\cdot 1)$; we denote by ${\bf 1}\in \fm(4)$
the image of $\Pi(1)$; so $f{\bf 1}$ denotes $M_{f}M_{1}$.  We
identify $\fg_{-1}$ with $\Pi(\Lambda (\xi))$ by setting
\vspace*{-1mm}
$$
\renewcommand{\arraystretch}{1.4}
\begin{array}{l}
        x_{0}= \Pi(1),\quad x_{i}=\sign (ijk)\Pi(\xi_j\xi_k)\text{ for
        }(i, j, k)\in S_3,\\
\eta_{0}=\Pi(\xi_1\xi_2\xi_3) ,\quad \eta_{i}=\Pi(\xi_i). 
\end{array}
$$
Let $V$ be an irreducible finite dimensional $\fg_{0}$-module with
highest weight $\Lambda$, and $v_{\Lambda}$ the corresponding vector;
let $f\in I(V)$ be a nonzero singular vector.

\noindent {\bf Theorem} (\,\cite{Ko0}) \quad In $I(V)$, there are only
the following singular vectors:

\noindent 1) $\Lambda=(\frac12, \frac12, \frac12, 3)$, and $m=1$ or $3$
\footnotesize 
$$
\renewcommand{\arraystretch}{1.4}
\begin{array}{l}
\eta_{0}^m\otimes v+\mathop{\sum}\limits_{1\leq i\leq
3}\eta_{0}^{m-1}x_{i}\otimes v_{i}+
\mathop{\sum}\limits_{1\leq i\leq 3}\eta_{0}^{m-1}\eta_{i}\otimes
w_{i}+\eta_{0}^{m-1}x_{0}\otimes v_{\Lambda}
\end{array}
$$
 \normalsize
\noindent 2) $\Lambda=(0, 0, 0; 2)$ and $m=3$ 
\vspace*{-1mm}
$$
\renewcommand{\arraystretch}{1.4}
\begin{array}{l}
\eta_{0}(\mathop{\sum}\limits_{0\leq i\leq 3}\eta_{i}x_{i})\otimes
v_{\Lambda}+2 \eta_{0}{\bf 1}\otimes v_{\Lambda};
\end{array}
$$

\noindent 3) $\Lambda=(0, 0, 0; a)$:
$\partial_{0}^{2}\otimes v_{\Lambda}$\,.\\[1mm]
%\noindent 
4) $\Lambda=(0, 0, -1; 0)$ for $m=1$ or $3$ 
\vspace*{-1mm}
$$
\partial_{0}^{m}\otimes
v_{\Lambda}+\mathop{\sum}\limits_{1\leq i\leq
3}\partial_{0}^{m-1}\delta_{i}\otimes w_{i}+ \mathop{\sum}\limits_{1\leq i\leq
3}\partial_{0}^{m-1}\partial_{j}\otimes
v_{j}, \text{ where $v_{3}=v_{\Lambda}$.}
$$
%\end{Theorem}
{\bf Remark} In Ref.\,\protect\citebk{Ko0} no description of the
$v_{j}$ and $w_{j}$ is given; now we can compare the above with our
latest result:

We give the weights with respect to the
following basis of Cartan subalgebra:
$$
\renewcommand{\arraystretch}{1.4}
\begin{array}{ll}
h_1 = \tau;&
h_2 = - q_0\,\xi_0 + q_1\,\xi_1,\\
h_3 = - q_0\,\xi_0 + q_2\,\xi_2,&
h_4 = - q_0\,\xi_0 + q_3\,\xi_3.
\end{array}
$$

For the negative elements of $\fg_0$ we take
$$
\renewcommand{\arraystretch}{1.4}
\begin{array}{lll}
    y_1 = q_2\,\xi_1,&
y_2 = q_3\,\xi_2,&
y_3 = - q_0\,q_1 + \xi_2\,\xi_3,\\
y_4 = - q_3\,\xi_1,&
y_5 = - q_0\,q_2 - \xi_1\,\xi_3,&
y_6 = - q_0\,q_3 + \xi_1\,\xi_2.
\end{array}
$$
The $m_{i}$ are the following elements of the irreducible
$\fg_{0}$-module $V$:
\tiny
$$
\begin{matrix}
\renewcommand{\arraystretch}{1.4}
\begin{array}{l}
m_{2} = y_{1}\,m_{1}, \cr
m_{3} = y_{2}\,m_{1}, \cr
m_{5} = y_{1}\,y_{2}\,m_{1}\cr
m_{7} = y_{3}\,m_{1}, \cr
m_{8} = y_{4}\,m_{1}, \cr
m_{12} = y_{1}\,y_{3}\,m_{1}\cr
\end{array}&
\renewcommand{\arraystretch}{1.4}\begin{array}{l}
m_{17} = y_{5}\,m_{1} \cr
m_{24} = y_{1}\, y_{2}\,y_{3}\,m_{1} \cr
m_{30} = y_{2}\,y_{5}\,m_{1}\cr
m_{31} = y_{3}\,y_{4}\,m_{1} \cr
m_{33} = y_{6}\,m_{1}\cr
\end{array}&
\renewcommand{\arraystretch}{1.4}\begin{array}{l}
m_{57} = y_{3}\,y_{5}\,m_{1} \cr
m_{91} = y_{2}\,y_{3}\,y_{5}\,m_{1}\cr
m_{94} = y_{3}\,y_{6}\,m_{1} \cr
m_{148} = y_{3}\,y_{4}\,y_{5}\,m_{1} \cr
m_{150} = y_{5}\,y_{6}\,m_{1} \cr
\end{array}\end{matrix}
$$
\normalsize
and $m_{320} = y_{3}\,y_{5}\,y_{6}\,m_{1}$.  Observe that our choice
of ordering obscures the fact that some of the vectors either are
proportional or represent zero.  %\normalsize

\noindent
{\bf Theorem} \quad In $I(V)$, there are only the following
singular vectors:

\noindent 1a) $\lambda \tto \lambda +(-1, 1,1,1)$, where $\lambda_1 = \lambda_2
+ \lambda_3 + \lambda_4$:\;\; 
$\xi_0\,m_{1}$.

\noindent 1b) $(-k+2l-2, k, l, l)\tto (-k+2l-3, k+1, l, l)$: $ -k q_1\,m_{1} +
\xi_0\,m_{7}$

\noindent 1c) $(k, k, k+1, k+1)\tto (k-1, k, k+2, k+1)$: 
$$
-k q_1\,m_{2} - kq_2\,m_{1} + \xi_0\,m_{12}
$$

\noindent 1d) $(2k+1, k, 0, k+1)\tto  (2k, k, 1, k+1)$
$$
q_1\,m_{2} - k q_2\,m_{1} + \xi_0\,(m_{12} - (k+1) m_{17})
$$

\noindent 1e) $(k+3, k, k, k-1)\tto  (k+2, k, k, k)$
$$
 -(k-3) q_1\,(m_{5} + m_{8}) + 2 (k-3) q_2\,m_{3} - 2 (k-3) q_3\,m_{1}
 + \xi_0\,(m_{24} - 3 m_{30} + m_{31} + 5 m_{33})
$$

\noindent 1f) $(2k+1, k, k, 1)\tto  (2k, k, k, 2)$
\footnotesize
$$
 q_1\,(m_{5} + (k-1) m_{8}) - k q_2\,m_{3} + k(k-1) q_3\,m_{1} +
 \xi_0\,(m_{24} - (k+1) m_{30} + (k-1) m_{31} + (k^2+1) m_{33})
$$
\normalsize

\noindent 1g) $(3, 1, 1, 1)\tto  (2, 0, 0, 0)$
\footnotesize 
$$
\renewcommand{\arraystretch}{1.4}
\begin{array}{l}
q_0\,m_{1} -q_1\,(m_{132} + m_{148}) 
+q_2\,m_{91} -q_3\,m_{57} -\xi_0\,m_{320}\\ 
+\xi_1\,m_{7} +\xi_2\,m_{12} 
-\xi_3\,(- m_{30} + m_{31})
\end{array}
$$ \normalsize

\noindent 2a) $\lambda\tto \lambda+(-2,2,2,2)$, where $\lambda_1 =
\lambda_2 + \lambda_3 + \lambda_4 + 2$:\;\;
$\xi_0^2\,m_{1}$.

\noindent 2b) $(2k, -2, k, k)\tto (2k-2, 0, k+1, k+1)$:\;\;$
\xi_0^2\,m_{7} + 2 q_1\,\xi_0\,m_{1}$.

\noindent 2c) $(2k, k-1, -1, k)\tto  (2k-2, k, 1, k+1)$
$$
\xi_0^2\,(m_{12} - (k+1) m_{17}) 
 + 2 q_1\,\xi_0\,m_{2} 
 - 2 k q_2\,\xi_0 \,m_{1}
$$

\noindent 2d) $(2k+2, k, k, 0)\tto  (2k, k+1, k+1, 2)$
\footnotesize
$$
\renewcommand{\arraystretch}{1.4}
\begin{array}{l}
 \xi_0^2\,(- m_{30} + m_{31} + (k+1) m_{33}) 
 + q_1\,\xi_0 \,(m_{5} + m_{8}) \\
 - 2 q_2\,\xi_0 \,m_{3} 
 + 2 k q_3\,\xi_0 \,m_{1} 
\end{array}
$$ \normalsize

\noindent 2e) $(2k, 0, 0, 0)\tto  (2k-2, 0, 0, 0)$
$$
((3 - k) + q_0\,\xi_0 + q_1\,\xi_1) + q_2\,\xi_2 + q_3\,\xi_3)\,m_{1}
$$

\noindent 3a) $\lambda\tto \lambda+(-3, 3,3,3)$, where $\lambda_1 =
\lambda_2 + \lambda_3 + \lambda_4 + 4$:\; \; 
$\xi_0\,m_{1}$.

\noindent 3b) $(2k+1, -3, k, k)\tto  (2k-2, 0, k+2, k+2)$:\; \;$
\xi_0^3\,m_{7} + 3 q_1\,\xi_0^2 \,m_{1}$.

\noindent 3c) $(2k+1, k+1, -2, k)\tto  (2k-2, k+4, 1, k+2)$
$$
\xi_0^3\,(m_{12} - (k+2) m_{17}) 
 + 3 q_1\,\xi_0^2 \,m_{2} 
 - 3 (k+1) q_2\,\xi_0^2 \,m_{1}
$$

\noindent 3d) $(2k+3, k, k, -1)\tto  (2k, k+2, k+2, 2)$
\footnotesize 
$$
\renewcommand{\arraystretch}{1.4}
\begin{array}{l}
\xi_0^3\,(m_{24} - 2 m_{30} + (3 + k) m_{33}) 
 + q_1\,\xi_0^2\,((2 - k) m_{5} + (1 + k) m_{8}) \\
 - 3 q_2\,\xi_0^2\,m_{3} 
 + 3 (1 + k) q_3\,\xi_0^2 \,m_{1}
\end{array}
$$ 
\normalsize

\noindent 3e) $(2, 0, 0, 0)\tto  (-1, 1, 1, 1)$
\footnotesize 
$$
\renewcommand{\arraystretch}{1.4}
\begin{array}{l}
(q_0\,\xi_0^2 + 3 \xi_0^3 + q_1\,\xi_0\,\xi_1 + q_2\,\xi_0\,\xi_2 +
q_3\,\xi_0\,\xi_3)\,m_{1}
\end{array}
$$ \normalsize

\noindent 3f) $(k+2, k, k, k) \tto (k-1, k+1, k+1, k+1)$
\vspace*{-1mm}
\footnotesize 
$$
\renewcommand{\arraystretch}{1.4}
\begin{array}{l}
\xi_0^3\,m_{320} + (-2 + k) 
 q_0\,\xi_0^2\,m_{1} 
+ q_1\,\xi_0^2\,(m_{132} + m_{148} + (1 - k) m_{150}) \\
+ q_2\,\xi_0^2\,(- m_{91} + (-1 + k) m_{94}) 
+ (2 - k) q_3\,\xi_0^2\,m_{57} \\
- (-3 + k) (-2 + k) \xi_0\,m_{1}
+ (-2 + k) \xi_0^2\,\xi_1\,m_{7}) 
+ (-2 + k) \xi_0^2\,\xi_2\,m_{12}\\ 
+ \xi_0^2\,\xi_3\,((-2 + k) m_{30} + (2 - k) m_{31}) 
- (-2 + k) (-1 + k) kq_1\,q_2\,q_3\,m_{1} \\
+ q_1\,q_2\,\xi_0\,(m_{24} - k m_{30} + (-2 + k) m_{31} 
+ (2 - 2 k + k^2) m_{33}) \\
+ q_1\,q_3\,\xi_0\, ((-2 + k) m_{12} - (-2 + k) k m_{17}) 
-(-2 + k) (-1 + k) q_1\,\xi_0\,\xi_1\,m_{1} \\
- (-2 + k) kq_1\,\xi_0\,\xi_2\,m_{2} 
+ q_1\,\xi_0\,\xi_3\,
((2 - k) m_{5} + (-2 + k) k m_{8}) \\
+ (-2 + k) (-1 + k) (q_2\,q_3\,\xi_0\,m_{7} 
- (-2 + k) (-1 + k) q_2\,\xi_0\,\xi_2\,m_{1} \\
- (-2 + k) (-1 + k) q_2\,\xi_0\,\xi_3\,m_{3} 
- (-2 + k) (-1 + k)q_3\,\xi_0\,\xi_3\,m_{1}
\end{array}
$$ \normalsize

\section{Singular vectors for {\boldmath $\fg=\fm\fb(3|8)$}}

We give the weights with respect to the
following basis of Cartan subalgebra:
$$
\renewcommand{\arraystretch}{1.4}
\begin{array}{l}
H_{1}=\frac12\tau+\frac32q_1\xi_1-\frac12q_2\xi_2-\frac12q_3\xi_3-\frac12q_4\xi_4;\\ 
H_{2}=-q_1\xi_1+q_2\xi_2,\quad H_{3}=-q_1\xi_1+q_3\xi_3,\quad H_{3}=-q_1\xi_1+q_4\xi_4.
\end{array}
$$
The basis elements of $\fg_-$ are denoted by
$$
\renewcommand{\arraystretch}{1.4}
\begin{array}{llll}
q_0&q_1&q_2&q_3, \\ 
I&\xi_1&\xi_2&\xi_3,
\end{array}
$$
and 
$$
A=-q_0q_1+ \xi_2\xi_3, \quad B=-q_0q_2 - \xi_1\xi_3, \quad C=-q_0q_3 +
\xi_1\xi_2.
$$

\noindent {\bf Theorem} \quad In $I(V)$, there are only the following
singular vectors calculated up to dualization (though some dual
vectors are also given):

\noindent 1a) $(k, 0, l, l) \tto (k+1, 0, l-1, l-1)$:
$A\otimes m_{1}$;
   
\vspace*{1mm}    
\noindent 1b) $(k, l, -k, l+1) \tto (k-1, l, -k+1, l+1)$, where $k\neq 0, -l$
 \vspace*{-1mm}
\footnotesize 
$$
\renewcommand{\arraystretch}{1.4}
\begin{array}{l}
k\, q_1\otimes \left( q_2\xi_1{\cdot} m_{1} \right) -(k + l)\,
B\otimes \left( \xi_0{\cdot}m_{1} \right) - k\,\left( k + l \right) \,
q_2\otimes m_{1} + 
A\otimes \left(
\xi_0{\cdot}q_2\xi_1{\cdot} m_{1} \right).
\end{array}
$$

\normalsize
\noindent 1c) $(k, l, l, 2) \tto (k+1, l-1, l-1, 2)$, where $l\neq 2$
\vspace*{-1mm}
$$
\renewcommand{\arraystretch}{1.4}
\begin{array}{l}
A\otimes \left( q_2\xi_1{\cdot}q_3 \xi_2{\cdot}m_{1}
\right) - B \otimes \left( q_3\xi_2{\cdot}m_{1}\right) -(2 - l)\, C\otimes m_{1}.
\end{array}
$$
       
\noindent 1d) $(k, l, l, 1-k) \tto (k-1, l, l, 2-k)$, where $k+l\neq 1$
\vspace*{-1mm}
\footnotesize 
$$
\renewcommand{\arraystretch}{1.4}
\begin{array}{l}
-k\, q_1\otimes \left( q_3\xi_1{\cdot}m_{1} \right)
- B\otimes \left(
\xi_0{\cdot}q_3\xi_2{\cdot} m_{1} \right) -k\,\left(1 - k - l
\right) \, q_3\otimes m_{1} \\-
(1 - k - l)\, C\otimes \left(
\xi_0{\cdot}m_{1} \right) - k\, q_2\otimes \left(
q_3\xi_2{\cdot} m_{1} \right) -A\otimes \left(
\xi_0{\cdot}q_3 \xi_1{\cdot}m_{1} \right) .
\end{array}
$$
\normalsize        
\noindent 1e) $(k, l, 1, l+1) \tto (k+1,  l-1, 1, l-1)$,  where $l\neq 1$
\vspace*{-1mm}
$$
\renewcommand{\arraystretch}{1.4}
\begin{array}{l}
(1 - l)\, B\otimes m_{1} + A\otimes \left(q_2\xi_1{\cdot}m_{1} \right).
        \end{array}
$$

\noindent 1f) $(k, -k-1, l, l) \tto (k-1,  -k+1, l, l)$,  where $k\neq 0$
\vspace*{-1mm}
$$
k\, q_1\otimes m_{1} + A\otimes \left( \xi_0{\cdot}m_{1} \right).
$$

\noindent 2a) $(0, -1, l, l) \tto (0, 0, l-1, l-1)$:
$\left( q_1{\cdot}A\right) \otimes m_{1}$;
   
\vspace*{1mm}
\noindent 2b) $(0, -1, 1, 1) \tto(0, -1, 1, 0)$
\vspace*{-1mm}\footnotesize 
$$
\renewcommand{\arraystretch}{1.4}
\begin{array}{l}
\left( q_2{\cdot}A + q_1{\cdot}B\right) \otimes m_{1} + \left( q_1{\cdot}
A\right) \otimes \left(q_2\xi_1{\cdot}m_{1} \right).
\end{array}
$$
\normalsize           
\noindent 2c) $(2, 0, 0,-1) \tto (2, 0, -1, -1)$ 
\vspace*{-1mm}
\footnotesize 
$$
\renewcommand{\arraystretch}{1.4}
\begin{array}{l}
2\, \xi_2\otimes m_{1} + 2\, \xi_3\otimes \left(q_3\xi_2{\cdot} m_{1} \right)
- 2\, \left( q_1{\cdot} A\right) \otimes \left( q_2\xi_1{\cdot}
q_3\xi_2{\cdot}m_{1} \right) + 2\, \left( q_2{\cdot} A\right) \otimes
\left( q_3\xi_2{\cdot}m_{1} \right) \\-
2\, \left( q_3{\cdot} A\right) \otimes m_{1} + \left( B{\cdot}A\right)
\otimes \left( \xi_0{\cdot}q_3\xi_2{\cdot}m_{1} \right) - \left( C{\cdot}
A\right) \otimes \left(\xi_0{\cdot}m_{1} \right)
\end{array}
$$
\normalsize

\noindent 2d)  $(-4, 3, 3, 2) \tto (-4, 3, 2, 2)$
\vspace*{-1mm}
\footnotesize 
$$
\renewcommand{\arraystretch}{1.4}
\begin{array}{l}
-2\, \left( q_1{\cdot} A\right)\otimes \left( q_2\xi_1{\cdot}q_3
\xi_2{\cdot}m_{1} \right) -\,\xi_2 \otimes m_{1} -\,\xi_3 \otimes \left(
q_3\xi_2{\cdot}m_{1} \right) -\,\left( C {\cdot} A \right) \otimes \left(
\xi_0{\cdot}m_{1} \right)\\- 
3\left(\,q_1{\cdot} C \right) \otimes m_{1} +\, \left( q_3{\cdot} A\right)
\otimes m_{1} + 3\, \left( q_1{\cdot} B\right) \otimes \left(
q_3\xi_2{\cdot}m_{1} \right) - \left( q_2{\cdot} A\right)\otimes \left(
q_3\xi_2{\cdot}m_{1} \right)\\ +
\left( B{\cdot} A \right) \otimes \left( \xi_0{\cdot}q_3\xi_2{\cdot}m_{1}
\right).
\end{array}
$$
\normalsize

\noindent 2e)  $(0, k, 0, k+1) \tto (0, k-1, 1, k)$
\vspace*{-1mm}
\footnotesize 
$$
\renewcommand{\arraystretch}{1.4}
\begin{array}{l}
\left( q_1{\cdot}A \right) \otimes \left( {{\left(
q_2\xi_1 \right) }^2}{\cdot}m_{1} \right) + (1 - k)\, \left(
q_1{\cdot} B\right) \otimes \left(
q_2\xi_1{\cdot}m_{1} \right)\\ + 
(1 - k)\, \left( q_2{\cdot}A \right) \otimes \left(q_2\xi_1{\cdot}m_{1}
\right) + \left( -1 + k \right) \,k\, \left( q_2{\cdot} B\right)
\otimes m_{1}.
\end{array}
$$ 
\normalsize     

\noindent 2f)  $(0, 1, 0, 2) \tto (0, 0, 0, 2)$
\vspace*{-1mm}\footnotesize 
$$
\renewcommand{\arraystretch}{1.4}
\begin{array}{l}
\left( q_1{\cdot}A\right) \otimes \left( {{\left( q_2\xi_1 \right)
}^2}{\cdot}q_3 \xi_2{\cdot}m_{1} \right) - \left( q_1{\cdot} B\right)
\otimes \left( q_2\xi_1{\cdot} q_3\xi_2{\cdot}m_{1} \right) \\-
2\,\left( q_1{\cdot} C\right) \otimes \left( q_2\xi_1{\cdot}m_{1} \right)
- \left( q_2{\cdot} A \right) \otimes \left( q_2\xi_1{\cdot}
q_3\xi_2{\cdot}m_{1} \right) \\+
2\, \left( q_2{\cdot} B\right) \otimes \left( q_3\xi_2{\cdot}m_{1} \right)
+ 2\, \left( q_2{\cdot} C\right) \otimes m_{1} - 2\, \left( q_3{\cdot} A
\right) \otimes \left( q_2\xi_1{\cdot}m_{1} \right) + 2\, \left(
q_3{\cdot} B\right) \otimes m_{1}.
\end{array}
$$
\normalsize
\noindent 3a) $(-3, 2, 2, 2) \tto (-2, 1, 1, 1)$ 
\vspace*{-1mm}\footnotesize 
$$
\renewcommand{\arraystretch}{1.4}
\begin{array}{l}
-\,q_0\otimes m_{1} - \left( \xi_1{\cdot} A\right) \otimes m_{1} - \left(
\xi_2{\cdot} B\right) \otimes m_{1} - \left( \xi_3{\cdot}C \right) \otimes
m_{1} - \left( q_1{\cdot} C{\cdot} B\right) \otimes m_{1}\\ +
\left( q_2{\cdot} C {\cdot} A\right) \otimes m_{1} - \left( q_3{\cdot} B
{\cdot} A\right) \otimes m_{1} + \left(C{\cdot} B{\cdot} A\right) \otimes
\left( \xi_0{\cdot}m_{1} \right).
\end{array}
$$
\normalsize

\noindent 3a$^*$) $(-2, 2, 2, 2) \tto (-3, 2, 2, 2)$ 
\vspace*{-1mm}\footnotesize 
$$
\renewcommand{\arraystretch}{1.4}
\begin{array}{l}
-4\, I\otimes m_{1} - 2\, q_0\otimes \left( \xi_0{\cdot}m_{1} \right) + 2\,
\left( \xi_1{\cdot}q_1 \right) \otimes m_{1} - 2\, \left( \xi_1{\cdot}A
\right) \otimes \left( \xi_0{\cdot}m_{1} \right)\\ +
2\, \left(
\xi_2{\cdot}q_2 \right) \otimes m_{1} - 
     2\, \left( \xi_2{\cdot}B \right) \otimes \left( \xi_0{\cdot}m_{1}
     \right) + 2\, \left( \xi_3{\cdot}q_3 \right) \otimes m_{1} - 2\,
     \left( \xi_3{\cdot}C \right) \otimes \left( \xi_0{\cdot}m_{1} \right)
     \\-
     2\, \left( q_1{\cdot}C{\cdot}B \right) \otimes \left(
     \xi_0{\cdot}m_{1} \right) + 2\, \left( q_2{\cdot}q_1{\cdot}C\right)
     \otimes m_{1} + 2\, \left( q_2{\cdot}C{\cdot}A \right) \otimes \left(
     \xi_0{\cdot}m_{1} \right)\\ -
     2\, \left( q_3{\cdot}q_1{\cdot}B \right)
     \otimes m_{1} +
     2\, \left( q_3{\cdot}q_2{\cdot}A \right) \otimes m_{1} - 2\, \left(
     q_3{\cdot}B{\cdot}A \right) \otimes \left( \xi_0{\cdot}m_{1} \right)
     +C{\cdot}B {\cdot}A\otimes \left( \xi_0^2{\cdot}m_{1} \right)
\end{array}
$$    
\normalsize    

\noindent 3b) $(1-k, k, k, k) \tto (-2-k, k+1, k+1, k+1)$, where
$k\neq -1$, $0$, $1$ \vspace*{-1mm}\footnotesize
$$
\renewcommand{\arraystretch}{1.4}
\begin{array}{l}
-2\,k\,\left( 1 + k \right) I\otimes \left( \xi_0{\cdot}m_{1} \right) -
\left( 1 + k \right) q_0\otimes \left( {{\xi_0}^2}{\cdot}m_{1} \right) +
k\,\left( 1 + k \right) \left( \xi_1{\cdot}q_1 \right) \otimes \left(
\xi_0{\cdot}m_{1} \right) \\- 
   \left( 1 + k \right) \left( \xi_1{\cdot}A \right) \otimes \left(
   {{\xi_0}^2}{\cdot}m_{1} \right) + k\,\left( 1 + k \right) \left(
   \xi_2{\cdot}q_2 \right) \otimes \left( \xi_0{\cdot}m_{1} \right)\\ -
   \left( 1 + k \right) \left( \xi_2{\cdot}B \right) \otimes \left(
   {{\xi_0}^2}{\cdot}m_{1} \right) + k\,\left( 1 + k \right)
   \left(\xi_3{\cdot}q_3 \right) \otimes \left( \xi_0{\cdot}m_{1} \right)
  \\ -
   \left( 1 + k \right) \left( \xi_3{\cdot}C \right) \otimes \left(
   {{\xi_0}^2}{\cdot}m_{1} \right) - \left( 1 + k \right) \left(
   q_1{\cdot}C{\cdot}B \right) \otimes \left( {{\xi_0}^2}{\cdot}m_{1}
   \right) \\ +
   k\,\left(1 + k \right) \left( q_2{\cdot}q_1{\cdot}C \right) \otimes
   \left( \xi_0{\cdot}m_{1} \right) + \left( 1 + k \right) \left(
   q_2{\cdot}C{\cdot}A \right) \otimes \left( {{\xi_0}^2}{\cdot}m_{1}
   \right) \\ - 
   k\,\left( 1 + k \right) \left( q_3{\cdot}q_1{\cdot}B \right)
   \otimes \left( \xi_0{\cdot}m_{1} \right)\\ -
   \left( -1 + k \right) \,k\,\left( 1 + k \right) \left(
   q_3{\cdot}q_2{\cdot}q_1 \right) \otimes m_{1}\\ + 
   k\,\left( 1 + k \right) \left( q_3{\cdot}q_2{\cdot}A \right)
   \otimes \left( \xi_0{\cdot}m_{1} \right) - \left( 1 + k \right) \left(
   q_3{\cdot}B{\cdot}A \right) \otimes \left( {{\xi_0}^2}{\cdot}m_{1}
   \right) \\+
   \left( C
   {\cdot}B {\cdot}A \right) \otimes \left( {{\xi_0}^3}{\cdot}m_{1}
   \right)
\end{array}
\vspace*{-1mm}
$$      
\normalsize
\noindent 3b$^*$) $(k, 2, 2, 2) \tto (k+3, 0, 0, 0)$: 
$C{\cdot}B{\cdot}A \otimes m_{1}$.

\vspace*{1mm}

\noindent 4a) $(0, 2, 2, 1) \tto (2, 0, 0, 0)$
\vspace{-1mm}\footnotesize 
$$
\renewcommand{\arraystretch}{1.4}
\begin{array}{l}
  \left( q_0{\cdot}A \right) \otimes \left(
  q_2\xi_1{\cdot}q_3\xi_2{\cdot}m_{1} \right) - \left( q_0{\cdot}B \right)
  \otimes \left( q_3\xi_2{\cdot}m_{1} \right) + \left( q_0{\cdot}C \right)
  \otimes m_{1}\\ +
    \left( \xi_1{\cdot} B{\cdot}A \right) \otimes \left(
    q_3\xi_2{\cdot}m_{1} \right) - \left( \xi_1{\cdot}C{\cdot}A \right)
    \otimes m_{1} + \left( \xi_2{\cdot}B{\cdot}A \right) \otimes \left(
    q_2\xi_1{\cdot}q_3\xi_2{\cdot}m_{1} \right)\\ + 
    \left( \xi_2{\cdot}C{\cdot}A \right) \otimes \left(
    q_2\xi_1{\cdot}m_{1} \right) - \left( \xi_2{\cdot}C {\cdot}B \right)
    \otimes m_{1} + \left( \xi_3{\cdot}C{\cdot}A \right) \otimes \left(
    q_2\xi_1{\cdot}q_3\xi_2{\cdot}m_{1} \right) \\- 
  \left( \xi_3{\cdot}C {\cdot}B \right) \otimes \left(
  q_3\xi_2{\cdot}m_{1} \right) + \left( q_1{\cdot}C {\cdot}B{\cdot}A
  \right) \otimes \left( q_2\xi_1{\cdot}q_3\xi_2{\cdot}m_{1} \right) \\- 
  \left( q_2{\cdot}C{\cdot}B{\cdot}A \right) \otimes \left(
  q_3\xi_2{\cdot}m_{1} \right) + \left( q_3{\cdot}C{\cdot} B{\cdot}A
  \right) \otimes m_{1}
\end{array}
$$  \normalsize
%\end{Theorem}

\section{Singular vectors for {\boldmath $\fg=\fk\fs\fle(5|10)$}}

We set: $\delta _{ij}=\pder{\theta_{ij}}+\mathop{\sum}\limits_{\text{even 
permutations } (ijklm)}\theta_{kl}\partial_{m}$\,; e.g.,  
\footnotesize 
$$
\renewcommand{\arraystretch}{1.4}
\begin{array}{l}
\delta
_{12}=\pder{\theta_{12}}+\theta_{34}\partial_{5}+\theta_{45}\partial_{3}-
\theta_{35}\partial_{4},\\
\delta _{13}=\pder{\theta_{13}}+\theta_{25}\partial_{4}-\theta_{24}\partial_{5}-
\theta_{45}\partial_{2},\\
\delta _{14}=\pder{\theta_{14}}+\theta_{23}\partial_{5}+\theta_{35}\partial_{2}-
\theta_{25}\partial_{3},\; \text{etc.}\end{array}
$$  
\normalsize
The $x$-part of the elements of $\fg_{0}=\fsl(5)$ is obvious. The 
negative elements are:
$$
y_{ij} = x_i \partial_j + \sum_k \theta_{jk} \delta_{ki}\; \text{ for 
}i<j
$$
and the basis of Cartan subalgebra is $h_i = y_{ii} - y_{i+1,i+1}$.

Let us estimate the possible degree of invariant operators.  Since
$\fg_{\ev}\simeq \fsvect(5|0)$, and the grading is consistent, wee see
that the degree of the singular vector can not exceed $2\times
2+10=14$: each element from $\fg_{-1}$ can only contribute once and
the degree of singular vector of the $\fsvect(5|0)$ modules can not
exceed 2; each counted with weight 2.  In reality, the degree of
singular vectors is much lower, even with infinite dimensional fibers. 
To compute the singular vectors directly is possible on modern
computers, but hardly on a workstation; the inbuilt {\it
Mathematica}'s restrictions aggravate the problem.

Still, even simple-minded direct calculations provide us with several
first and second order operators.  The only ``known'' operator, the
exterior differential, is inhomogeneous in the consistent grading and
consists of parts of degree 1 and parts of degree 2.  To match these
parts with our operators is a problem.

The $m_{i}$ are the following elements of the irreducible
$\fg_{0}$-module $V$: 
\tiny
$$
\begin{matrix}
\renewcommand{\arraystretch}{1.4}
\begin{array}{l}
m_{1}\text{ is the highest weight vector}\cr
m_{2} = {y_{21}}\cdot\,m_{1}
  \cr
m_{3} = {y_{32}}\cdot\,m_{1}
  \cr
m_{4} = {y_{43}}\cdot\,m_{1}
  \cr
m_{5} = {y_{54}}\cdot\,m_{1}
  \cr
m_{7} = {y_{21}}\cdot\,
    {y_{32}}\cdot\,m_{1}
    \cr
m_{8} = {y_{21}}\cdot\,
    {y_{43}}\cdot\,m_{1}
    \cr
m_{9} = {y_{21}}\cdot\,
    {y_{54}}\cdot\,m_{1}
    \cr
m_{11} = {y_{32}}\cdot\,
    {y_{43}}\cdot\,m_{1}
    \cr
m_{12} = {y_{32}}\cdot\,
    {y_{54}}\cdot\,m_{1}
    \cr
m_{14} = {y_{43}}\cdot\,
    {y_{54}}\cdot\,m_{1}
    \cr
m_{16} =  -\,
     {y_{31}}\cdot\,m_{1}\cr
m_{17} = -\,
     {y_{42}}\cdot\,m_{1}\cr
\end{array}&
\renewcommand{\arraystretch}{1.4}\begin{array}{l}
m_{18} = -\,
     {y_{53}}\cdot\,m_{1}\cr
m_{24} = {y_{21}}\cdot\,
    {y_{32}}\cdot\,
     {y_{43}}\cdot\,m_{1}
     \cr
m_{25} = {y_{21}}\cdot\,
    {y_{32}}\cdot\,
     {y_{54}}\cdot\,m_{1}
     \cr
m_{30} =  -\,
     {y_{21}}\cdot\,
      {y_{42}}\cdot\,m_{1}
        \cr
m_{31} =  -\,
     {y_{21}}\cdot\,
      {y_{53}}\cdot\,m_{1}
        \cr
m_{36} = {y_{32}}\cdot\,
    {y_{43}}\cdot\,
     {y_{54}}\cdot\,m_{1}
     \cr
m_{40} =  -\,
     {y_{32}}\cdot\,
      {y_{53}}\cdot\,m_{1}
        \cr
m_{44} =  -\,
     {y_{43}}\cdot\,
      {y_{31}}\cdot\,m_{1}
        \cr
m_{48} =  -\,
     {y_{54}}\cdot\,
      {y_{31}}\cdot\,m_{1}\cr
m_{49} =  -\,
     {y_{54}}\cdot\,
      {y_{42}}\cdot\,m_{1}
        \cr
        m_{51} =  -\,
     {y_{41}}\cdot\,m_{1}
      \cr
m_{52} =  -\,
     {y_{52}}\cdot\,m_{1}
      \cr
\end{array}\end{matrix}
$$
$$
\begin{matrix}
\renewcommand{\arraystretch}{1.4}
\begin{array}{l}
m_{70} = {y_{21}}\cdot\,
    {y_{32}}\cdot\,
     {y_{43}}\cdot\,
      {y_{54}}\cdot\,m_{1}\cr
m_{73} =  -\,
     {y_{21}}\cdot\,
      {y_{32}}\cdot\,
       {y_{42}}\cdot\,m_{1}\cr
m_{74} =  -\,
     {y_{21}}\cdot\,
      {y_{32}}\cdot\,
       {y_{53}}\cdot\,m_{1}\cr
m_{83} =  -\,
     {y_{21}}\cdot\,
      {y_{54}}\cdot\,
       {y_{42}}\cdot\,m_{1}\cr
m_{86} =  -\,
     {y_{21}}\cdot\,
      {y_{52}}\cdot\,m_{1}\cr
m_{101} =  -\,
     {y_{32}}\cdot\,
      {y_{43}}\cdot\,
       {y_{53}}\cdot\,m_{1}\cr
m_{115} =  -\,
     {y_{43}}\cdot\,
      {y_{54}}\cdot\,
       {y_{31}}\cdot\,m_{1}\cr
m_{124} =  -\,
     {y_{54}}\cdot\,
      {y_{41}}\cdot\,m_{1}\cr
m_{127} = {y_{31}}\cdot\,
    {y_{42}}\cdot\,m_{1}\cr
m_{128} = {y_{31}}\cdot\,
    {y_{53}}\cdot\,m_{1}\cr
m_{130} = {y_{42}}\cdot\,
    {y_{53}}\cdot\,m_{1}\cr
m_{132} =  -\,
     {y_{51}}\cdot\,m_{1}\cr
m_{171} = {y_{21}}\cdot\,
    {y_{32}}\cdot\,
     {y_{32}}\cdot\,
      {y_{43}}\cdot\,
       {y_{54}}\cdot\,m_{1}\cr
\end{array}&
\renewcommand{\arraystretch}{1.4}\begin{array}{l}
m_{175} =  -\,
     {y_{21}}\cdot\,
      {y_{32}}\cdot\,
       {y_{32}}\cdot\,
        {y_{53}}\cdot\,m_{1}\cr
m_{181} =  -\,
     {y_{21}}\cdot\,
      {y_{32}}\cdot\,
       {y_{43}}\cdot\,
        {y_{53}}\cdot\,m_{1}      \cr
m_{184} =  -\,
     {y_{21}}\cdot\,
      {y_{32}}\cdot\,
       {y_{54}}\cdot\,
        {y_{42}}\cdot\,m_{1}        \cr
m_{187} =  -\,
     {y_{21}}\cdot\,
      {y_{32}}\cdot\,
       {y_{52}}\cdot\,m_{1}\cr
m_{241} =  -\,
     {y_{32}}\cdot\,
      {y_{43}}\cdot\,
       {y_{54}}\cdot\,
        {y_{31}}\cdot\,m_{1}        \cr
m_{250} =  -\,
     {y_{32}}\cdot\,
      {y_{54}}\cdot\,
       {y_{41}}\cdot\,m_{1}\cr
m_{254} = {y_{32}}\cdot\,
    {y_{31}}\cdot\,
     {y_{53}}\cdot\,m_{1}\cr
m_{258} =  -\,
     {y_{32}}\cdot\,
      {y_{51}}\cdot\,m_{1}\cr
m_{279} = {y_{43}}\cdot\,
    {y_{31}}\cdot\,
     {y_{53}}\cdot\,m_{1}
     \cr
m_{291} = {y_{54}}\cdot\,
    {y_{31}}\cdot\,
     {y_{42}}\cdot\,m_{1}
     \cr
m_{298} = {y_{31}}\cdot\,
    {y_{52}}\cdot\,m_{1}
    \cr
m_{301} = {y_{53}}\cdot\,
    {y_{41}}\cdot\,m_{1}
    \cr
m_{397} =  -\,
     {y_{21}}\cdot\,
      {y_{32}}\cdot\,
       {y_{32}}\cdot\,
        {y_{43}}\cdot\,
         {y_{53}}\cdot\,m_{1}
         \cr
m_{539} = {y_{32}}\cdot\,
    {y_{43}}\cdot\,
     {y_{31}}\cdot\,
      {y_{53}}\cdot\,m_{1}
      \cr
\end{array}\end{matrix}
$$
\normalsize
The Cartan subalgebra is spanned by
\tiny
$$
\renewcommand{\arraystretch}{1.4}
\begin{array}{l}
h_1={x_1}\partial_{1} - \theta _{12}\delta _{12} - \theta _{13}\delta
_{13} - \theta _{14}\delta _{14} - \theta _{15}\delta _{15} -
{x_2}\partial_{2} - \theta _{12}\delta _{12} - \theta _{23}\delta
_{23} - \theta _{24}\delta _{24} - \theta _{25}\delta _{25}\\
      h_2={x_2}\partial_{2} - \theta _{12}\delta _{12} - \theta
      _{23}\delta _{23} - \theta _{24}\delta _{24} - \theta
      _{25}\delta _{25} - {x_3}\partial_{3} - \theta _{13}\delta _{13}
      - \theta _{23}\delta _{23} - \theta _{34}\delta _{34} - \theta
      _{35}\delta _{35}\\
h_3={x_3}\partial_{3} - \theta _{13}\delta _{13} - \theta _{23}\delta
_{23} - \theta _{34}\delta _{34} - \theta _{35}\delta _{35} -
{x_4}\partial_{4} - \theta _{14}\delta _{14} - \theta _{24}\delta
_{24} - \theta _{34}\delta _{34} - \theta _{45}\delta _{45}\\
h_4={x_4}\partial_{4} - \theta _{14}\delta _{14} - \theta _{24}\delta
_{24} - \theta _{34}\delta _{34} - \theta _{45}\delta _{45} - {{\xi
}_1}\partial_{5} - \theta _{15}\delta _{15} - \theta _{25}\delta _{25}
- \theta _{35}\delta _{35} - \theta _{45}\delta _{45}
\end{array}
$$

\normalsize

\noindent
{\bf Theorem} \quad In $I(V)$ in degree $d)$, there are only the
following singular vectors (computed for degree $2$ up to dualization):

\vspace*{1mm}
\noindent 1a) $(k, l, 0, 0) \tto  (k, l+1, 0, 0)$:
$
\delta _{12}\otimes m_{1}
$;

\noindent 1a$^*$) $(0, 0, k, l)\tto (0, 0, k-1, l)$, where $k\neq 0$
and $k+l+1\neq 0$\,\footnote{~Hereafter in similar statements the reader
can check our restrictions: the coefficient of $\otimes \, m_{1}$ must
not vanish.}\footnotesize 
$$
\renewcommand{\arraystretch}{1.4}
\begin{array}{l}
2k\left( 1 + k + l \right) \delta _{45}\otimes \, m_{1} 
 -2\left( 1 + k + l \right) \delta _{35}\otimes \, m_{4}\\ +
  2\delta _{25}\otimes \left(  \left( 1 + k - l \right) \, m_{11} 
+ 2\,l \, m_{17} \right)\\ + 
  2\delta _{15}\otimes \left( 3\,
\left( -1 + k - l \right) \, m_{24} 
- 2\,\left( -1 + 2\,k - 2\,l \right) \, 
m_{44} + 2\,\left( -1 + 2\,k - l \right)
\, m_{51} \right) \\ + 
  2\delta _{34}\otimes \left( m_{14} 
      - \left( 1 + k \right) \, m_{18} \right)
    + 2\delta _{24}\otimes \left(m_{36} 
      + \left( 1 + k \right)  \, m_{40} 
      - 2\, m_{49} + 2\, m_{52} \right)\\
  -2\delta _{23}\otimes \left( m_{101} 
      + 2\, m_{130} \right)  
 \\ +
  2\delta _{14}\otimes\! \left( 3\, m_{70} 
     \! + 3\left( -1 \!+\! k \right) m_{74} 
      \!- 4\, m_{115}\! + 
          2\, m_{124}  
     \! - 2\left( -1 \!+\! 2\,k \right) m_{128} 
      \! - 2\, m_{132} \right) \\  
  -2\delta _{13}\otimes \left( 3\, m_{181} 
      -4\, m_{279} +2\, m_{301} \right)  
  + \delta _{12}\otimes \left( m_{397} - 4\, m_{539} \right)
\end{array}
$$  
\normalsize
\noindent 1b) $(k, l, -k-1, 0) \tto  (k+1, l-1, -k, 0)$,  where $l\neq 0$
$$
-l \,\delta _{13}\otimes m_{1} 
  + \delta _{12}\otimes m_{3}
$$

\noindent
1b$^*$) $(0, k, l, -k-1)\tto (0, k-1, l+1, -k)$, where $k\neq 0,
-1$ and $l\neq -1$\footnotesize 
$$
\renewcommand{\arraystretch}{1.4}
\begin{array}{l}
k(1 + k)  (1 + l)\delta _{35}\otimes  \, m_{1}  
 -(1 + k)(1 + l )\delta _{25}\otimes \, m_{3}\\+ 
          \delta _{15}\otimes 
    \left( \left( 1 + k \right) \,
          \left( 1 + k - l \right)  \, m_{7} 
      - \left( 1 + k \right) \,
           \left( k - 2\,l \right)   \, m_{16} \right) \\ +
  k(1 + l)  \delta _{34}\otimes m_{5} 
 -(1 + l)\delta _{24}\otimes \,  m_{12}\\  +
\delta _{23}\otimes 
    \left( m_{36} -\left( 1 +k \right) \, m_{40} 
      -\left( 1 + k \right) \, m_{49} 
      - k\,\left( 1 + k \right)  \, m_{52} \right)\\ +
        \delta _{14}\otimes \left( \left( 1 + k - l \right) \, 
        m_{25} -\left(k -2\,l \right) \, m_{48} \right) \\ +
  \delta _{13}\otimes \left( m_{70} - m_{74} 
      + \left( 1 + k \right) \, m_{83} + k\, m_{86} - 2\, m_{115} 
      + 2\, m_{124} + \left( 2 + k \right) \, m_{128} \right.\\
\left.      -\left(2 +k^2 \right) \, m_{132} \right)  -
  \delta _{12}\otimes \left(m_{184} - k\, m_{258} 
-2\, m_{291} -k\, m_{298} \right)
\end{array}
$$
 \normalsize
\noindent 1c) $(k, 0, l, -l-1) \tto  (k+1, 0, l-1, -l)$, where $l\neq 0$
$$
l\, \delta _{14}\otimes m_{1} 
  - \delta _{13}\otimes m_{4} + \delta _{12}\otimes m_{17}
$$

\noindent 1c$^*$) $(k, -k-1, 0, l)\tto   (k-1, -k, 0, l-1)$, where $k\neq 0$ 
and $l\neq 0$\footnotesize 
$$
\renewcommand{\arraystretch}{1.4}
\begin{array}{l}
        -k\,l \delta _{25}\otimes m_{1}
  + l\, \delta _{15}\otimes m_{2} 
  + k\, \delta _{24}\otimes m_{5} 
  - k\, \delta _{23}\otimes m_{18} 
  - \delta _{14}\otimes m_{9} \\+
  \delta _{13}\otimes m_{31} 
  + \delta _{12}\otimes \left( m_{86} 
      + \left( 1 + k \right) \, m_{132} \right)
\end{array}
$$     \normalsize

\noindent 1d) $(k, l, -k-l-2, 0) \tto (k-1, l, -k-l-1, 0)$, where
$k\neq 0$ and $k+l+1\neq 0$
\footnotesize 
$$
\renewcommand{\arraystretch}{1.4}
\begin{array}{l}
k(1 + k + l) \delta _{23}\otimes \, m_{1} 
-(1 + k + l) \delta _{13}\otimes \, m_{2}\\ +
\delta _{12}\otimes 
\left( m_{7} + \left( -1 - k \right) \, m_{16} \right) 
\end{array}
$$
\normalsize
\noindent 1d$^*$) $(0, k, l, -k-l-2) \tto (0, k-1, l, -k-l-1)$, where
$k\neq 0$ and $k+l+1\neq 0$
\footnotesize 
$$
\renewcommand{\arraystretch}{1.4}
\begin{array}{l}
k(1 + k + l) \delta _{34}\otimes \, m_{1} 
-(1 + k + l) \delta _{24}\otimes  m_{3} + \delta _{23}\otimes 
    \left( m_{11} -(1 +k) \, m_{17} \right)  
 \\ +
  \delta _{14}\otimes \left( \left( 1 \!+\! k\! -\! l \right)  m_{7} 
      + 2\,l \, m_{16} \right)  
  + \delta _{13}\otimes \left( m_{24} 
      + \left( 1 + k \right)  m_{30} - 2\, m_{44} + 2\, m_{51}
      \right) \\ -
          \delta _{12}\otimes 
    \left(m_{73} - 2\, m_{127} \right)
\end{array}
$$ 
\normalsize
\noindent 1e) $(k, 0, 0, l) \tto  (k+1, 0, 0, l-1)$, where $l\neq 0$
$$
  l\, \delta _{15}\otimes m_{1} - \delta _{14}\otimes m_{5} + \delta
  _{13}\otimes m_{18} + \delta _{12}\otimes m_{52}
$$  

\noindent 1f) $(k, -k-1, l, -l-1) \tto (k-1, -k, l-1, -l)$, where $k\neq
0$ and $l\neq 0, -1$
\footnotesize 
$$
\renewcommand{\arraystretch}{1.4}
\begin{array}{l}
( k\,l\,  -1 - l)\delta _{24}\otimes \, m_{1} 
    + k(1 + l) \delta _{23}\otimes \, m_{4} 
  + l(1 + l) \delta _{14}\otimes \, m_{2} \\
  -\left( 1 +l \right)\delta _{13}\otimes \, m_{8}
      + \delta _{12}\otimes 
    \left( m_{24} + l\, m_{30} -\left(1 + k \right) \, m_{44} 
      + \left( 1 + k \right) \,\left( 1 + l \right) \, m_{51} \right)
\end{array}
$$ 
 \normalsize
\noindent 2a) $(k, 0, 0, 1) \tto (k+1, 1, 0, 0)$
$$
\renewcommand{\arraystretch}{1.4}
\begin{array}{l}
        \delta _{15} \delta _{12}m_{1} 
- \delta _{14} \delta _{12}m_{5}
+ \delta _{13} \delta _{12}m_{18}
\end{array}
$$

\noindent 2b) $(k, -k-1, 0, 1) \tto (k-1, -k+1, 0, 0)$, where $k\neq 0$
$$
\renewcommand{\arraystretch}{1.4}
\begin{array}{l}
-k \, \delta _{25} \delta _{12}m_{1} \!+\! \delta _{15} \delta _{12}m_{2}
\!+ \! k\, \delta _{24} \delta _{12} m_{5}\! -\!
\delta _{14} \delta _{12}m_{9}\! - \! k\, \delta _{23}
\delta _{12} m_{18} \!+\! \delta _{13} \delta _{12}m_{31}
\end{array}
$$ 

%\end{Theorem}

\section{Singular vectors for {\boldmath $\fg=\fk\fs\fle(9|11)$}}

%\vspace*{-2mm}

Consider the following negative operators from $\fg_{0}$:
\tiny
$$
\begin{matrix}
\renewcommand{\arraystretch}{1.4}
\begin{array}{l}
 {y_1} = x_{2}\partial_{1} 
 - \theta _{13}\delta _{23} 
 - \theta _{14}\delta _{24} 
 - \theta _{15}\delta _{25} 
 \cr
{y_2} = x_{3}\partial_{2} 
 - \theta _{12}\delta _{13} 
 - \theta _{24}\delta _{34} 
 - \theta _{25}\delta _{35} 
 \cr
{y_3} = \delta _{12} 
 + \theta _{34}\partial_{5} 
 + \theta _{45}\partial_{3} 
 - \theta _{35}\partial_{4} 
 \cr
{y_4} = x_{5}\partial_{4} 
 - \theta _{14}\delta _{15} 
 - \theta _{24}\delta _{25} 
 - \theta _{34}\delta _{35} 
 \cr
\end{array}&\renewcommand{\arraystretch}{1.4}\begin{array}{l}
{y_5} = -x_{3}\partial_{1} 
-\theta _{12}\delta _{23} 
+\theta _{14}\delta _{34} 
+\theta _{15}\delta _{35} 
 \cr
{y_6} = -\delta _{13} 
-\theta _{25}\partial_{4} 
+\theta _{24}\partial_{5} 
+\theta _{45}\partial_{2} 
\cr
{y_7} = -\delta _{23} 
-\theta _{14}\partial_{5} 
-\theta _{45}\partial_{1} 
+\theta _{15}\partial_{4}
 \cr
\end{array}\end{matrix}
$$
\normalsize
and the operators from $\fg_{-}$:
\tiny
$$
\begin{matrix}
\renewcommand{\arraystretch}{1.4}
\begin{array}{l}
n_{1} = \partial_{4}\cr
n_{2} = \partial_{5}\cr
n_{3} = \delta _{14} 
 + \theta _{23}\partial_{5} 
 + \theta _{35}\partial_{2} 
 - \theta _{25}\partial_{3} 
 \cr
n_{4} = \delta _{15} 
 + \theta _{24}\partial_{3} 
 - \theta _{23}\partial_{4} 
 - \theta _{34}\partial_{2} 
 \cr
n_{5} = \delta _{24} 
 + \theta _{15}\partial_{3} 
 - \theta _{13}\partial_{5} 
 - \theta _{35}\partial_{1} 
 \cr
n_{6} = \delta _{25} 
 + \theta _{13}\partial_{4} 
 + \theta _{34}\partial_{1} 
 - \theta _{14}\partial_{3} 
 \cr
n_{7} = \delta _{34} 
 + \theta _{12}\partial_{5} 
 + \theta _{25}\partial_{1} 
 - \theta _{15}\partial_{2} 
 \cr
n_{8} = \delta _{35} 
 + \theta _{14}\partial_{2} 
 - \theta _{12}\partial_{4} 
 - \theta _{24}\partial_{1} 
 \cr
n_{9} = {x_4}\partial_{1} 
 + \theta _{12}\delta _{24} 
 + \theta _{13}\delta _{34} 
 - \theta _{15}\delta _{45} 
 \cr
\end{array}&\renewcommand{\arraystretch}{1.4}
\begin{array}{l}
n_{10} = {{\xi }_1}\partial_{1} 
 + \theta _{12}\delta _{25} 
 + \theta _{13}\delta _{35} 
 + \theta _{14}\delta _{45} \cr
n_{11} = {x_4}\partial_{2} 
 - \theta _{12}\delta _{14} 
 + \theta _{23}\delta _{34} 
 - \theta _{25}\delta _{45} 
 \cr
n_{12} = {{\xi }_1}\partial_{2} 
 - \theta _{12}\delta _{15} 
 + \theta _{23}\delta _{35} 
 + \theta _{24}\delta _{45} 
 \cr
n_{13} = {x_4}\partial_{3} 
 - \theta _{13}\delta _{14} 
 - \theta _{23}\delta _{24} 
 - \theta _{35}\delta _{45} 
 \cr
n_{14} = {{\xi }_1}\partial_{3} 
 - \theta _{13}\delta _{15} 
 - \theta _{23}\delta _{25} 
 + \theta _{34}\delta _{45} 
 \cr
n_{17} = \partial_{1}\cr
n_{18} = \partial_{2}\cr
n_{19} = \partial_{3}\cr
n_{20} = \delta _{45} 
 + \theta _{12}\partial_{3} 
 + \theta _{23}\partial_{1} 
 - \theta _{13}\partial_{2} 
\end{array}\end{matrix}
$$
$$
\renewcommand{\arraystretch}{1.4}
\begin{array}{l}
 n_{15} = -\, {{\xi }_1} \theta _{12}\partial_{3} + {{\xi }_1} \theta
 _{13}\partial_{2} - {{\xi }_1} \theta _{23} \partial_{1} + 2\, \theta
 _{12} \theta _{13} \theta _{23}\partial_{4} - \theta _{12} \theta
 _{13} \theta _{24} \partial_{3} + \cr 
 \theta _{12} \theta _{13}
 \theta _{34} \partial_{2} - \theta _{12} \theta _{14} \theta _{23}
 \partial_{3} - \theta _{12} \theta _{23} \theta _{34} \partial_{1} +
 \theta _{13} \theta _{14} \theta _{23} \partial_{2} + \theta _{13}
 \theta _{23} \theta _{24} \partial_{1} - {{\xi }_1}\delta _{45} +\cr
 2\,
 \theta _{12} \theta _{13} \delta _{15} + 2\, \theta _{12} \theta
 _{23} \delta _{25} - \theta _{12} \theta _{34} \delta _{45} + 2\,
 \theta _{13} \theta _{23} \delta _{35} + \theta _{13} \theta _{24}
 \delta _{45} + \theta _{14} \theta _{23} \delta _{45} \cr \cr 
 n_{16} =
 -\, {x_4} \theta _{12}\partial_{3} + {x_4} \theta _{13}\partial_{2} -
 {x_4} \theta _{23}\partial_{1} + 2\, \theta _{12} \theta _{13} \delta
 _{14} + 2\, \theta _{12} \theta _{23} \delta _{24} + \theta _{12}
 \theta _{35} \delta _{45} + \cr 2\, \theta _{13} \theta _{23} \delta
 _{34} - \theta _{13} \theta _{25} \delta _{45} - \theta _{15} \theta
 _{23} \delta _{45} -
 2\, \theta _{12} \theta _{13} \theta
 _{23}\partial_{5} + \theta _{12} \theta _{13} \theta _{25}
 \partial_{3} - \theta _{12} \theta _{13} \theta _{35} \partial_{2} +
 \theta _{12} \theta _{15} \theta _{23} \partial_{3} + \cr 
 \theta _{12}
 \theta _{23} \theta _{35} \partial_{1} - \theta _{13} \theta _{15}
 \theta _{23} \partial_{2} - \theta _{13} \theta _{23} \theta _{25}
 \partial_{1} - {x_4}\delta _{45} \cr
\end{array}
$$
\normalsize

The $m_{i}$ are the following elements of the irreducible
$\fg_{0}$-module $V$: \tiny
$$
\begin{matrix}
\renewcommand{\arraystretch}{1.4}
\begin{array}{l}
 m_{1}\text{ is the highest weight vector}\cr
m_{2} = {y_1}\cdot m_{1} \cr
m_{3} = {y_2}\cdot m_{1} \cr
m_{4} = {y_3}\cdot m_{1} \cr
m_{5} = {y_4}\cdot m_{1} \cr
m_{7} = {y_1}\cdot {y_2}\cdot m_{1} \cr
m_{9} = {y_1}\cdot {y_4}\cdot m_{1} \cr
m_{10} = {y_2}\cdot {y_2}\cdot m_{1} \cr
m_{11} = {y_2}\cdot {y_3}\cdot m_{1} \cr
m_{12} = {y_2}\cdot {y_4}\cdot m_{1} \cr
m_{13} = {y_3}\cdot {y_4}\cdot m_{1} \cr
\end{array}&
\renewcommand{\arraystretch}{1.4}
\begin{array}{l}
m_{15} = {y_5}\cdot m_{1} \cr
m_{21} = {y_1}\cdot 
 {y_2}\cdot {y_2}\cdot m_{1} \cr
m_{22} = {y_1}\cdot 
 {y_2}\cdot {y_3}\cdot m_{1} \cr
m_{27} = {y_1}\cdot {y_6}\cdot m_{1} \cr
m_{31} = {y_2}\cdot 
 {y_3}\cdot {y_4}\cdot m_{1} \cr
m_{36} = {y_3}\cdot {y_5}\cdot m_{1} \cr
m_{39} = {y_4}\cdot {y_5}\cdot m_{1} \cr
m_{41} = {y_7}\cdot m_{1} \cr 
m_{56} = {y_1}\cdot 
 {y_2}\cdot {y_3}\cdot {y_4}\cdot m_{1}\cr
m_{65} = {y_1}\cdot 
 {y_4}\cdot {y_6}\cdot m_{1} \cr
m_{82} = {y_3}\cdot 
 {y_4}\cdot {y_5}\cdot m_{1} \cr
\end{array}&\renewcommand{\arraystretch}{1.4}
\begin{array}{l}
m_{88} = {y_4}\cdot {y_7}\cdot m_{1} \cr
m_{91} = {y_1}\cdot 
 {y_1}\cdot {y_1}\cdot {y_1}\cdot {y_1}\cdot m_{1}
 \cr
m_{92} = {y_1}\cdot 
 {y_1}\cdot {y_1}\cdot {y_1}\cdot {y_2}\cdot m_{1}
 \cr
m_{93} = {y_1}\cdot 
 {y_1}\cdot {y_1}\cdot {y_1}\cdot {y_3}\cdot m_{1}
 \cr
m_{94} = {y_1}\cdot 
 {y_1}\cdot {y_1}\cdot {y_1}\cdot {y_4}\cdot m_{1}
 \cr
m_{95} = {y_1}\cdot 
 {y_1}\cdot {y_1}\cdot {y_2}\cdot {y_2}\cdot m_{1}
 \cr
m_{96} = {y_1}\cdot 
 {y_1}\cdot {y_1}\cdot {y_2}\cdot {y_3}\cdot m_{1}
 \cr
m_{97} = {y_1}\cdot 
 {y_1}\cdot {y_1}\cdot {y_2}\cdot {y_4}\cdot m_{1}
 \cr
m_{98} = {y_1}\cdot 
 {y_1}\cdot {y_1}\cdot {y_3}\cdot {y_4}\cdot m_{1}
 \cr
m_{99} = {y_1}\cdot 
 {y_1}\cdot {y_1}\cdot {y_4}\cdot {y_4}\cdot m_{1}
 \cr
m_{100} = {y_1}\cdot 
 {y_1}\cdot {y_1}\cdot {y_5}\cdot m_{1} \cr
\end{array}\end{matrix}
$$
\normalsize

\noindent
{\bf Theorem} \quad In $I(V)$ in degree $1)$ (higher degrees were 
not considered), there are only the
following singular vectors:

\vspace*{1mm}
\noindent 1a) $\lambda \tto \lambda + (0, 0, -2, 1)$:
$n_{16}\otimes m_{1}$ for ANY $\lambda$;

\vspace*{1mm}
\noindent 1b) $\lambda \tto \lambda + (0, 0, -1, -1)$:
$-\lambda_4 n_{15}\otimes m_{1} + n_{16}\otimes m_{5}$ for ANY $\lambda$;

\vspace*{1mm}
\noindent 1c) $(k, l, 1, 1) \tto (k, l+1, 0, 0)$ 
\vspace*{-1mm}
$$
 -\,n_{15}\otimes m_{4} + 2\, n_{14}\otimes m_{1} + n_{16}\otimes
 m_{13} - 2\, n_{13}\otimes m_{5}
$$

\noindent 1d) $(k, l, 2, 0) \tto (k, l+1, 0, 1)$:
$n_{16}\otimes m_{4} - 2\, n_{13}\otimes m_{1}$;

\vspace*{1mm} 
\noindent 1e) $(k, -1, 1, 1) \tto (k, 0, 1, 0)$ 
\vspace*{-1mm}
$$
 -\,n_{15}\otimes m_{11} - 2\, n_{12}\otimes m_{1} + 2\, n_{14}\otimes
 m_{3} + n_{16}\otimes m_{31} + 2\, n_{11}\otimes m_{5} - 2\,
 n_{13}\otimes m_{12}
$$ 

\noindent 1f) $(-1, k, k+2, 1) \tto (-2, k, k+2, 0)$, where $k\neq 0$
\vspace*{-1mm}
\footnotesize 
$$
\renewcommand{\arraystretch}{1.4}
\begin{array}{l}
 -\,n_{15}\otimes (m_{22} - ( 1 + k ) \, m_{27} ) - 2\,k \,
 n_{10}\otimes m_{1} + 2\,k \, n_{12}\otimes m_{2} + 2\, n_{14}\otimes
 m_{7} \\+ n_{16}\otimes (m_{56} - ( 1 + k ) \, m_{65} ) + 2\,k \,
 n_{9}\otimes m_{5} - 2\,k \, n_{11}\otimes m_{9} - 2\, n_{13}\otimes
 m_{23}\end{array}
$$ 
\normalsize 

\noindent 1g) $(k, 0, -k-1, 1) \tto (k-1, 0, -k-1, 0)$, where $k\neq 0, -1$
\vspace*{-1mm}
\footnotesize 
$$
\renewcommand{\arraystretch}{1.4}
\begin{array}{l} 
 -\,n_{15}\otimes(m_{22} -(k+2)\,m_{27}-(k+1)\,m_{36} - (k+1)^2
 \,m_{41})\\ +
 2\,k\,(k + 1) \, n_{10}\otimes m_{1} + 2\, (k + 1) \, n_{12}\otimes
 m_{2} \\+ 
 n_{14}\otimes ( -2( 3 + 2\,k)\, m_{7} - 2\,( 1 + k) \, m_{15}) 
 \\+ 
 n_{16}\otimes (m_{56} - (k+2) \, m_{65}- (k+1) \, m_{82} - (k+1)^2 \,
 m_{88}) \\- 
 2\,k\,( 1 + k)\, n_{9}\otimes m_{5} 
 - 2\,( 1 + k) \, n_{11}\otimes m_{9} 
 + 2\,(1 + k)\, n_{13}\otimes ( m_{23} + m_{39})\end{array}
$$ 
\normalsize

\noindent 1h) $(k, -k-2, 1, 1) \tto (k-1, -k-2, 1, 0)$, where $k\neq 0$
\vspace*{-1mm}
\footnotesize 
$$
\renewcommand{\arraystretch}{1.4}
\begin{array}{l}
 -\,n_{15}\otimes (m_{22} - ( 1 + k) \, m_{36} + ( 1 + k) \, m_{41}) -
 2\,k \, n_{10}\otimes m_{1} - 2\, n_{12}\otimes m_{2} \\
 + 2\,
 n_{14}\otimes ( m_{7} - ( 1 + k )\, m_{15}) + n_{16}\otimes ( m_{56}
 - ( 1 + k) \, m_{82} + ( 1 + k ) \, m_{88}) \\ + 2\,k \, n_{9}\otimes
 m_{5} + 2\, n_{11}\otimes m_{9} - 2\, n_{13}\otimes ( m_{23} - ( 1 +
 k ) \, m_{39})\end{array}
$$ 
\normalsize

%\vspace*{-2mm}
\section{Singular vectors for {\boldmath $\fg=\fk\fs\fle(11|9)$}}
%\vspace*{-2mm} 

Here we realize the elements of $\fg$, as in Ref.\,\protect\citebk{KR1}, as
divergence-free vector fields and closed 2-forms with shifted parity. 
We consider the following negative operators from $\fg_{0}$:
%\vspace*{-2mm}
\tiny
$$
\renewcommand{\arraystretch}{1.4}
\begin{array}{lll}
y_1 = x_2\partial_1&
y_2 = x_4\partial_3&
y_3 = x_5\partial_4\\
y_4 = \pi dx_1dx_2&
y_5 = - x_5\partial_3&\cr
\end{array}
$$
\normalsize
and the elements of Cartan subalgebra
\tiny
$$
\renewcommand{\arraystretch}{1.4}
\begin{array}{ll}
h_1 = x_1\partial_1 - x_2\partial_2        &
h_2 = -\frac12(x_1\partial_1 +x_2\partial_2) + x_3\,\partial_3,\\
h_3 = -\frac12(x_1\partial_1 +x_2\partial_2)+ x_4\,\partial_4,&
h_4 = -\frac12(x_1\partial_1 +x_2\partial_2) + x_5\,\partial_5\cr
\end{array}
$$
\normalsize The $m_{i}$ are the following elements of the irreducible
$\fg_{0}$-module $V$: \tiny
$$
\begin{matrix}
\renewcommand{\arraystretch}{1.4}
\begin{array}{l}
    m_{2} = y_1\,m_{1}\cr
m_{3} = y_2\,m_{1}\cr
m_{4} = y_3\,m_{1}\cr
m_{6} = y_1\,y_2\,m_{1}\cr
m_{7} = y_1\,y_3\,m_{1}\cr
m_{8} = (y_2)^2\,m_{1}\cr
m_{9} = y_2\,y_3\,m_{1}\cr
m_{10} = (y_3)^2\,m_{1}\cr
m_{11} = y_5\,m_{1}\cr
m_{15} = y_1\,(y_2)^2\,m_{1}\cr
\end{array}&\renewcommand{\arraystretch}{1.4}
\begin{array}{l}
m_{16} = y_1\,y_2\,y_3\,m_{1}\cr
m_{18} = y_1\,y_5\,m_{1}\cr
m_{20} = (y_2)^2\,y_3\,m_{1}\cr
m_{21} = y_2\,(y_3)^2\,m_{1}\cr
m_{22} = y_2\,y_5\,m_{1}\cr
m_{24} = y_3\,y_5\,m_{1}\cr
m_{25} = y_4\,m_{1}\cr
m_{34} = y_1\,(y_2)^2\,y_3\,m_{1}\cr
m_{36} = y_1\,y_2\,y_5\,m_{1}\cr
m_{39} = y_1\,y_4\,m_{1}\cr
\end{array}&\renewcommand{\arraystretch}{1.4}
\begin{array}{l}
m_{42} = (y_2)^2\,(y_3)^2\,m_{1}\cr
m_{45} = y_2\,y_3\,y_5\,m_{1}\cr
m_{49} = y_3\,y_4\,m_{1}\cr
m_{50} = (y_5)^2\,m_{1}\cr
m_{74} = y_1\,y_3\,y_4\,m_{1}\cr
m_{85} = y_2\,y_3\,y_4\,m_{1}\cr
m_{89} = (y_3)^2\,y_4\,m_{1}\cr
m_{126} = y_1\,y_2\,y_3\,y_4\,m_{1}\cr
m_{146} = y_2\,(y_3)^2\,y_4\,m_{1}\cr
m_{231} = (y_2)^2\,(y_3)^2\,y_4\,m_{1}\cr
\end{array}\end{matrix}
$$
\normalsize

\noindent {\bf Theorem} \quad In $I(V)$ in degrees $d)$, there are only
the following singular vectors:

\vspace*{1mm}
\noindent 1a) $(2k, -k, l, m) \tto (2k+1, -k+\frac32, l+1/2, m+\frac12)$:
$x_3\,\partial_2 \,m_{1}$;

\vspace*{1mm}\noindent 
1b) $(2k, l, 1-k, m) \tto (2k+1, l+\frac12, -k+\frac52, m+\frac12)$:
$$
(x_3\,\partial_2)\,m_{3} + (1 - k - l) (x_4\,\partial_2)\,m_{1};
$$

\vspace*{1mm} \noindent 1c) $(2k, l, m, 2-k) \tto (2k+1, l+\frac12,
m+\frac12, -k+\frac72)$\vspace*{-2mm}\footnotesize 
$$
(x_3\,\partial_2)\,(m_{11} + (-2 + k + m) m_{15}) 
 + (1 - k - l) (x_4\,\partial_2)\,m_{4} 
 + (-1 + k + l) (-2 + k + m) (x_5\,\partial_2)\,m_{1}
$$ 
\normalsize
\vspace*{1mm} \noindent 1d) $(2k, 3-k, 3-k, 2-k) \tto (2k+1,
\frac52-k, \frac52-k, \frac52-k)$\vspace*{-2mm}\footnotesize 
$$
(\pi dx_1 dx_3)\,m_{11} - (\pi dx_1 dx_4)\,m_{4}+
(\pi dx_1 dx_5)\,m_{1}) - (x_3 \partial_2)\,m_{85} +
(x_4 \partial_2)\,m_{49}) - (x_5 \partial_2)\,m_{25}
$$ 
\normalsize

\vspace*{1mm} \noindent 1e) $(2k, k+1, l, m) \tto (2k-1, k+\frac52,
l+\frac12, m+\frac12) $: $2 k (x_3\,\partial_1)\,m_{1} +
(x_3\,\partial_2)\,m_{2}$

\vspace*{1mm} \noindent 1f) $(2k, l, 2+k, m) \tto (2k-1, l+\frac12,
k+\frac72, m+\frac12)$\vspace*{-2mm}
$$
2 k (x_3\,\partial_1)\,m_{3}) + (x_3\,\partial_2)\,m_{6}
 + 2 k (2 + k - l) (x_4\,\partial_1)\,m_{1} 
 + (2 + k - l) (x_4\,\partial_2)\,m_{2}
$$

\vspace*{1mm} \noindent 1g) $(2k, l, m, 3+k) \tto 
(2k-1, l+\frac12, m+\frac12, k+\frac92)$\vspace*{-2mm}\footnotesize 
$$
\renewcommand{\arraystretch}{1.4}
\begin{array}{l}
    (x_3\,\partial_1)\,(2 k m_{9} - 2 k (3 + k - m) m_{11})
 + (x_3\,\partial_2)\,(m_{16} + (-3 - k + m) m_{18}) \\
 + 2 k (2 + k - l) (x_4\,\partial_1)\,m_{4}
 + (2 + k - l) (x_4\,\partial_2)\,m_{7}\\
 + 2 k (2 + k - l) (3 + k - m) (x_5\,\partial_1)\,m_{1}
 + (2 + k - l) (3 + k - m) (x_5\,\partial_2)\,m_{2}
\end{array}
$$
\normalsize

\vspace*{1mm} \noindent 1h) $(2k, 4+k, 4+k, 3+k) \tto 
(2k-1, \frac72-k, \frac72-k, \frac72-k)$\vspace*{-2mm}\footnotesize 
$$
\renewcommand{\arraystretch}{1.4}
\begin{array}{l}
    (\pi \,dx_1 dx_3)\,m_{18} -(\pi \,dx_1 dx_4)\,m_{7}
(\pi \,dx_1dx_5)\,m_{2} -2 k (\pi \,dx_2dx_3)\,m_{11}
2 k (\pi \,dx_2dx_4)\,m_{4}\\
-2 k (\pi \,dx_2dx_5)\,m_{1}
 -2 k (x_3\,\partial_1)\,m_{85} -(x_3\,\partial_2)\,m_{126} 
 +2 k (x_4\,\partial_1)\,m_{49} +(x_4\,\partial_2)\,m_{74})\\
 -2 k (x_5\,\partial_1)\,m_{25} -(x_5\,\partial_2)\,m_{39}
\end{array}
$$
 \normalsize

\vspace*{1mm} \noindent 2a) $(2k, -k-1, l, m) \tto (2k+2, -k+2, l+1,
m+1)$: $(x_3\,\partial_2)^2\,m_{1}$
 
\vspace*{1mm} \noindent 2b) $(2k, -k-1, -k+1, l) \tto (2k+2, -k+1,
-k+3, l+1)$: 
$$
(x_3\,\partial_2)^2\,m_{3} + 2
x_3\partial_2\,x_4\partial_2\,m_{1}
$$

\vspace*{1mm} \noindent 2c) $(2k, -k-1, l, -k+2) \tto 
(2k+2, -k+1, l+1, -k+4)$\vspace*{-2mm}
$$
(x_3\partial_2)^2\,(m_{9} + (-2 + k + l) m_{11}) 
 + 2 x_3\partial_2\,x_4\partial_2\,m_{4}
 - 2 (-2 + k + l) x_3\partial_2\,x_5\partial_2 \,m_{1}
$$

\vspace*{1mm} \noindent 2d) $(2k, l, -k, m) \tto (2k+2, l+1, -k+3, m+1)$\vspace*{-2mm}
$$
(x_3\partial_2)^2\,m_{8} 
 + (-1 + k + l) (k + l) (x_4\partial_2)^2\,m_{1} 
 - 2 (-1 + k + l) x_3\partial_2\,x_4\partial_2 \,m_{3}
$$

\vspace*{1mm} \noindent 2e) $(2k, l, -k, -k+2) \tto 
(2k+2, l+1, -k+2, -k+4)$\vspace*{-2mm}\footnotesize 
$$
\renewcommand{\arraystretch}{1.4}
\begin{array}{l}
    (x_3\partial_2)^2\,(m_{20} - (2 m_{22})) 
 + (-1 + k + l) (k + l) (x_4\partial_2)^2\,m_{4} \\
 - 2 (-1 + k + l) x_3\partial_2\,x_4\,\partial_2\,(m_{9} - m_{11}) 
 - 2 (-1 + k + l) x_3\,\partial_2\,x_5\partial_2 \,m_{3} \\
 + 2 (-1 + k + l) (k + l) x_4\partial_2\,x_5\partial_2 \,m_{1}
\end{array}
$$ \normalsize
In particular, 
\noindent 2ea) $l=1-k$:
$$
(x_3\partial_2)^2\,m_{22} 
 + (x_4\partial_2)^2\,m_{4} 
 - x_3\partial_2\,x_4\partial_2 \,m_{9} 
 - 2 x_3\partial_2\,x_5\partial_2 \,m_{3} 
 + 2 x_4\partial_2\,x_5\partial_2 \,m_{1}
$$

\vspace*{1mm} \noindent 2f) $(2k, l, m, -k+1) \tto 
(2k+2, l+1, m+1, -k+4)$\vspace*{-2mm}\footnotesize 
$$
\renewcommand{\arraystretch}{1.4}
\begin{array}{l}
 (x_3\partial_2)^2\,(m_{42} + 2 (-2 + k + m) m_{45} + (-2 + k + m) (-1
 + k + m) m_{50})\\
 + (-1 + k + l) (k + l) (x_4\partial_2)^2\,m_{10} +
 (-1 + k + l) (k + l) (-2 + k + m) (-1 + k + m)
 (x_5\partial_2)^2\,m_{1}\\
 - 2 (-1 + k + l)
 x_3\partial_2\,x_4\partial_2 \,(m_{21} + (-2 + k + m) m_{24})\\
 + 2 (-1+ k + l) (-2 + k + m) x_3\partial_2\,x_5\partial_2 \, (m_{9} + (-1 +
 k + m) m_{11}) - \\
 2 (-1 + k + l) (k + l) (-2 + k + m)
 x_4\partial_2\,x_5\partial_2 \,m_{4}
\end{array}
$$ \normalsize
In particular, 
\noindent 2fa) $(2k, -k+1, -k, -k+1) \tto  (2k+2, -k+2, -k+1, -k+4)$:
\footnotesize 
$$
\renewcommand{\arraystretch}{1.4}
\begin{array}{l}
    (x_3\partial_2)^2\,(m_{42} - 4 m_{50}) 
 + 2 (x_4\partial_2)^2\,m_{10} 
 + 4 (x_5\partial_2)^2\,m_{1} \\
 + 4 x_3\partial_2\,x_4\partial_2\,(- m_{21} + 2 m_{24}) 
 + 8 x_3\partial_2\,x_5\partial_2)\,(- m_{9} + m_{11}) 
 + 8 x_4\partial_2\,x_5\partial_2 \,m_{4} \end{array}
$$ \normalsize
2fb) $(2k, -k, -k, -k+1) \tto  (2k+2, -k+1, -k+1, -k+4)$:
\footnotesize 
$$
(x_3\partial_2)^2\,m_{45} 
 + (x_4\partial_2)^2\,m_{10} 
 + 2 (x_5\partial_2)^2\,m_{1} 
 - x_3\partial_2\,x_4\partial_2 \,m_{21} 
 - 4 x_3\partial_2\,x_5\partial_2 \,m_{9} 
 + 4 x_4\partial_2\,x_5\partial_2 \,m_{4}
$$

\normalsize
\vspace*{1mm} \noindent 2g) $(2k, 3-k, 3-k, 1-k) \tto 
(2k+2, 3-k, 3-k, 3-k)$\vspace*{-2mm}\footnotesize 
$$
\renewcommand{\arraystretch}{1.4}
\begin{array}{l}
    (x_3\partial_2)^2\,m_{231} 
 + 2 (x_4\partial_2)^2\,m_{89}
 + 4 (x_5\,\partial_2)^2\,m_{25}\\
 - 4 \pi dx_1dx_3\,x_3\partial_2 \,m_{45} 
 + 4 \pi dx_1dx_3\,x_4\partial_2 \,m_{24}
 - 4 \pi dx_1dx_3\,x_5\partial_2 \,m_{11}\\
 + 4 \pi dx_1dx_4\,x_3\partial_2 \, (m_{21} - m_{24}) 
 - 4 \pi dx_1dx_4\,x_4\partial_2 \,m_{10} 
 + 4 \pi dx_1dx_4\,x_5\partial_2 \,m_{4}\\
 + 4 \pi dx_1dx_5\,x_3\partial_2 \, (-3 m_{9} + 2 m_{11}) 
 + 4 \pi dx_1dx_5\,x_4\partial_2 \,m_{4}
 - 8 \pi dx_1dx_5\,x_5\partial_2 \,m_{1}\\ 
 - 2 x_3\partial_2\,x_4\partial_2 \,m_{146} 
 + 4 x_3\partial_2\,x_5\partial_2 \,m_{85} \\
 - 4 x_4\partial_2\,x_5\partial_2 \,m_{49}
\end{array}
$$
\normalsize

\vspace*{1mm} \noindent 2h) $(-2, 0, k, l) \tto (-2, 3, k+1, l+1)$: 
$(x_3\partial_2)^2\,m_{2} 
 - 2 x_3\partial_1\,x_3\partial_2 \,m_{1}$

\vspace*{1mm} \noindent 2i) $(2k, -k, l, m) \tto  
(2k, -k+2, l+2, m+1)$\vspace*{-2mm}\footnotesize 
$$
(x_3\partial_2)^2\,m_{6} 
 + 2 k x_3\partial_1\,x_3\partial_2 \,m_{3} 
 + 4 k (1 + k) x_3\partial_2\,x_4\partial_1 \, m_{1}
 + 2 (1 + k) x_3\partial_2\,x_4\partial_2 \,m_{2}
$$
\normalsize

\vspace*{1mm} \noindent 2j) $(2k, k+1, l, m) \tto 
(2k, k+3, l+2, m+1)$\vspace*{-2mm}\footnotesize 
$$
\renewcommand{\arraystretch}{1.4}
\begin{array}{l}
    (x_3\partial_2)^2\,m_{6} 
 + 2 k x_3\partial_1\,x_3\partial_2 \,m_{3} 
 - 2 k (1 + 2 k) x_3\partial_1\,x_4\partial_2 \,m_{1}\\ 
 + 2 k x_3\partial_2\,x_4\partial_1 \,m_{1} 
 - 2 k x_3\partial_2\,x_4\partial_2 \,m_{2}
 - 2 (-1 + k + l) (k + l) (-2 + k + m) x_4\partial_2\,x_5\partial_2 \,m_{4}
\end{array}
$$ \normalsize
In particular, 
\noindent 2ja) $(0, 1, 1, m) \tto (0, 3, 3, m+1)$:
$$
- x_3\partial_1\,x_4\partial_2 \,m_{1} 
 + x_3\partial_2\,x_4\partial_1 \,m_{1}
$$

\vspace*{1mm} \noindent 2k) $(2k, -k, l, k+3) \tto  
(2k, -k+2, l+1, k+4)$\vspace*{-2mm}\footnotesize 
$$
\renewcommand{\arraystretch}{1.4}
\begin{array}{l}
    (x_3\partial_2)^2\,(m_{16} + (-3 - k + l) m_{18}) 
 + 2 k x_3\partial_1\,x_3\partial_2 \, (m_{9} - (3 + k - l) m_{11}) \\
 + 4 k (1 + k) x_3\partial_2\,x_4\partial_1 \, m_{4}
 + 2 (1 + k) x_3\partial_2\,x_4\partial_2 \,m_{7} \\
 + 4 k (1 + k) (3 + k - l) x_3\partial_2\,x_5\partial_1 \,m_{1} 
 + 2 (1 + k) (3 + k - l) x_3\partial_2\,x_5\partial_2 \,m_{2} 
 \end{array}
$$\normalsize

\vspace*{1mm} \noindent 2l) $(2k, k+1, l, 2-k) \tto (2k, k+3, l+1,
4-k)$\vspace*{-2mm}\footnotesize 
$$
\renewcommand{\arraystretch}{1.4}
\begin{array}{l}
    (x_3\partial_2)^2\,(m_{16} + (-2 + k + l) m_{18}) 
 + 2 k x_3\partial_1\,x_3\partial_2 \, (m_{9} + (-2 + k + l) m_{11}) \\
 - 2 k (1 + 2 k) x_3\partial_1\,x_4\partial_2 \,m_{4} 
 + 2 k (1 + 2 k) (-2 + k + l) x_3\partial_1\,x_5\partial_2 \,m_{1} \\
 + 2 k x_3\partial_2\,x_4\partial_1 \,m_{4} 
 - 2 k x_3\partial_2\,x_4\partial_2 \,m_{7}\\
 - 2 k (-2 + k + l) x_3\partial_2\,x_5\partial_1 \,m_{1} 
 + 2 k (-2 + k + l) x_3\partial_2\,x_5\partial_2 \,m_{2}
\end{array}
$$\normalsize
In particular, 
\noindent 2la) $(0, 1, l, 2) \tto  (0, 3, l+1, 4)$:
$$
x_3\partial_1\,x_4\partial_2 \,m_{4} 
 -(l-2) x_3\partial_1\,x_5\partial_2 \, m_{1}
-x_3\partial_2\,x_4\partial_1 \, m_{4}
+(l-2) x_3\partial_2\,x_5\partial_1 \,m_{1}
$$

\vspace*{1mm} \noindent 2m) $(-2, k, 1, l) \tto  
(-2, k+1, 4, l+1)$\vspace*{-2mm}\footnotesize 
$$
\renewcommand{\arraystretch}{1.4}
\begin{array}{l}
    (x_3\partial_2)^2\,m_{15} 
 + (-2 + k) (-1 + k) (x_4\partial_2)^2\,m_{2} 
 - 2 x_3\partial_1\,x_3\partial_2 \,m_{8} 
 + 2 (-2 + k) x_3\partial_1\,x_4\partial_2 \, m_{3}\\
 + 2 (-2 + k) x_3\partial_2\,x_4\partial_1 \,m_{3} 
 - 2 (-2 + k) x_3\partial_2\,x_4\partial_2 \, m_{6}
 - 2 (-2 + k) (-1 + k) x_4\partial_1\,x_4\partial_2 \,m_{1}\end{array}
$$
 \normalsize

\vspace*{1mm} \noindent 2n) $(2k, l, -k+1, k+3) \tto  
(2k, l+1, -k+3, k+5)$\vspace*{-2mm}\footnotesize 
$$
\renewcommand{\arraystretch}{1.4}
\begin{array}{l}
 (x_3\partial_2)^2\,(m_{34} - 2 (1 + k) m_{36}) - (2 + k - l) (-1 + k
 + l) (x_4\partial_2)^2\,m_{7} + \\
 2 k x_3\partial_1\,x_3\partial_2 \,
 (m_{20} - 2 (1 + k) m_{22}) +  
 2 k (-1 + k + l)
 x_3\partial_1\,x_4\partial_2 \, (- m_{9}\\
 + 2 (1 + k) m_{11}) +
 2 k
 x_3\partial_2\,x_4\partial_1 \, ((3 + k - l) m_{9} - 2 (1 + k)
 m_{11}) + \\ 
 2 x_3\partial_2\,x_4\partial_2 \, (-(-2 + l) m_{16} + (1 +
 k) (-2 + k + l) m_{18}) + 4 k (1 + k) (2 + k - l)
 x_3\partial_2\,x_5\partial_1 \,m_{3}\\
 + 2 (1 + k) (2 + k - l)
 x_3\partial_2\,x_5\partial_2 \,m_{6} - 2 k (2 + k - l) (-1 + k + l)
 x_4\partial_1\,x_4\partial_2 \,m_{4} - \\ 
 4 k (1 + k) (2 + k - l) (-1 +
 k + l) x_4\partial_2\,x_5\partial_1 \,m_{1} - 2 (1 + k) (2 + k - l)
 (-1 + k + l) x_4\partial_2\,x_5\partial_2 \,m_{2}\end{array}
$$
\normalsize

\vspace*{1mm} \noindent 2o) $(2k, l, k+2, -k+2) \tto  
(2k, l+1, k+4, -k+4)$\vspace*{-2mm}\footnotesize 
$$
\renewcommand{\arraystretch}{1.4}
\begin{array}{l}
    (x_3\partial_2)^2\,(m_{34} + (2 k) m_{36}) 
 - (2 + k - l) (-1 + k + l) (x_4\partial_2)^2\,m_{7} 
 + 2 k x_3\partial_1\,x_3\partial_2 \, (m_{20} + 2 k m_{22}) \\ 
 + 2 k (-1 + k + l) x_3\partial_1\,x_4\partial_2 \, (- m_{9} + m_{11})
 + 2 k (1 + 2 k) (-1 + k + l) x_3\partial_1\,x_5\partial_2 \,m_{3}\\ 
 + 2 k x_3\partial_2\,x_4\partial_1 \, ((3 + k - l) m_{9} + (1 + 5 k +
 2 k^2 - l - 2 k l) m_{11})\\
 + 2 x_3\partial_2\,x_4\partial_2 \, (- (-2
 + l) m_{16} + k (3 + k - l) m_{18})  
 - 2 k (-1 + k + l)
 x_3\partial_2\,x_5\partial_1 \,m_{3}\\ + 2 k (-1 + k + l)
 x_3\partial_2\,x_5\partial_2 \,m_{6}
 - 2 k (2 + k - l) (-1 + k + l)
 x_4\partial_1\,x_4\partial_2 \,m_{4}\\
 + 2 k (1 + 2 k) (2 + k - l) (-1
 + k + l) x_4\partial_1\,x_5\partial_2 \,m_{1} - 2 k (2 + k - l) (-1 +
 k + l) x_4\partial_2\,x_5\partial_1 \,m_{1} \\
 + 2 k (2 + k - l) (-1 + k + l) x_4\partial_2\,x_5\partial_2 \,m_{2}
\end{array}
$$
\normalsize

\section{Singular vectors for {\boldmath $\fg=\fk\fa\fs$} and
{\boldmath $\fg=\fk(1|n)$}}

The coordinates of the weights are given with respect to the 
following basis of $\fg_{0}$:
$$
(K_{t}, K_{\xi_{1}\eta_{1}}, \ldots, K_{\xi_{s}\eta_{s}}), \quad
s=[\frac{n}{2}].
$$

\noindent
{\bf Theorem}\quad In $I(V)$, there are only the following singular for
$\fk(1|3)$ vectors:

\noindent 1a) $(k,  -k) \tto (k-1,  -k+1)$:
$\xi _{1}\otimes m_{1}$;

\noindent 1a$^*$) $(k+1, k) \tto (k, k-1)$
$$
\xi _{1}\otimes \left( {{\left( \eta _{1} \theta _{1}
\right)}^2}{\cdot}m_{1} \right) - k\, \left( -1 + 2\, k \right) \, \eta
_{1}\otimes m_{1} + (1 - 2\, k)\, \theta _{1}\otimes \left( \eta
_{1}\theta _{1}{\cdot}m_{1} \right)
$$

\noindent 1b) $(1, -1) \tto (0, -1)$: $\xi _{1}\otimes \left( \eta _{1} \theta
_{1}{\cdot}m_{1} \right) + \theta _{1}\otimes m_{1}$;

\noindent 2a) $\frac12(3,  1) \tto\frac12(-1,  1)$
$$
I\otimes m_{1} -2\,  \left( \xi _{1} {\cdot}\theta _{1} \right)
\otimes \left( \eta _{1}\theta _{1}{\cdot}m_{1} \right) + \left(
\eta _{1} {\cdot}\xi _{1} \right) \otimes m_{1}
$$      
%\end{Theorem} 

\noindent
{\bf Theorem} \quad In $I(V)$, there are only the following singular for
$\fk(1|4)$ vectors:

\noindent 1a) $(k, -k, 0) \tto (k-1, -k+1, 0)$: $ \xi _{1}\otimes m_{1}$;

\vspace*{1mm}
\noindent 1a$^*$) $(k+1, -1, k) \tto (k, -1, k-1)$,  where $k\neq1$
\vspace*{-1mm}
$$
\xi _{1}\otimes \left( \eta _{1} \eta _{2}{\cdot}m_{1} \right) +
(1 - k)\,  \eta _{2}\otimes m_{1}
$$  
  
\noindent 1b) $(k, -1, 1-k) \tto (k-1, -1, 2-k)$,   where
$k\neq 2$
\vspace*{-1mm}
$$
(-2 + k)\,  \xi _{2}\otimes m_{1} + \xi _{1}\otimes \left( \xi
_{2}\eta _{1}{\cdot}m_{1} \right)
$$

\noindent 1c) $(k+2, k, 0) \tto (k+1, k-1, 0)$,  where $k\neq 0$
\vspace*{-1mm}
$$
\xi _{1}\otimes \left( \eta _{1} \eta _{2}{\cdot}m_{1} \right) -
k\,  \eta _{2}\otimes m_{1}
$$

\noindent 2a) $(1, -1, -1) \tto (-1, 0, 0)$:
$\left( \xi _{2}{\cdot}\xi _{1} \right) \otimes m_{1}$;

\vspace*{1mm}
\noindent 2b) $(2, -1, 1) \tto (0,  0,  0)$:
$\left( \eta _{2}{\cdot}\xi _{1} \right) \otimes m_{1}$;

\vspace*{1mm}
\noindent 2c) $(2, -1, -1) \tto (0, -1, -1)$
\vspace*{-1mm}
$$
I\otimes m_{1} -\left( \xi _{2} {\cdot}\xi _{1} \right) \otimes
\left( \eta _{1}\eta _{2}{\cdot}m_{1} \right) + \left( \eta
_{1}{\cdot}\xi _{1} \right) \otimes m_{1} + \left( \eta _{2}{\cdot}
\xi _{2} \right) \otimes m_{1}
$$      
%\end{Theorem} 

\noindent
{\bf Theorem} \quad In $I(V)$, there are only the following singular for
$\fk(1|6)$ and $\fk\fa\fs$ vectors:

\noindent 1a)  $\lambda =(k, -k, l, l)\tto \lambda +(-1, 1, 0, 0)$;
$\fk\fa\fs$ and $\fk(1|6)$: $\xi_{1}\otimes m_{1}$ (for $\fk(1|6)$ only if
$l=0$; for $\fk\fa\fs$ without restrictions);

\vspace*{1mm}
\noindent 1a$^*$) $\lambda =(k, l, 1-k, l+1)\tto \lambda +(-1, 0,  1,  0)$
\vspace*{-1mm}
$$
\renewcommand{\arraystretch}{1.4} 
\begin{array}{ll}
\fk\fa\fs:&(-1 + k + l)\, \xi _{2}\otimes m_{1} +
\xi _{1}\otimes \left( \xi _{2}\eta _{1}{\cdot}m_{1} \right)\\
\fk(1|6): &\text{ the above for $l=-1$}
\end{array}
$$

\noindent 1b) $\lambda =(k,  l,  l,  2-k)\tto \lambda +(-1, 0, 0, 1)$, where $l+k \neq  2$ 
$$
\renewcommand{\arraystretch}{1.4} 
\begin{array}{ll}
\fk\fa\fs: & \xi _{1}\otimes \left( \xi _{2}\eta_{1}{\cdot}\xi_{3}\eta
_{2}{\cdot}m_{1} \right) + \xi_{2}\otimes \left( \xi_{3}\eta _{2}{\cdot}m_{1}
\right) + (-2 + k + l)\, \xi_{3}\otimes m_{1}\\
\fk(1|6): &\text{ the above for $l=-1$}
\end{array}
$$       

\noindent 1c)  $\lambda =(k,  l,  -l-2,  k-2)\tto\lambda + (-1, 0, 0, -1)$, where $l+k\neq 1$ 
and $k-l\neq 4$ \vspace*{-1mm}
$$
\renewcommand{\arraystretch}{1.4} 
\begin{array}{ll}
\fk\fa\fs: &\xi _{1}\otimes
\left( \xi _{2}\eta _{1}{\cdot}\eta _{2}\eta_{3}{\cdot} m_{1} \right) + (-4 + k - l)\, \xi _{1}\otimes \left( \eta
_{1}\eta _{3}{\cdot}m_{1} \right)\\ 
&+(-1 + k+ l)\,  \xi _{2}\otimes \left( \eta _{2}\eta _{3}{\cdot}m_{1}
\right) - ( -4 + k - l ) \,  ( -1 + k + l) \,  \eta _{3}\otimes m_{1}\\
\fk(1|6): &\text{ the above for $l=-1$}
\end{array}
$$

The singular vectors of degree $2$ for $\fk\fa\fs$ and $\fk(1|6)$ are
the same:
\vspace*{1mm}

\noindent 2a) $(3, -1, -1, -1) \tto (1,  -1, -1, -1)$\vspace*{-1mm}
\footnotesize 
$$
\renewcommand{\arraystretch}{1.4}
\begin{array}{l}
-2\,  I\otimes m_{1} + \left( \xi _{2}{\cdot}\xi _{1} \right) \otimes
\left( \xi _{3}\eta _{2}{\cdot} \eta _{1}\eta_{3}{\cdot}m_{1} \right) + 
\left( \xi _{3}{\cdot}\xi _{1} \right)
\otimes \left( \eta _{1}\eta _{3}{\cdot} m_{1} \right)\\ +
\left(\xi _{3}{\cdot}\xi _{2} \right) \otimes \left( \eta_{2}\eta _{3}{\cdot}m_{1} \right) - 
\left( \eta _{1}{\cdot}\xi_{1} \right) \otimes m_{1} - \left( \eta _{2}{\cdot}\xi _{1} \right)
\otimes \left( \xi _{2}\eta _{1}{\cdot}m_{1} \right)\\ -
\left(\eta _{2}{\cdot}\xi _{2} \right) \otimes m_{1} + \left( \eta _{3}{\cdot}
\xi _{1} \right) \otimes \left( \xi _{2}\eta _{1}{\cdot}
\xi _{3}\eta _{2}{\cdot}m_{1} \right)\\ - 
\left( \eta_{3}{\cdot}\xi _{1} \right) \otimes \left( \xi _{3}\eta_{1}{\cdot}m_{1} \right) - 
\left( \eta _{3}{\cdot}\xi _{3} \right)\otimes m_{1}
\end{array}
$$
\normalsize

\noindent 2a$^*$)  $(3, -1, 0, 0) \tto (1, -1, 0, 0)$\vspace*{-1mm}
\footnotesize 
$$
\renewcommand{\arraystretch}{1.4}
\begin{array}{l}
        -\, I\otimes m_{1} + \left( \xi _{2}{\cdot}\xi _{1} \right) \otimes
\left( \xi _{3}\eta _{2}{\cdot} \eta _{1}\eta_{3}{\cdot}m_{1} \right) + 
\left( \xi _{3}{\cdot}\xi _{1} \right)
\otimes \left( \eta _{1}\eta _{3}{\cdot} m_{1} \right) \\+
\left(\xi _{3}{\cdot}\xi _{2} \right) \otimes 
\left(\eta_{2}\eta _{3}{\cdot}m_{1} \right) - \left( \eta _{1}{\cdot}\xi_{1} \right) \otimes m_{1} - 
\left(\eta _{2}{\cdot}\xi _{1} \right)\otimes \left( \xi _{2}\eta _{1}{\cdot}m_{1}
\right)\\ +
\left(\eta _{3}{\cdot}\xi _{1} \right) \otimes \left( \xi _{2}
\eta _{1}{\cdot}\xi _{3}\eta _{2}{\cdot}m_{1} \right) - 
\left(\eta _{3}{\cdot}\xi _{1} \right) \otimes\left( \xi_{3}\eta _{1}{\cdot}m_{1} \right) - 
\left( \eta _{3}{\cdot}\xi_{2} \right) \otimes \left(\xi _{3}\eta _{2}{\cdot}m_{1} \right)
\end{array}
$$
\normalsize
\noindent 2b) $(3, -1, -1, 1) \tto (1, -1, -1, 1)$\vspace*{-1mm}
\footnotesize 
$$
\renewcommand{\arraystretch}{1.4}
\begin{array}{l}
-\, I\otimes m_{1} + \left( \xi _{2}{\cdot}\xi _{1} \right) \otimes \left(
\eta _{2}\eta _{3}{\cdot} \xi _{3}\eta_{1}{\cdot}m_{1} \right) - \left(
\xi _{3}{\cdot}\xi _{1} \right) \otimes \left( \xi _{2}\eta
_{1}{\cdot} \eta _{2}\eta_{3}{\cdot}m_{1} \right)\\ + 
\left( \xi _{3}{\cdot}\xi _{1} \right) \otimes \left( \eta _{1}\eta
_{3}{\cdot} m_{1} \right) - \left(\eta _{1}{\cdot}\xi_{1} \right) \otimes
m_{1} - \left( \eta_{2}{\cdot}\xi _{1} \right) \otimes \left( \xi
_{2}\eta_{1}{\cdot}m_{1} \right)\\ - 
\left( \eta _{2} {\cdot}\xi _{2} \right) \otimes m_{1} - \left( \eta
_{3}{\cdot}\xi _{1} \right) \otimes \left( \xi _{3}\eta_{1}{\cdot}m_{1}
\right)\\ - 
\left( \eta_{3}{\cdot} \xi _{2} \right) \otimes \left( \xi _{3}\eta
_{2}{\cdot}m_{1} \right) + \left( \eta _{3}{\cdot}\xi _{3} \right) \otimes
m_{1}
\end{array}
$$
\normalsize

\section{Singular vectors for {\boldmath $\fg=\fk\fas(1|6; 3\xi)$}}
Set \tiny
$$
\begin{matrix}
\renewcommand{\arraystretch}{1.4}
\begin{array}{l}
 m_{1}\text{ is the highest weight vector}\cr
m_{2} = \xi_1\eta_2\,m_{1}\cr
m_{3} = \xi_3\eta_1\,m_{1}\cr
m_{4} = \eta_3\,m_{1}\cr
m_{6} = \xi_1\eta_2\,\xi_3\eta_1\,m_{1}\cr
m_{7} = \xi_1\eta_2\,\eta_3\,m_{1}\cr
m_{9} = \xi_3\eta_1\,\eta_3\,m_{1}\cr
m_{10} = \xi_3\eta_2\,m_{1}\cr
m_{11} = \eta_1\,m_{1}\cr
\end{array}&\renewcommand{\arraystretch}{1.4}\begin{array}{l}
m_{15} = \xi_1\eta_2\,(\xi_3\eta_1)^2\,m_{1}\cr
m_{16} = \xi_1\eta_2\,\xi_3\eta_1\,\eta_3\,m_{1}\cr
m_{18} = \xi_1\eta_2\,\eta_1\,m_{1}\cr
m_{23} = \eta_3\,\xi_3\eta_2\,m_{1}\cr
m_{25} = - (\eta_2\,m_{1})\cr
m_{34} = \xi_1\eta_2\,(\xi_3\eta_1)^2\,\eta_3\,m_{1}\cr
m_{36} = \xi_1\eta_2\,\xi_3\eta_1\,\eta_1\,m_{1}\cr
m_{46} = - (\xi_3\eta_1\,\eta_2\,m_{1})\cr
m_{49} = \xi_3\eta_2\,\eta_1\,m_{1}\cr
\end{array}\end{matrix}
$$\normalsize {\bf Theorem} \quad In $I(V)$ in  degrees
$d)$, there are only the following singular vectors:

\vspace*{1mm}
\noindent 1a) $\lambda \tto \lambda+(1,1,1,1)$ for ANY $\lambda$:
$\xi_1\xi_2\xi_3\,m_{1}$;

\vspace*{1mm}
\noindent 1b) $%[0,1,1,0]
(k, l-1, -l-1, k) \tto (k, l, -l, k)$:
$\xi_1\xi_2\,m_{1} + \xi_1\xi_2\xi_3\,m_{4}$;

\vspace*{1mm}
\noindent 1c) $%[0,1,0,1]
(k, l, k-1, -l-1) \tto (k, l, k, -l)$, where $k+l\neq 0$
\vspace*{-2mm}
$$
\xi_1\xi_2\,m_{3} - (k + l) \xi_2\xi_3\,m_{1} 
 + \xi_1\xi_2\xi_3\,(m_{9} - (k + l + 1) m_{11})
$$

\noindent 1d) $%[0,0,1,1] 
(k+1, k-1, l, -l) \tto (k+1, k-1, l+1, -l+1)$, where
$k+l\neq 0$, $k\neq l+1$
\vspace*{-2mm}
\footnotesize 
$$
\renewcommand{\arraystretch}{1.4}
\begin{array}{l}
 \xi_1\xi_2\,(m_{6} + (-k + l) m_{10}) 
 - (-1\! +\! k\! -\! l) (k\! +\! l) \xi_1\xi_3\,m_{1} 
 \\- 
 (k\! +\! l) \xi_2\xi_3\,m_{2} 
 + \xi_1\xi_2\xi_3\,(m_{16} - (k\! +\! l \!+ \!1) m_{18} + (-k \!+\! l)
m_{23} - 
 (k\! -\! l) (k \!+ \! l) m_{25})
\end{array}
$$
\normalsize

\noindent 1e) $(k+1, k-1, 0, 0) \tto (k+1, k-1, 1, 1)$, where $k\neq
1$ (new for $k=0$ only) \vspace*{-2mm} \footnotesize
$$
\renewcommand{\arraystretch}{1.4}
\begin{array}{l}
 \xi_1\xi_2\,(m_{6} -2m_{10})
 + 2(1 -k) \xi_1\xi_3\,m_{1} 
\\- 
 2\xi_2\xi_3\,m_{2}
-2\xi_1\xi_2\xi_3\,(m_{18} + m_{23} + k m_{25})\end{array}
$$
\normalsize

\noindent 2a) $%[1,2,2,1]
\lambda \tto \lambda+(1,2,2,1)$, where $\lambda_4=-2-\lambda_3$:
$\xi_1\xi_2\cdot\xi_1\xi_2\xi_3\,m_{1}$;

\vspace*{1mm}
\noindent 2b) $%[1,2,1,2] 
\lambda \tto \lambda+(1,2,1,2)$, where
$\lambda_4=-1-\lambda_2$, $\lambda_2+\lambda_3\neq -1$
\vspace*{-2mm}
$$
\xi_1\xi_2\cdot\xi_1\xi_2\xi_3\,m_{3} 
 + (-1 - \lambda_2 - \lambda_3) (\xi_2\xi_3\cdot\xi_1\xi_2\xi_3\,m_{1})
$$

\noindent 2c) $%[1,1,2,2] 
\lambda \tto \lambda+(1,1,2,2)$, where $\lambda_4=-\lambda_3$,
$\lambda_2 \neq \lambda_3, \lambda_2+\lambda_3\neq -1$ \vspace*{-2mm}
\footnotesize 
$$
\renewcommand{\arraystretch}{1.4}
\begin{array}{l}
 (\xi_1\xi_2\cdot\xi_1\xi_2\xi_3)\, (m_{6} + (-1 - \lambda_2 +
\lambda_3) m_{10})  -(\lambda_2 - \lambda_3) (1 + \lambda_2 +
\lambda_3)\\
\times (\xi_1\xi_3\cdot\xi_1\xi_2\xi_3)\,m_{1}
 +
 (-1 - \lambda_2 - \lambda_3)(\xi_2\xi_3\cdot\xi_1\xi_2\xi_3)\,m_{2}
\end{array}
$$
\normalsize

\noindent 2d) $(k, l, 0, 0) \tto (k+1, l+1, 2, 2)$, where $l\neq 0$ 
(new for $l=-1$ only)
\vspace*{-2mm}
$$
(\xi_1\xi_2\cdot\xi_1\xi_2\xi_3)\,m_{10} 
 + l (\xi_1\xi_3\cdot\xi_1\xi_2\xi_3)\,m_{1}
 + (\xi_2\xi_3\cdot\xi_1\xi_2\xi_3)\,m_{2}
$$

\noindent 2e) $%[0,2,1,1]
(k, -k, k-3, k-1) \tto (k, -k+2, k-2, k)$
\vspace*{-2mm}
$$
(\xi_1\xi_2)^2 \,m_{3}
 - 2 \xi_2\cdot\xi_1\xi_2\xi_3 \,m_{1}
 + 2 \xi_1\xi_2\cdot\xi_2\xi_3 \,m_{1}
 + \xi_1\xi_2\cdot\xi_1\xi_2\xi_3 \,(m_{9} - m_{11}) 
 + 2 \xi_2\xi_3\cdot\xi_1\xi_2\xi_3 \,m_{4} 
$$

\noindent 2f) $(k,-k,l,l) \tto (k, -k+2, l+1, l+1)$, where $k\neq l$
\vspace*{-2mm}
$$
\renewcommand{\arraystretch}{1.4}
\begin{array}{l}
(\xi_1\xi_2)^2 \,m_{3}-4(k - l) \xi_2\cdot\xi_1\xi_2\xi_3 \,m_{1}
-4\xi_1\xi_2\cdot\xi_1\xi_2\xi_3 \,m_{11} 
+4\xi_2\xi_3\cdot\xi_1\xi_2\xi_3 \,m_{4}\end{array}
$$

\noindent 2g) $%[0,1,2,1]
(1 + k, l, - k, 1 + l) \tto (1+k, l+1, 2-k, 2+l)$
\vspace*{-2mm}\footnotesize  
$$
\renewcommand{\arraystretch}{1.4}
\begin{array}{l}
(k - l) (1 + k - l) (k + l) \xi_1\cdot\xi_1\xi_2\xi_3 \,m_{1} + (k - l)
(1 + k - l) \xi_2\cdot\xi_1\xi_2\xi_3 \,m_{2} \\+
 \xi_1\xi_2\cdot\xi_1\xi_2\xi_3 \, (m_{16} + (k - l) m_{18} + (-1 - k -
 l) m_{23} + (1 + k - l)(1 + k + l) m_{25}) \\- 
 (1 + k - l) (k + l) \xi_1\xi_3\cdot\xi_1\xi_2\xi_3 \,m_{4} +(-1 - k + l)
 \xi_2\xi_3\cdot\xi_1\xi_2\xi_3 \,m_{7}\end{array}
$$
\normalsize

\noindent 2h) $%[0,1,1,2]
(k+2, k-2, k, -k) \tto (k+2, k-1, k+1, -k+2)$
\vspace*{-2mm}\footnotesize 
$$
\renewcommand{\arraystretch}{1.4}
\begin{array}{l}
 (\xi_1\xi_2)^2 \,m_{15} + 2 k (-1 + 2 k) (\xi_2\xi_3)^2 \,m_{2} - 2
 (-1 + 2 k) \xi_1\cdot\xi_1\xi_2\xi_3 \,m_{3} \\- 2
 \xi_2\cdot\xi_1\xi_2\xi_3 \,(m_{6} + 2 k m_{10}) - 4 k (-1 + 2 k)
 \xi_3\cdot\xi_1\xi_2\xi_3 \,m_{1} + 2 (-1 + 2 k)
 \xi_1\xi_2\cdot\xi_1\xi_3 \,m_{3}\\ - 2 (-1 + 2 k)
 \xi_1\xi_2\cdot\xi_2\xi_3 \,m_{6} + \xi_1\xi_2\cdot\xi_1\xi_2\xi_3 \,
 (m_{34} - 2 (1 + k) m_{36} - 2 k m_{46} - 2 k m_{49})\\- 4 k (-1 + 2
 k) \xi_1\xi_3\cdot\xi_2\xi_3 \,m_{1} + 2 (-1 + 2 k)
 \xi_1\xi_3\cdot\xi_1\xi_2\xi_3 \,( m_{9} - (1 + 2 k) m_{11})\\ + 2
 \xi_2\xi_3\cdot\xi_1\xi_2\xi_3 \, (- (-1 + k) m_{16} + (-1 + 2 k^2 )
 m_{18} + k m_{23} + 2 (-1 + 2 k) m_{25})
\end{array}
$$
\normalsize

\noindent 2i) $(-l+1, k, k, l+1) \tto (-l+2, k+2, k+2, l+2)$
\vspace*{-2mm}\footnotesize 
$$
\renewcommand{\arraystretch}{1.4}
\begin{array}{l}
 (k - l)(\xi_2\xi_3)^2 \,m_{2}
 - 2(k + l) \xi_1\cdot\xi_1\xi_2\xi_3 \,m_{3} 
 + 2\xi_2\cdot\xi_1\xi_2\xi_3 \,(- m_{6} + (1 - k - l) m_{10}) \\- 
 2(-1 + k - l) (k + l) \xi_3\cdot\xi_1\xi_2\xi_3 \,m_{1} 
 + 2\xi_1\xi_2\cdot\xi_2\xi_3 \,(- m_{6} + m_{10}) 
\\ + 
 2\xi_1\xi_2\cdot\xi_1\xi_2\xi_3 \,(m_{46} + m_{49}) 
 +2\xi_1\xi_3\cdot\xi_1\xi_2\xi_3 \,(- m_{9} + (k - l) m_{11}) 
 \\+ 
 2\xi_2\xi_3\cdot\xi_1\xi_2\xi_3 \,(m_{18} - m_{23} + (2 - k + l) m_{25})
\end{array}
$$
\normalsize 

\section{Singular vectors for {\boldmath $\fg=\fk\fas(1|6; 3\eta)$}}

The $m_{i}$ are the following elements of the irreducible $\fg_{0}$-module
$V$: \tiny
$$
\begin{matrix}
\renewcommand{\arraystretch}{1.4}
\begin{array}{l}
m_{2} = \xi_2\,m_{1}\cr
m_{3} = \xi_1\eta_2\,m_{1}\cr
m_{4} = \xi_3\eta_1\,m_{1}\cr
m_{5} = \xi_2\,\xi_1\eta_2\,m_{1}\cr
m_{6} = \xi_2\,\xi_3\eta_1\,m_{1}\cr
m_{7} = (\xi_1\eta_2)^2\,m_{1}\cr
m_{8} = \xi_1\eta_2\,\xi_3\eta_1\,m_{1}\cr
m_{9} = (\xi_3\eta_1)^2\,m_{1}\cr
m_{10} = \xi_1\,m_{1}\cr
\end{array}&
\renewcommand{\arraystretch}{1.4}\begin{array}{l}
m_{11} = \xi_3\eta_2\,m_{1}\cr
m_{13} = \xi_2\,\xi_1\eta_2\,\xi_3\eta_1\,m_{1}\cr
m_{16} = \xi_2\,\xi_3\eta_2\,m_{1}\cr
m_{18} = (\xi_1\eta_2)^2\,\xi_3\eta_1\,m_{1}\cr
m_{19} = \xi_1\eta_2\,(\xi_3\eta_1)^2\,m_{1}\cr
m_{21} = \xi_1\eta_2\,\xi_3\eta_2\,m_{1}\cr
m_{23} = \xi_3\eta_1\,\xi_1\,m_{1}\cr
m_{24} = \xi_3\eta_1\,\xi_3\eta_2\,m_{1}\cr
m_{25} = - (\xi_3\,m_{1})\cr
m_{27} = \xi_2\,(\xi_1\eta_2)^2\,\xi_3\eta_1\,m_{1}\cr
\end{array}&
\renewcommand{\arraystretch}{1.4}\begin{array}{l}
m_{30} = \xi_2\,\xi_1\eta_2\,\xi_3\eta_2\,m_{1}\cr
m_{34} = - (\xi_2\,\xi_3\,m_{1})\cr
m_{37} = (\xi_1\eta_2)^2\,(\xi_3\eta_1)^2\,m_{1}\cr
 m_{41} = \xi_1\eta_2\,\xi_3\eta_1\,\xi_1\,m_{1}\cr
m_{42} = \xi_1\eta_2\,\xi_3\eta_1\,\xi_3\eta_2\,m_{1}\cr
m_{43} = - (\xi_1\eta_2\,\xi_3\,m_{1})\cr
m_{48} = \xi_1\,\xi_3\eta_2\,m_{1}\cr
m_{49} = (\xi_3\eta_2)^2\,m_{1}\cr
m_{86} = - (\xi_1\,\xi_3\,m_{1})\cr
\end{array}\end{matrix}
$$
\normalsize

\noindent {\bf Theorem} \quad In $I(V)$ in degrees $d)$, there are only
the following singular vectors:

\vspace*{1mm}
\noindent 1a) $%[0,0,-1,-1]
(k, l, m, -m) \tto (l, l, m-1,-m-1)$:
$\eta_1\eta_3\,m_{1}$;

\vspace*{1mm}
\noindent 1b) $%[0,-1,0,-1]
(k, l, m, -l-1) \tto (k, l-1, m, -l-2)$:
$\eta_1\eta_3\,m_{3} + (-l + m) \eta_2\eta_3\,m_{1}$;

\vspace*{1mm}
\noindent 1c) $%[0,-1,-1,0]
(k, l, -l-2, m) \tto (k, l-1, -l-3, m)$\vspace*{-2mm}
$$
(1 + l - m) (2 + l + m) \eta_1\eta_2\,m_{1} 
 + \eta_1\eta_3\,(m_{8} + (-2 - l - m) m_{11}) 
 + (-1 - l + m) \eta_2\eta_3\,m_{4}
$$ 

\noindent 1d) $%[-1,0,-1,0]
(k+3, -k-2, k, k-1) \tto (k+2, -k-2, k-1, k-1)$\vspace*{-2mm}
$$
2 k \eta_1\,m_{1} - 2 k \eta_3\,m_{4} 
 + 2 k \eta_1\eta_2\,m_{2}
 + \eta_1\eta_3\,(m_{13} + m_{16}) 
 + 2 k \eta_2\eta_3\,m_{6}
$$

\noindent 1e) $%[-1,-1,0,0]
(k+3, k-1, -k-1, k-1) \tto (k+2, k-2, -k-1, k-1)$
\vspace*{-2mm}\footnotesize 
$$
\renewcommand{\arraystretch}{1.4}
\begin{array}{l}
 2 k \eta_1\,m_{3} - 4 k^2 \eta_2\,m_{1} + 2 k \eta_3\,(- m_{8} + (1 +
 2 k) m_{11}) \\+ 
 2 k \eta_1\eta_2\,( m_{5} + 2 k m_{10}) + \eta_1\eta_3\,(m_{27} - 2 k
 m_{30} + 2 k m_{41} - 2 k m_{43} - 4 k^2 m_{48}) \\+ 
 2 k \eta_2\eta_3\, (- m_{16} - m_{23} + m_{25})
\end{array}
$$
\normalsize

\noindent 1f) $(4, 0, -1, -1) \tto (3, -1, -1, -1)$
\vspace*{-2mm}\footnotesize  
$$
\renewcommand{\arraystretch}{1.4}
\begin{array}{l}
 -2 \eta_1\,m_{3}) + 2 \eta_2\,m_{1} + \eta_3\,(m_{8} - 2 m_{11}) -
 \eta_1\eta_2\,( m_{5} + m_{10})\\+ 
 \eta_1\eta_3\,(m_{27} + m_{43} + m_{48}) - \eta_2\eta_3\,(2 m_{13} -
 m_{16} + m_{23})\end{array}
$$ 
\normalsize

\noindent 1g) $%[-2,0,0,0]
(4, 0, 0, 0) \tto (2, 0, 0, 0)$
\vspace*{-2mm}
$$
 6 \,m_{1} + \eta_1\,(m_{5} + 3 m_{10}) + 3 \eta_2\,m_{2} - \eta_3\,(
 m_{16} + m_{23}) + \eta_1\eta_2\,m_{15} + \eta_1\eta_3\,m_{86} +
 \eta_2\eta_3\,m_{34}
$$ 

\noindent 2a) $%[0,0,-2,-2]
(k, l, m, 2-m) \tto (k, l, m-2, -m)$:
$(\eta_1\eta_3)^2\,m_{1}$;

\vspace*{1mm}
\noindent 2b) $%[0,-1,-1,-2]
(k, l, l+2, -l-1) \tto (0, l-1, l+1, -l-3)$
\vspace*{-2mm}
$$
(\eta_1\eta_3)^2\,m_{3}
 + 2 \eta_1\eta_3\cdot\eta_2\eta_3 \,m_{1}
$$

\noindent 2c) $%[0,-1,-2,-1]
(k, l-2, -l, l+1) \tto (k, l-3, -l-2, l)$, where $l\neq -\frac12$
$$
(\eta_1\eta_3)^2\,(m_{8} - (1 + 2 l) m_{11}) 
 - 2 (1 + 2 l) \eta_1\eta_2\cdot\eta_1\eta_3 \,m_{1} 
 + 2 \eta_1\eta_3\cdot\eta_2\eta_3 \,m_{4}
$$ 

\noindent 2d) $%[0,-2,0,-2]
(k, l, m, -l) \tto (k, l-2, m, -l-2)$
$$
(\eta_1\eta_3)^2\,m_{7}
 + (-1 + l - m) (l - m) (\eta_2\eta_3)^2\,m_{1} 
 - 2 (-1 + l - m) \eta_1\eta_3\cdot\eta_2\eta_3 \,m_{3}
$$

\noindent 2e) $%[0,-2,-1,-1]
(k, l-1, -l-1, -l+1) \tto (k, l-3, -1-2, -l)$
\vspace*{-2mm}\footnotesize 
$$
\renewcommand{\arraystretch}{1.4}
\begin{array}{l}
 (\eta_1\eta_3)^2\,(m_{18} - 2 m_{21})
 + 2 l (-1 + 2 l) (\eta_2\eta_3)^2\,m_{4} 
 + 2 (-1 + 2 l) \eta_1\eta_2\cdot\eta_1\eta_3 \,m_{3}\\- 
 4 l (-1 + 2 l) \eta_1\eta_2\cdot\eta_2\eta_3 \,m_{1} 
 + 2 (-1 + 2 l) \eta_1\eta_3\cdot\eta_2\eta_3 \,(-2 m_{8} + m_{11})
\end{array}
$$
\normalsize
2ea) Particular solution for $l=\frac12$:
$$
(\eta_1\eta_3)^2\,m_{21}
 + (\eta_2\eta_3)^2\,m_{4} 
 + 2 \eta_1\eta_2\cdot\eta_1\eta_3 \,m_{3}
 - 2 \eta_1\eta_2\cdot\eta_2\eta_3 \,m_{1}
 - \eta_1\eta_3\cdot\eta_2\eta_3 \,m_{8}
$$

\noindent 2f) $%[0,-2,-2,0]
(k, l, -l-1, m) \tto (k, l-2, -l-3, m)$, where $m\neq l, l+1, -l-1, -l-2$
\vspace*{-2mm}\footnotesize 
$$
\renewcommand{\arraystretch}{1.4}
\begin{array}{l}
(l - m) (1 + l - m) (1 + l + m) (2 + l + m) (\eta_1\eta_2)^2\,m_{1}\\ +
(\eta_1\eta_3)^2\,(m_{37} - 
2 (2 + l + m) m_{42}+ 
(1 + l + m) (2 + l + m) m_{49})\\
+ (l - m) (1 + l - m)
(\eta_2\eta_3)^2\,m_{9} \\+ 
2 (l - m) (2 + l + m) \eta_1\eta_2\cdot\eta_1\eta_3 \, (m_{8} - 
(1 + l + m) m_{11})\\
- 2 (l - m) (1 + l - m) (2 + l + m)
\eta_1\eta_2\cdot\eta_2\eta_3 \,m_{4} \\+ 2 (l - m)
\eta_1\eta_3\cdot\eta_2\eta_3\, (- m_{19} + (2 + l + m) m_{24})
\end{array}
$$
\normalsize 
2fa) Particular solution for $l=m=0$:\footnotesize 
$$
\renewcommand{\arraystretch}{1.4}
\begin{array}{l}
 4(\eta_1\eta_2)^2\,m_{1} + (\eta_1\eta_3)^2\,(m_{37} - 4 m_{49}) +
 2 (\eta_2\eta_3)^2\,m_{9}\\+
 8 \eta_1\eta_2\cdot\eta_1\eta_3 \,(m_{8} - m_{11}) - 8
 \eta_1\eta_2\cdot\eta_2\eta_3 \,m_{4} + 4 \eta_1\eta_3\cdot\eta_2\eta_3
 \,(- m_{19} + 2 m_{24})
\end{array}
$$ \normalsize
2fb) Particular solution for $l=-\frac12$, $m=\frac12$:\footnotesize 
$$
\renewcommand{\arraystretch}{1.4}
\begin{array}{l}
2 (\eta_1\eta_2)^2\,m_{1} + (\eta_1\eta_3)^2\,m_{42} +
(\eta_2\eta_3)^2\,m_{9} + 4 \eta_1\eta_2\cdot\eta_1\eta_3 \,m_{8}\\ -
4 \eta_1\eta_2\cdot\eta_2\eta_3 \,m_{4} - \eta_1\eta_3\cdot\eta_2\eta_3
\,m_{19}
\end{array}
$$ \normalsize

\noindent 2g) $%[-1,0,-2,-1]
(4, -3, 1, 0) \tto (3, -3, -1, -1)$
\footnotesize 
$$
(\eta_1\eta_3)^2\,(m_{13} + m_{16})
 + 2 \eta_1\cdot\eta_1\eta_3 \,m_{1} 
 - 2 \eta_3\cdot\eta_1\eta_3 \,m_{4}
 + 2 \eta_1\eta_2\cdot\eta_1\eta_3 \,m_{2} 
 + 2 \eta_1\eta_3\cdot\eta_2\eta_3 \,m_{6}
$$
\normalsize

\vspace*{-4mm}
\section{Singular vectors for {\boldmath $\fg=\fv\fas(4|4)$}}
%\vspace*{-2mm}
We consider the following negative operators from $\fg_{0}$:
\vspace*{-2mm}\tiny
$$
\begin{matrix}
\renewcommand{\arraystretch}{1.4}
\begin{array}{l}
a_{5}= {x_2} {{\delta }_3}+{x_3} {{\delta }_2} \cr 
a_{6}= {x_3}{{\delta }_3}\cr 
a_{8} = {x_2} {{\delta }_4}+{x_4} {{\delta}_2} \cr 
a_{9} = {x_3} {{\delta }_4}+{x_4} {{\delta}_3} \cr 
a_{10} = {x_4}{{\delta }_4}\cr 
a_{12} =
-{x_2}{{\partial}_1} +{{\xi }_1} {{\delta }_2}\cr 
a_{13} = -{x_3}{{\partial}_1}+{{\xi }_1} {{\delta }_3} 
\cr\end{array}&\renewcommand{\arraystretch}{1.4}
\begin{array}{l}
a_{14} =-{x_4}{{\partial}_1}+{{\xi }_1}{{\delta }_4}\cr
a_{15} =-2{x_3}{{\delta }_4}-{{\xi }_1}{{\partial}_2}+{{\xi }_2}
 {{\partial}_1}\cr
a_{18} = -{x_3}{{\partial}_2}+{{\xi }_2}{{\delta }_3}\cr
a_{19} = -{x_4}{{\partial}_2}+{{\xi }_2} {{\delta }_4} \cr
a_{20} = 2{x_2}{{\delta }_4}-{{\xi }_1}{{\partial}_3}+{{\xi }_3}
 {{\partial}_1}\cr
a_{25} = -{x_4}{{\partial}_3}+{{\xi }_3}{{\delta }_4}\cr
a_{26} = -2{x_2}{{\delta }_3}-{{\xi }_1}{{\partial}_4}+{{\xi }_4}
 {{\partial}_1}\cr
\end{array}\end{matrix}
$$
%\vspace*{-3mm}
\normalsize
For the basis of Cartan subalgebra we take \tiny
$$
\renewcommand{\arraystretch}{1.4}
\begin{array}{l}
a_{11} =
-\frac{1}{2}{x_1}{{\partial}_1}+\frac{1}{2}{x_2}{{\partial}_2}+
\frac{1}{2}{x_3}{{\partial}_3}+\frac{1}{2}{x_4}{{\partial}_4}+{{\xi
}_1} {{\delta }_1}) \cr 
a_{17}
=\frac{1}{2}{x_1}{{\partial}_1}-\frac{1}{2}{x_2}{{\partial}_2}+\frac{1}{2}{x_3}
{{\partial}_3}+\frac{1}{2}{x_4}{{\partial}_4}+{{\xi}_2} {{\delta }_2}\cr
a_{24} = \frac{1}{2}{x_1}
{{\partial}_1}+\frac{1}{2}{x_2}{{\partial}_2}-
\frac{1}{2}{x_3}{{\partial}_3}\frac{1}{2}{x_4}{{\partial}_4}+{{\xi
}_3} {{\delta }_3}\cr 
a_{32} = \frac{1}{2}{x_1}
{{\partial}_1}+\frac{1}{2}{x_2}{{\partial}_2}+
\frac{1}{2}{x_3}{{\partial}_3}-\frac{1}{2}{x_4}{{\partial}_4}+{{\xi
}_4} {{\delta }_4}) \cr 
\end{array} 
$$
\normalsize 
The $m_{i}$ are the
following elements of the irreducible $\fg_{0}$-module $V$: \tiny
$$
\begin{matrix}
\renewcommand{\arraystretch}{1.4}
\begin{array}{l}
m_{1}\text{ is the highest weight vector}\cr
m_{2} = a_{5}\cdot m_{1}\cr
m_{3} = a_{25}\cdot m_{1}\cr
m_{4} = a_{26}\cdot m_{1}\cr
m_{8} = a_{25}\cdot 
 a_{26}\cdot m_{1}\cr
\end{array}&
\renewcommand{\arraystretch}{1.4}\begin{array}{l}
    m_{10} = a_{8}\cdot m_{1}\cr
m_{11} = a_{20}\cdot m_{1}\cr
m_{24} = a_{26}\cdot  a_{20}\cdot m_{1}\cr
m_{27} = -\, a_{12}\cdot m_{1}\cr
\end{array}\end{matrix}
$$
\normalsize
%\vspace*{-2mm}
\noindent
{\bf Theorem}\quad In $I(V)$, there are only the following singular
vectors:

\vspace*{1mm}
\noindent 1a) $(k, l, l, l) \tto (k+1, l, l, l)$:
${{\delta }_1}\otimes m_{1}$;

\vspace*{1mm}
\noindent 1b) $(-1, 0, 0, 0) \tto \frac12(-1, 1, 1, -1)$:
${{\partial}_4}m_{1} + {{\delta }_1}m_{4}$

%\newpage
%\vspace*{1mm}
\noindent 1c) $\frac12(-1, 1, 1, -1) \tto (0, 1, 0, 0)$:
\qquad
$
{{\partial}_3}m_{1} 
 - {{\partial}_4}m_{3} 
 + {{\delta }_1}m_{11}
$

\vspace*{1mm} \noindent 1d) $(l, k+l, l, l) \tto (l, k+l+1, l, l)$;
two particular cases:

\vspace*{1mm}
\noindent 1da) $l\neq 0 \Longrightarrow k\neq -1$:
\vspace*{-2mm}
$$
- {{\partial}_3}\left( 4\,l \, m_{2} + m_{4} \right) 
 - {{\partial}_4} \left( -\,m_{8} + 4\,l \, m_{10} \right) 
 + {{\delta }_1} \left( m_{24} - 4\,l \, m_{27} \right) 
 - 4\,( 1 + k ) \, l \, {{\delta }_2}m_{1}
$$

\vspace*{-2mm}

\noindent
1db) $l=0\Longrightarrow k\neq 0$:
\qquad
$
{{\partial}_3}m_{2} + {{\partial}_4}m_{10} 
 + {{\delta }_1}m_{27} + k\, {{\delta }_2}m_{1}$

\vspace*{-2mm}

\section*{Acknowledgments}
\addcontentsline{toc}{section}{\numberline{}Acknowledgments}

We are thankful for partial financial support: P.G. to TBSS, I.Shch. 
to RFBR grants 99-01-00245 and 01-01-00490a; D.L.: to NFR. D.L. is also
thankful to Bernstein, Kirillov and Rudakov for teaching and kindness
in times of need.  Larsson's conjectures and generous discussion of
unpublished results\,\cite{La} were helpful to us.\\
PS. After the deadline for submission to the Marinov Memorial Volume
there appeared a preprint\,\cite{KR3} with some singular vectors for
$\fm\fb(3|8)$ and $\fm\fb(5|10)$ (finite dimensional fibers).  We use
the opportunity to add the reference; comparison will be done
elsewhere.

\vspace*{-2mm}

\section*{References}
\addcontentsline{toc}{section}{\numberline{}References}

\end{document}